\documentclass[]{amsart}
\usepackage{amsmath}
\usepackage{amssymb,amsfonts}
\usepackage{xypic}
\xyoption{all}

\newtheorem{theorem}{Theorem}[section]
\newtheorem{definition}[theorem]{Definition}
\newtheorem{proposition}[theorem]{Proposition}

\newtheorem{lemma}[theorem]{Lemma}
\newtheorem{remark}[theorem]{Remark}

\newtheorem{fact}[theorem]{Fact}
\def\m{\ensuremath{\mathcal{M}}}
\def\a{\ensuremath{\mathcal{A}}}
\def\b{\ensuremath{\mathcal{B}}}
\def\c{\ensuremath{\mathcal{C}}}
\def\d{\ensuremath{\mathcal{D}}}
\def\cat{\ensuremath{\mathbf{Cat}}}
\def\mcat{\ensuremath{\mathbf{MonCat}}}
\def\ner{\ensuremath{\mathrm{N}}}
\def\gner{\ensuremath{\Delta}}
\def\set{\ensuremath{\mathbf{Set}}}

\def\diag{\ensuremath{\mathrm{diag}}}
\def\cdco{\ensuremath{Z}^2_{\hspace{-2pt}_{^\mathrm{cat}}}\hspace{-2pt}}
\def\ctco{\ensuremath{Z}^3_{\hspace{-1.5pt}_{^\mathrm{bicat}}}\hspace{-2pt}}
\newcommand{\sset}{\ensuremath{\mathbf{Simpl.Set}}}
\newcommand{\bicat}{\ensuremath{\mathbf{Bicat}}}

\newcommand{\class}{\ensuremath{\mathrm{B}}}
\newcommand{\f}{\ensuremath{\mathcal{F}}}
\newcommand{\g}{\ensuremath{\mathcal{G}}}
\newcommand{\T}{\ensuremath{\mathcal{T}}}
\newcommand{\h}{\ensuremath{\mathcal{H}}}
\newcommand{\ct}{\ensuremath{\mathrm{ct}}}
\newcommand{\lfunc}{\ensuremath{\mathrm{LaxFunc}}}

\newcommand{\hoco}{\ensuremath{\mathrm{hocolim}}}

\begin{document}
\title[{\em Realizations of braided monoidal categories}]
{{Classifying spaces for braided monoidal categories and lax diagrams of bicategories}}

\author{P. Carrasco}
\author{A.M. Cegarra}
\author{A. R. Garz\'{o}n}
\thanks{The authors are much indebted to the referee, whose useful observations greatly improved our exposition. This work has been partially supported by DGI of Spain and FEDER (Project: MTM2007-65431); Consejer\'{i}a de Innovaci\'{o}n
de J. de Andaluc\'{i}a (P06-FQM-1889); MEC de Espa\~{n}a, `Ingenio Mathematica(i-Math)' No. CSD2006-00032
(consolider-Ingenio 2010).}

\address{\newline
Departamento de \'Algebra, Facultad de Ciencias, Universidad de Granada.
\newline 18071 Granada, Spain \newline mcarrasc@ugr.es\newline
acegarra@ugr.es\newline agarzon@ugr.es}

\begin{abstract} This work contributes to clarifying several relationships between certain higher categorical structures and the homotopy type of their classifying spaces. Bicategories (in particular monoidal categories) have well understood simple geometric realizations,  and we here deal with homotopy types represented by lax diagrams of bicategories, that is, lax functors to the tricategory of bicategories.  In this paper, it is proven that, when a certain bicategorical Grothendieck construction is performed on a lax diagram of bicategories, then the classifying space of the resulting bicategory can be thought of as the homotopy colimit of the classifying spaces of the bicategories that arise from the initial input data given by the lax diagram. This result is applied to produce bicategories whose classifying space has a double loop space with the same homotopy type, up to group completion, as the underlying category  of any given (non-necessarily strict) braided monoidal category. Specifically, it is proven that these double delooping spaces, for categories enriched with a braided monoidal structure, can be explicitly realized by means of certain genuine simplicial sets characteristically associated to any braided monoidal categories,  which we refer to as their (Street's) geometric nerves.
\end{abstract}


\maketitle


\section{Introduction and Summary} Higher-dimensional categories provide a suitable setting for the treatment of an extensive list of subjects with recognized mathematical interest. The construction of nerves and classifying spaces of higher categorical structures reveals ways to transport categorical coherence to homotopical coherence  and it has shown its relevance as a tool in algebraic topology, algebraic geometry, algebraic $K$-theory, string field theory,  conformal field theory, and in the study of geometric structures on low-dimensional manifolds. In particular, {\em braided monoidal categories} \cite{joyal} have been playing a key role in recent developments in  quantum theory and its related topics, mainly thanks to the following result, which was the starting point for this paper:

\begin{quote} {\em ``The group completion of the classifying space of a braided monoidal category is a
double loop space"}
\end{quote}

\noindent as  was noticed by J. D. Stasheff in \cite{sta}, but originally proven by Z. Fiedorowicz in \cite[Theorem 2]{fie} (some other proofs can be found in \cite[Theorem 1.2]{berger} or in \cite[Theorem  2.2]{b-f-s-v}, for example).
More precisely, given any braided monoidal category
$$(\m,\otimes,\boldsymbol{c})=(\m,\otimes,\text{I},\boldsymbol{a},\boldsymbol{l},\boldsymbol{r},\boldsymbol{c}),$$
Stasheff-Fiedorowicz's theorem implies the existence of a path-connected, simply connected space, uniquely defined up to homotopy equivalence, $\class(M,\otimes,\boldsymbol{c})$, and a homotopy-natural map $\class\m\to\Omega^2 \class(M,\otimes,\boldsymbol{c}),$
where $\class\m$ is the classifying space of the underlying category $\m$, which is, up to group completion, a homotopy equivalence. Hereafter, we shall refer to $\class(M,\otimes,\boldsymbol{c})$ both as the {\em classifying space of the braided monoidal category} and as the {\em double delooping} of $\class\m$, induced by the  braided monoidal structure given on $\m$.

However, there is a problem with the space $\class(M,\otimes,\boldsymbol{c})$  since its existence is proven as an application of May's theory of $E_2$-operads \cite{may74} and, therefore, its various known constructions are based on some complicated and irritating processes of rectifying homotopy coherent diagrams. In fact, the double delooping construction is provided by  May's bar-construction that only takes place after replacing $(\m,\otimes,\boldsymbol{c})$ by an equivalent {\em strict} braided monoidal category $(\m',\otimes',\boldsymbol{c}')$, and then by carrying out a substitution of  $\class\m'$ by a homotopy equivalent space upon which the little square operad of Boardman-Vogt acts \cite{b-v}, which depends on an explicit equivalence of operads between the  braided operad used and the little 2-cube one. The resulting CW-complex thus obtained has many cells with little apparent intuitive connection with the data of the original monoidal category, and this leads one to search for any simplicial set, say ``nerve of the braided monoidal category", realizing the space $\class(M,\otimes,\boldsymbol{c})$ and whose cells give a logical geometric meaning to the data of the braided monoidal category.

A natural response for that nerve was postulated in the nineties  by J. Dolan and R. Street (probably among others) and it is as follows: since a braided monoidal category can be regarded as a one-object,
one-arrow  tricategory  \cite[Corollary 8.7]{g-p-s} and each category as a tricategory
 whose 2-cells and 3-cells are all identities, one can consider strictly unitary lax
functors from the categories $[p]=\{0<1<\cdots<p\}$ to the tricategory  $\Omega^{^{-2}}\hspace{-3pt}\m$ that the braided monoidal
category $(\m,\otimes,\boldsymbol{c})$ defines. Then, its {\em geometric nerve} is the simplicial set
$$ Z^3(\m,\otimes,\boldsymbol{c}): [p]\mapsto \mathrm{NorLaxFunc}([p],\Omega^{^{-2}}\hspace{-5pt}\m),$$
whose $p$-simplices are all strictly unitary lax functors $[p]\to \Omega^{^{-2}}\hspace{-5pt}\m$ (also called 3{\em -cocycles with coefficients in the braided monoidal category}, \cite{cegarra3}). This geometric nerve of the braided monoidal category is a $4$-coskeletal simplicial set whose simplices have a pleasing interpretation: there is only one 0-simplex, there is only one 1-simplex, the  $2$-simplexes $x$ are the objects of $\m$, the $3$-simplexes $\zeta$ with 2-faces (in order) $x_0,x_1,x_2,x_3$ are morphisms $\zeta:x_0\otimes x_2\to x_3\otimes x_1$, and so on.
 The most striking instance is for  $(\m,\otimes,\boldsymbol{c})=(A,+,0)$, the strict braided monoidal category with only one object defined by an abelian group $A$, where both composition and tensor product are given by the addition $+$ in $A$; in this case, $Z^3(A,+,0)=K(A,3)$, the minimal Eilenberg-Mac Lane complex. Geometric nerves of {\em braided categorical groups} \cite[\S 3]{joyal} were studied in \cite{Ca-Ceg}, where it was proven that the mapping  $(\m,\otimes,\boldsymbol{c})\mapsto |Z^3(\m,\otimes,\boldsymbol{c})|$ induces an equivalence between the homotopy category of braided categorical groups and the homotopy category of pointed 1-connected 3-types (a fact due to A. Joyal and M. Tierney \cite{j-t}, see also \cite[Theorem 3.3]{berger}).

A main goal of this article is to prove the following result for which, as far as we know, no proof has yet
appeared in the literature:

\begin{quote} {\em ``For any braided monoidal category $(\m,\otimes,\boldsymbol{c})$, there is a  homotopy equivalence
$\class(\m,\otimes,\boldsymbol{c})\simeq |Z^3(\m,\otimes,\boldsymbol{c})|$."}\end{quote}

Our proof for this theorem requires a long preliminary  discussion on the notion of {\em realization} or {\em classifying space} for lax diagrams of bicategories $I^{\mathrm{op}}\to \bicat$, where $I$ is any small category. This requirement is due to the fact that, in a first approach, we show that the space $\class(\m,\otimes,\boldsymbol{c})$ can be realized by means of the pseudo-simplicial bicategory
$$\ner(\m,\otimes,\boldsymbol{c}):\Delta^{\!^{\mathrm{op}}} \to\bicat,\hspace{0.6cm}[p]\mapsto \Omega^{^{-1}}\hspace{-5pt}\m^p,$$
defined (thanks to the braiding) by the familiar bar construction; here $\Omega^{^{-1}}\hspace{-5pt}\m$ denotes the one-object bicategory delooping of the underlying monoidal category $(\m,\otimes)$, that is, that obtained forgetting the braiding. Then, the proof we give of the claimed above result reduces to show the existence of a homotopy equivalence between the realization of the simplicial set geometric nerve $Z^3(\m,\otimes,\boldsymbol{c})$ (viewed as a simplicial discrete bicategory) and the  realization of the pseudo-simplicial bicategory $\ner(\m,\otimes,\boldsymbol{c})$.

Hence, much of our work here is dedicated to establishing and proving the most basic results needed concerning the homotopy theory of lax diagrams of bicategories, paralleling   corresponding facts for lax diagrams of categories as  stated and proven by G. Segal \cite{segal74} and R. W. Thomason \cite{thomason}, following the methods of A. Grothendieck. The resulting theory is in itself of independent interest and yields, as an added benefit, the foundation for other future developments,
for example in the homotopy theory of monoidal bicategories or arbitrary tricategories. Although this subject will not be treated here, let us say that the classifying space of any monoidal bicategory  $(\b,\otimes)$ is precisely the realization, in the sense studied here, of the  pseudo-simplicial bicategory
$\ner(\b,\otimes): \Delta^{\!^{\mathrm{op}}}\to \bicat$, $[p]\mapsto \b^p$, which it defines by the reduced bar construction.

After this introductory Section 1, the paper is organized in six sections. Section 2 is an attempt to make the paper as self-contained as possible; hence, at the same time as we fix notations and terminology, we review in it some necessary aspects from the background of bicategories by briefly describing $\mathbf{Bicat}$, the tricategory of bicategories,  homomorphisms, pseudo natural transformations, and modifications. This material is quite standard, so the expert reader may skip most of it, but note that some notations may be idiosyncratic.
Also, we
describe the kind of lax diagrams of bicategories we are going to treat in this paper: lax morphisms of tricategories in the sense of \cite{g-p-s}, $\f: I^\mathrm{op}\to \bicat$, where $I$ is any small category, all of whose coherence 3-cells are invertible. For any given category $I$, the lax diagrams of bicategories are the
objects of a  tricategory, denoted $\bicat^{ I^{\mathrm{op}}}$.
The following two sections, 3 and 4, are very technical, but crucial to our discussions. Section 3 is mainly dedicated to study a bicategorical Grothendieck construction  \cite{grothendieck,thomason}. More precisely, the aim there is to prove the following:
\begin{quote}{\em
`` \!There is a Grothendieck construction on lax diagrams of bicategories defining  a trihomomorphism of tricategories $$\xymatrix{\int_I :\bicat^{ I^{\mathrm{op}}}\rightarrow \bicat}$$ which, moreover, is left triadjoint to the diagonal trihomomorphism $\bicat\to \bicat^{ I^{\mathrm{op}}}$."}
\end{quote}

Hence, the function on objects of the Grothendieck construction assembles any lax diagram $\f:I^\mathrm{op}\to \bicat$ into a large bicategory
$\xymatrix{\int_I\f}$,
which is a lax colimit of the bicategories $\f_i$, $i\in \mbox{Ob}I$,
 and, as we shall detail later, it can be thought as its  homotopy colimit.  Section 4 is dedicated to proving, following Giraud and Street's methods \cite{giraud,street72}, that
\begin{quote}{\em
`` \!There exists a  rectifying  trihomomorphism ${(\
)^{^\mathrm{r}} :\bicat^{ I^{\mathrm{op}}}\to \bicat^{ I^{\mathrm{op}}}}$\!"
}\end{quote}
through which any lax diagram of bicategories $\f: I^\mathrm{op}\to \bicat$ has naturally associated a genuine $I$-diagram of
bicategories, that is, a functor $\f^{^\mathrm{r}}:I^\mathrm{op}\to \bicat $ that, as we will show later, represents
the same homotopy type as the original $\f$.

 Heavily dependent on the results in \cite{ccg}, where nerves and classifying spaces of bicategories are studied, in Section 5 we introduce and study realizations for lax diagrams of bicategories. The classifying space of the lax diagram of bicategories $\f: I^\mathrm{op}\to \bicat$, denoted $\class \f$, is defined to be $\class\!\! \int_I\!\f$, the classifying space of the bicategory obtained by the Grothendieck construction on $\f$, and the more basic and relevant properties of this construction $\f\mapsto\class\f$ are stated and proven throughout the section. Namely, we prove the following two results:

\begin{quote}{\em ``If $F:\f\to\g$ is a lax $I$-homomorphism between lax $I$-diagrams $\f,\g:  I^{{\mathrm{op}}}\to\bicat$, such that the induced maps $\class F_i:\class \f_i\to \class\g_i$ are homotopy equivalences, for all objects $i$ of $I$, then the induced map $\class F:\class\f \to \class\g$ is a homotopy equivalence."}
\end{quote}
\vspace{0.2cm}
\begin{quote}{\em ``Let $\f:I^{^{\mathrm{op}}}\to \bicat$ be a lax diagram of bicategories such that the induced map
$\class a^*:\class \f_{i}\to \class \f_{j}$, for each morphism $a:j\to i$ in  $I$,  is a homotopy equivalence. Then, for every object $i$ of $I$, there is a homotopy fibre sequence  $
\class\f_i \hookrightarrow\class \f \to \class I
$."
}
\end{quote}

In Section 6, the facts demonstrated on realizations for lax diagrams of bicategories are mainly applied to state and prove several facts concerning the classifying space construction $(\m,\otimes,\boldsymbol{c})\mapsto \class (\m,\otimes,\boldsymbol{c})$, for braided monoidal categories. Specifically, we give here a new proof of the above-mentioned Stasheff-Fiedorowicz theorem that, as an added value, includes the following more explicit fact:

\begin{quote}
``{\em For any braided monoidal category, the double loop space of the realization of its geometric nerve is a group completion of the classifying space of the underlying category." }
\end{quote}

And finally,  Section 7 simply collects the expression of various coherence conditions concerning  definitions in Section 2 and used throughout the paper.

\vspace{0.3cm}

\subsection*{Some Frequently Used Notations}~

\vspace{0.2cm}
To help the reader we list below the following notations used along the paper, with indication of their meaning and first  appearance.

\vspace{0.2cm}
\noindent\begin{tabular}{llc}
  $\bicat$ &\small{tricategory of bicategories}& (\ref{bicat})\\
  $\mathbf{Hom}$ &\small{category of bicategories and homomorphisms}& (\ref{hom}) \\
  $\bicat^{ I^{\mathrm{op}}}$ &\small{tricategory of lax $I$-diagrams of bicategories}& (\ref{bicatiop})  \\
  $\xymatrix{\int_I\f}$ &\small{Grothendieck construction on a lax $I$-diagram of bicategories}& (\ref{gc})  \\
  $\f^{^\mathrm{r}}$ &\small{rectification construction on a lax $I$-diagram of bicategories}& (\ref{recti}) \\
  $\ner\c$ &\small{pseudo-simplicial nerve of a bicategory}& (\ref{ps1.1}) \\
  $\class \c$ &\small{classifying space of a bicategory}& (\ref{clasi})  \\
 $\Delta^{\hspace{-2pt}^\mathrm{u}}\!\c$&\small{unitary geometric nerve of a bicategory}&(\ref{ngn})\\
$\Delta\c$&\small{geometric nerve of a bicategory}&(\ref{geb})\\
  $\class \f $ &\small{classifying space of a lax $I$-diagram of bicategories}&  \ref{defclass} \\
  $\Omega^{^{-1}}\hspace{-5pt}\m$ &\small{delooping bicategory of a monoidal category}& (\ref{ome})\\
 $\ner(\m,\otimes)$ &\small{pseudo-simplicial nerve of a monoidal category}&  (\ref{psnm})\\
 $\class(\m,\otimes)$& \small{classifying space of a monoidal category}& (\ref{clmc}) \\
 $\Omega^{^{-2}}\hspace{-5pt}\m$ & \small{double delooping tricategory of a braided monoidal category} & (\ref{omenosdos}) \\
 $\ner(\m,\otimes,\boldsymbol{c})$ &\small{pseudo-simplicial nerve of a braided monoidal category}& (\ref{psnbmc})\\
 $\class(\m,\otimes,\boldsymbol{c})$ &\small{classifying space of a braided monoidal category}& (\ref{cbmc})\\[3pt]
 $Z^2(\m,\otimes)$ &\small{geometric nerve of a monoidal category}& (\ref{gnmc})\\[3pt]
 $\cdco(\m,\otimes)$ &\small{categorical geometric nerve of a monoidal category}& (\ref{cgnm})\\[3pt]
 $Z^3(\m,\otimes,\boldsymbol{c})$ &\small{geometric nerve of a braided monoidal category}& (\ref{gnbmc})\\[3pt]
$\ctco(\m,\otimes,\boldsymbol{c})$ &\small{bicategorical geometric nerve of a braided monoidal category}& (\ref{cgnbmc})\\
\end{tabular}

\section{Bicategorical preliminaries: Lax diagrams of bicategories.}
 We shall begin by reviewing some necessary facts concerning the  tricategory of bicategories. Also, we will
describe the kind of lax diagrams of bicategories we are going to treat in this paper.

\subsection{The tricategory of bicategories}
We refer to \cite{benabou,g-p-s,gurski} and \cite{street} for  background on bicategories and tricategories. For definiteness or emphasis, we state the following:

 In any small {\em bicategory} $\a$,  its set of objects (or 0-cells) is denoted by $\mbox{Ob}\a$ and, for each ordered pair of objects $(y,x)$,  $\a(y,x)$ is the category whose objects $u:y\to x$ are the 1-cells (or morphisms) of $\a$ with source $y$ and target $x$, and whose arrows ${\alpha:u\Rightarrow u'}$  are the 2-cells (or deformations) of $\a$. The composition of deformations in each category $\a(y,x)$, that is, the vertical composition of 2-cells,  is denoted by $\beta\cdot \alpha$, while the symbol $\circ$ is used to denote the horizontal composition functors $\circ:\a(y,x)\times\a(z,y)\to\a(z,x)$.
The identity of an object is written as $1_x:x\to x$, and we shall use the letters $\boldsymbol{a}$, $\boldsymbol{r}$, and $\boldsymbol{l}$ to denote the associativity, right unit, and left unit constraints of the bicategory, respectively.

A {\em lax functor} is usually written as a pair $ F=(F,\widehat{F}):\a \to \b$ since we will generically denote its structure constraints by $\widehat{F}_{u,v}:Fu\circ Fv\Rightarrow F(u\circ v)$ and $\widehat{F}_x:1_{Fx}\Rightarrow F1_x$, or merely by  $\widehat{F}:Fu\circ Fv\Rightarrow F(u\circ v)$ and $\widehat{F}:1_{Fx}\Rightarrow F1_x$ since the source and target of this constraint make it clear what kind of constraint deformation it is.
The lax functor is termed a {\em pseudo-functor} or {\em homomorphism}
whenever all the structure constraints $\widehat{F}$ are invertible.  If the unit constraints $\widehat{F}_x$ are all identities, then the lax functor is qualified as (strictly)  {\em unitary} or {\em normal} and if, moreover, the constraints $\widehat{F}_{u,v}$ are also identities, then $F$ is called a $2$-{\em functor}.

We will use  pasting diagrams of 2-cells inside bicategories. A diagram of the form
\begin{equation}\label{e1.1}\xymatrix@C=10pt@R=3pt{
&x_1\ar[r]&\cdots\ar[r]&x_n\ar[dr]^{\textstyle u_n}&\\x\ar[ur]^{\textstyle u_0}\ar[rd]_{\textstyle v_0}&&\Downarrow\varphi&&y\\
&y_1\ar[r]&\cdots\ar[r]&y_m\ar[ur]_{\textstyle v_m}&}\end{equation} will represent a deformation $\varphi$
whose source (resp. target) is obtained by horizontal composition of the morphisms in the string
$u_0,\dots,u_n$ (resp. $v_0,\dots,v_m)$ following a particular given association.  By the
bicategorical coherence theorem, such a deformation  uniquely determines another when any other
particular bracketing is used for computing the source and the target morphism from the given
strings of morphisms. Therefore, diagram (\ref{e1.1}) is not ambiguous once a choice of
association has been made for the source and target of the deformation.
When $F:\a\to \b$ is a homomorphism and diagram (\ref{e1.1}) is given in $\a$, then we will
denote by
\begin{equation}\label{e1.2}\xymatrix@C=10pt@R=2pt{
&Fx_1\ar[r]&\cdots\ar[r]&Fx_n\ar[dr]^{\textstyle Fu_n}&\\Fx\ar[ur]^{\textstyle Fu_0}\ar[rd]_{\textstyle Fv_0}&&\Downarrow F\varphi&&Fy\\
&Fy_1\ar[r]&\cdots\ar[r]&Fy_m\ar[ur]_{\textstyle Fv_m}&}\end{equation} the diagram in $\b$ in which the
deformation is obtained by appropriately composing the original $F\varphi$ with constraints $\widehat{F}$
of $F$.
That diagram (\ref{e1.2})
is well defined from diagram (\ref{e1.1}) is a consequence of the coherence theorem for
homomorphisms of bicategories \cite[Theorem 1.6]{g-p-s}.
A diagram such as (\ref{e1.2}), with the symbol $\cong$ inside instead of $\Downarrow\!F\varphi$, means that the deformation is obtained only by composition of the structure constraints of the homomorphism $F$ and the bicategories  involved.

If $F,F':\a \to \b$ are lax functors, then we follow the convention of \cite{g-p-s} in what is meant by a {\em lax transformation}
${\alpha=(\alpha,\widehat{\alpha}):F\Rightarrow F'}$. Thus, $\alpha$ consists of morphisms ${\alpha x:Fx\to F'x}$, $x\in
\mbox{Ob}\a$, and of deformations $\widehat{\alpha}_u:\alpha y\circ Fu\Rightarrow F'\!u\circ \alpha x$
 that are natural on morphisms $u:x\to y$, subject to the usual two axioms.  When the deformations
$\widehat{\alpha}_u$ are all invertible, we say that $\alpha$ is a {\em pseudo transformation}.
In accordance with the orientation of the naturality deformations chosen, if
${\alpha,\alpha':F\Rightarrow F'}$ are two lax transformations, then a {\em modification}
$\varphi:\alpha\Rrightarrow\alpha'$ will consist of deformations
$\varphi x:\alpha x\Rightarrow\alpha' x$, $x\in \mbox{Ob}\a$, subject to the commutativity
condition $(1_{F'u}\circ \varphi x)\cdot \widehat{\alpha}_u=\widehat{\alpha}'_u\cdot(\varphi y\circ 1_{Fu})$, for each
morphism $u:x\to y$ of $\a$.

Next, we shall briefly  describe the most striking example of {\em tricategory}: the tricategory of bicategories,  homomorphisms, pseudo-natural transformations and modifications, denoted by \begin{equation}\label{bicat}\mathbf{Bicat}.\end{equation} We refer the reader to
\cite[\S 5]{g-p-s} and \cite[\S 6.3]{gurski} for more details.

For bicategories $\a$, $\b$, $\bicat(\a,\b)$ denotes the bicategory whose objects are the
homomorphisms $F:\a\to \b$, 1-cells the pseudo-transformations $\alpha:F\Rightarrow F'$, and
2-cells the modifications $\varphi:\alpha \Rrightarrow \alpha'$. Let us briefly recall  that
a modification $\varphi:\alpha\Rrightarrow\alpha'$ composes vertically with a modification
$\varphi':\alpha'\Rrightarrow \alpha''$ yielding the modification
$\varphi'\cdot\varphi:\alpha\Rrightarrow \alpha''$, such that
$(\varphi'\cdot\varphi)x=\varphi'x\cdot\varphi x$, $x\in \mbox{Ob}\a$. The horizontal composition of 1-cells in $\bicat(\a,\b)$ is given by the ``vertical composition" of pseudo-transformations:  for $\alpha:F\Rightarrow F'$ and
$\alpha':F'\Rightarrow F''$, where ${F,F',F'':\a\to \b}$, the composite
$\alpha'\circ\alpha:F\Rightarrow F''$  is defined by putting
$(\alpha'\circ\alpha)x=\alpha'x\circ\alpha x$ for any object $x$ of $\a$,  the component of
$\alpha'\circ\alpha$ at a morphism $u:x\to y$  being the deformation obtained by pasting the diagram
$$
\xymatrix{Fx\ar[r]^{\textstyle \alpha x}\ar[d]_{\textstyle Fu}\ar@{}[dr]|{\textstyle \overset{\textstyle \widehat{\alpha}_u}\Rightarrow}&F'\!x\ar[r]^{\textstyle \alpha'\!x}\ar[d]|{F'\!u}\ar@{}[dr]|{\textstyle \overset{\textstyle \widehat{\alpha}'_u}\Rightarrow}&F''\!x\ar[d]^{\textstyle F''\!u}\\Fy\ar[r]_{\textstyle \alpha y}&F'\!y\ar[r]_{\textstyle \alpha' \!y}&F''\!y\,.}
$$
The horizontal composition of a modification
$\psi:\alpha\Rrightarrow \beta:F\Rightarrow F'$ with a modification
$\psi':\alpha'\Rrightarrow\beta':F'\Rightarrow F''$ is the modification
$\psi'\circ\psi:\alpha'\circ\alpha\Rrightarrow\beta'\circ\beta$ such that
$(\psi'\circ\psi)x=\psi'x\circ\psi x$, for each object $x$ of $\a$. The structure constraints in
$\bicat(\a,\b)$ are canonically derived, in a pointwise way, from those of $\b$; thus, for
example, the associativity modifications $\boldsymbol{a}:\alpha''\circ (\alpha'\circ
\alpha)\Rrightarrow (\alpha''\circ\alpha')\circ\alpha$ are given by
$\boldsymbol{a}x=\boldsymbol{a}_{\alpha''\!x,\alpha'\!x,\alpha x}$, $x\in \mbox{Ob}\a$.

The composition of lax functors $F:\a \to\b$ and $G: \b\to \c$ will be denoted by juxtaposition,
that is, $GF:\a \to \c$. And recall that its constraints are obtained from those of
$F$ and $G$ by the rules $\widehat{GF}_{u,v}=G\widehat{F}_{u,v}\cdot \widehat{G}_{Fu,Fv}$ and  $\widehat{GF}_{x}=G\widehat{F}_{x}\cdot \widehat{G}_{Fx}$.
This composition of lax functors is associative and unitary, so that the category of bicategories and lax functors is defined. Following \cite[Notation 4.9]{g-p-s},  the category of bicategories with homomorphisms between them will be denoted \begin{equation}\label{hom}\mathbf{ Hom}.\end{equation}

The composition of homomorphisms gives the function
on objects of a homomorphism of bicategories
\begin{equation}\label{e1.3}\bicat(\b,\c)\times \bicat(\a,\b)\to
\bicat(\a,\c),\end{equation} which on
$\xymatrix @C=2pt{\a  \ar@/^0.5pc/[rr]^{F} \ar@/_0.5pc/[rr]_{ F'} & {\Downarrow\!\alpha} &\b\ar@/^0.5pc/[rr]^{ G} \ar@/_0.5pc/[rr]_{G'} & {\Downarrow\!\beta} &\c }$,
is given by $\beta \alpha=\beta F'\circ G\alpha$,  where the
 pseudo-transformations $G\alpha:GF\Rightarrow GF'$ and $\beta F':GF'\Rightarrow G'F'$ are those whose respective components at an object $x$ of $\a$ are the morphisms $G\alpha x$ and $\beta F'x$, and
at a morphism $u$ are $\widehat{G\alpha}_u=G\widehat{\alpha}_u$ and $\widehat{\beta F'}\!_u=\widehat{\beta}_{F'\!u}$.
Similarly, the composition ${\psi\varphi:\beta\alpha\Rrightarrow \beta'\alpha'}$, of modifications
$\varphi:\alpha\Rrightarrow \alpha'$ and $\psi:\beta\Rrightarrow \beta'$, is  given by the formula
$\psi\varphi=\psi F'\circ G\varphi$, that is,  the modification whose component at an object $x\in
\a$ is $(\psi\varphi) x=\psi{F'\!x}\circ G\varphi x$. Moreover, given homomorphisms $\a\overset{F}\to \b\overset{G}\to\c$ and pseudo transformations
$F\overset{\alpha}\Rightarrow F'\overset{\alpha'}\Rightarrow F'':\a\to \b$ and
$G\overset{\beta}\Rightarrow G'\overset{\beta'}\Rightarrow G'':\b\to \c$, the structure constraints
of the homomorphism (\ref{e1.3}) at them are provided by the invertible modifications
\begin{equation}\label{e1.5} 1_{GF}\Rrightarrow 1_G1_F,\hspace{0.3cm} \beta'\alpha'\circ\beta\alpha\Rrightarrow
(\beta'\circ\beta)(\alpha'\circ\alpha)\,,\end{equation} whose respective components at an object $x\in \mbox{Ob}\a$ are
given by pasting the diagrams
$$
\xymatrix@C=20pt{&GFx\ar[rd]^{\textstyle 1_{GFx}}&\\GFx\ar@/^1pc/[ru]^{\textstyle G1_{Fx}}_(0.45){\ \textstyle \cong}\ar@/_1pc/[ru]| {\ \ 1_{GFx}}\ar[rr]_{\textstyle 1_{GFx}}&\ar@{}[u]_(0.4){\textstyle \cong}&GFx,}
\xymatrix{\ar@{}[drr]|(0.67){\textstyle\cong}&&GF'\!x\ar@{}[dr]|{\textstyle\overset{(\ref{4})}\cong}\ar[r]^{\textstyle \beta F'\!x}\ar[d]|{ G\alpha'\!x}&G'F'\!x\ar[d]|{ G'\alpha'\!x}\ar[rd]^{\textstyle \beta'\!F''\!x\circ G'\!\alpha'\!x}& \\
\ar[rru]^{ \textstyle G\alpha x}GFx\ar[rr]_{\textstyle G(\alpha'\!x\circ \alpha x)}&&GF''\!x\ar[r]_{\textstyle\beta F''\!x}&G'\!F''\!x\ar[r]_{\textstyle \beta'F''\! x}\ar@{}[ru]|(0.3){\textstyle =}&G''\!F''\!x,
}
$$
where, for any horizontally composable pseudo transformations $\alpha$ and $\beta$ as above, the invertible modification
\begin{equation}\label{4}\beta F'\circ G\alpha \Rrightarrow
G'\alpha\circ \beta F\,,\end{equation}  at an object $x$ of $\a$, is $\widehat{\beta}_{\alpha x}$, the
component of $\beta$ at the morphism $\alpha x$.

The composition of homomorphisms is associative and unitary as we have remarked before. Besides, the unit constraints for the compositions (\ref{e1.3})  are the pseudo-natural equivalences $\boldsymbol{l}$ and $\boldsymbol{r}$, whose components at any homomorphism ${F:\a\to \b}$ are both the identity transformations on it, and at a pseudo transformation $\alpha:F\Rightarrow F'$ are the modifications \begin{equation}\label{unibicat} \widehat{\boldsymbol{l}}:1_{F'}\circ 1_{1_\b}\alpha\Rrightarrow\alpha\circ 1_F,\hspace{0.8cm} \widehat{\boldsymbol{r}}:1_{F'}\circ \alpha
1_{1_\a}\Rrightarrow\alpha\circ 1_F,\end{equation} canonically obtained from the modifications
$1_{1_\b}\alpha\Rrightarrow\alpha$ and $ \alpha
1_{1_\a}\Rrightarrow\alpha$, respectively defined by the 2-cells of $\b$, $x\in \mbox{Ob}\a$,
$$1_{F'x}\circ \alpha x
\overset{\textstyle \boldsymbol{l}_{\alpha x}}\Rightarrow\alpha x,\hspace{0.6cm}
\alpha x \circ
F1_x \hspace{-5pt}\overset{\textstyle 1\!\circ\! \widehat{F}^{-1}}\Longrightarrow \hspace{-5pt}\alpha x\circ 1_{Fx}\overset{\textstyle  \boldsymbol{r}}\Rightarrow\alpha x.$$ Also, for any homomorphisms $\a\overset{F}\to \b \overset{G}\to\c\overset{H}\to \d$, the associativity pseudo-natural equivalence $\boldsymbol{a}:H(GF)\Rightarrow (HG)F$ is the identity on the composite homomorphism $HGF$, and its component at a morphism $(\gamma,\beta,\alpha):(H,G,F)\Rightarrow
(H',G',F')$ is the modification \begin{equation}\label{asso}\widehat{\boldsymbol{a}}:1_{_{H'\!G'\!F'}}\circ \gamma(\beta\alpha)\Rrightarrow
(\gamma\beta)\alpha\circ 1_{_{HGF}},\end{equation} canonically obtained from the invertible modification  $\gamma(\beta\alpha)\Rrightarrow
(\gamma\beta)\alpha$ associating to each object $x$ of $\a$ the 2-cell of $\d$ given by the composition
$$
\gamma{G'\!F'\!x}\circ H(\beta {F'\!x}\circ G\alpha x)\overset{ 1\circ\widehat{H}^{-1}}\Rightarrow
\gamma{G'\!F'\!x}\circ (H\beta{F'\!x}\circ HG\alpha x)\overset{ \boldsymbol{a}}\Rightarrow
(\gamma{G'\!F'\!x}\circ H\beta{F'\!x})\circ HG\alpha x
\,.
$$

In $\bicat$, the structure invertible modifications $\pi$ and $\mu$, as in the definition of a tricategory \cite{g-p-s}, at any homomorphisms $\a\overset{F}\to \b \overset{G}\to\c\overset{H}\to \d\overset{K}\to \mathcal{E}$,

\begin{equation}\label{pimu}
\xymatrix{ K(H(GF))\ar@{}[rrd]|{\textstyle \overset{\textstyle \pi_{_{K,H,G,F}}}\Rrightarrow}\ar@{=>}[d]_{\textstyle 1_K\boldsymbol{a}}\ar@{=>}[rr]^-(0.6){\textstyle \boldsymbol{a}}&&(KH)(GF)\ar@{=>}[d]^{\textstyle \boldsymbol{a}}\\K((HG)F)\ar@{=>}[r]^{\textstyle \boldsymbol{a}}&
(K(HG))F\ar@{=>}[r]^{\textstyle \boldsymbol{a}1_F}&((KH)G)F\\
}\hspace{0.2cm}
\xymatrix@C=10pt{G(1_\b F)\ar@{=>}[rr]^{\textstyle \boldsymbol{a}}\ar@{=>}[rd]_{\textstyle 1_G\boldsymbol{l}}&
\ar@{}[d]|(0.4){\textstyle \overset{\textstyle \mu_{_{G,F}}}\Lleftarrow}&(G1_\b)F\ar@{=>}[ld]^{\textstyle \boldsymbol{r}1_F}\\ &GF& }
\end{equation}
are respectively given by the unique coherence 2-cells, $x\in \mbox{Ob}\a$,
$$
\xymatrix@C=25pt{KHGFx \ar[d]^{\textstyle K1_{\!_{HGFx}}}\ar[r]^{\textstyle 1}&KHGFx\ar@{}[rd]_(0.3){\textstyle \cong}\ar[r]^{\textstyle 1}&KHGFx&\\KHGFx\ar[r]^{\textstyle 1}&KHGFx\ar[r]^{\textstyle 1}&KHGFx\ar[r]_{\textstyle KHG1_{\!_{Fx}}}&KHGFx,\ar[lu]_{\textstyle 1}}
\xymatrix@C=8pt{GFx\ar@{}[rd]^(0.6){\textstyle \cong}\ar[r]^{\textstyle 1}\ar[d]_{\textstyle 1}&GFx\ar[r]^{\textstyle 1}&GFx \ar[d]^{\textstyle G1_{\!_{Fx}}}\\GFx \ar[rr]^(0.55){\textstyle G1_{\!_{Fx}}}&&GFx.}
$$
The structure modifications $\lambda$ and $\rho$  can be defined in a similar fashion.

\vspace{0.2cm}
\subsection{The tricategory of lax diagrams of bicategories}

Throughout the paper, a  {\em lax diagram of bicategories}, with the shape of a small category $I$, means a lax functor of tricategories \cite[Definition 3.1]{g-p-s}
$$\f=(\f,\chi,\iota, \omega,\gamma,\delta): I^\mathrm{op}\to \bicat,$$ from
$I^\mathrm{op}$, regarded as a tricategory in
which the 2-cells and 3-cells are all identities, to the tricategory
$\mathbf{Bicat}$ of small bicategories, {\em all of whose coherence $3$-cells are invertible and such that each homomorphism $I(j,i)\to \bicat(\f_i,\f_j)$ is normal} (cf. \cite{g-g}, where they are called {\em lax homomorphisms}). The homomorphism $\f a$ attached at an arrow $a:j\to i$ of $I$ is usually written  as $$a^*:\f_i \to \f_j,$$ so that the remaining data of the lax diagram $\f$ provide us with pseudo transformations $$\xymatrix @C=7pt{\f_i  \ar@/^0.9pc/[rr]^{\textstyle b^*a^*} \ar@/_0.9pc/[rr]_{\textstyle (ab)^*} & {}_{\textstyle \Downarrow\!\chi\!=\!\chi_{a,b}} &\f_k},\hspace{0.4cm}\xymatrix @C=7pt{\f_i  \ar@/^0.9pc/[rr]^{\textstyle 1_{\f_i}} \ar@/_0.9pc/[rr]_{\textstyle 1_i^*} & {}_{\textstyle \Downarrow\!\iota\!=\!\iota_{i}} &\f_i},$$ respectively associated to pairs of composible arrows  $k\overset{b}\to j\overset{a}\to i$ and objects $i$ of $I$, and invertible modifications
$$\xymatrix@C=40pt{c^*b^*a^*\ar@{=>}[d]_{\textstyle \chi a^*}\ar@{}[rd]|{\textstyle \overset{\textstyle \omega\!=\!\omega_{a,b,c}}\Rrightarrow}\ar@{=>}[r]^{\textstyle c^*\chi}&c^*(ab)^*\ar@{=>}[d]^{\textstyle \chi}\\(bc)^*a^*\ar@{=>}[r]_{\textstyle \chi}&(abc)^*\,,}\hspace{0.4cm}\xymatrix@C=70pt{a^*\ar@{=>}[d]_{\textstyle \iota a^*}\ar@{=>}[r]^{\textstyle a^*\iota}\ar@{=>}[dr]|{\textstyle  1}&a^*1_i^*\ar@{=>}[d]^{\textstyle \chi}\\ 1_j^*a^*\ar@{=>}[r]_{\textstyle \chi}\ar@{}[ru]|(0.3){\textstyle \overset{\textstyle \gamma\!=\!\gamma_a}\Rrightarrow}
\ar@{=>}[r]_{\textstyle \chi}\ar@{}[ru]|(0.7){\textstyle
\overset{\textstyle \delta\!=\!\delta_a}\Rrightarrow}& a^*\,,
}$$

\noindent respectively associated to triplets of composible arrows
$\ell \overset{c}\to k\overset{b}\to j \overset{a}\to i$ and arrows $j\overset{a}\to i$
 of the category $I$, subject to the  two coherence axioms {\bf (CC1)} and {\bf (CC2)}, as stated in Section 7.

The lax diagram is termed {\em normal} or {\em unitary} whenever the following conditions
hold: i) for each object $i$ of $I$, $1_i^*=1_{\f_i}$ and $\iota_i=1_{1_{\f_i}}$; ii) for each arrow $a:j\to i$ of $I$,
$\chi_{a,1_j}=1_{a^*}=\chi_{1_i,a}$ and the modifications $\gamma_a$ and $\delta_a$ are the unique coherence isomorphisms.

Note that a lax functor $\f:I^{\mathrm{op}}\to \bicat$ consists of
the same data as above, for a lax diagram, but with the difference that the modifications
 $\omega$, $\gamma$, and $\delta$ are no longer required to be invertible. However, we need
 lax diagrams of bicategories as above in order for the Grothendieck construction on
  them, as  shown in the next Section 3,  to give rise to  bicategories.

  A {\em diagram of bicategories} is a functor $\f:I^{\mathrm{op}}\to \mathbf{ Hom}\subset \bicat$  to the category {\bf  Hom} of
bicategories and homomorphisms, that is, a lax diagram where each of the pseudo transformations $\chi$ and $\iota$
are identities and the modifications $\omega$, $\gamma$, and $\delta$ are given by
the unique coherence isomorphisms.

A {\em pseudo-diagram of bicategories} is a
 trihomomorphism, or pseudo functor, $\f:I^{\mathrm{op}}\to \bicat$, that is, a lax diagram whose data $\chi$ and $\iota$ are pseudo natural equivalences.

 A {\em lax diagram of categories}, that is,  a lax functor $\f:I^{\mathrm{op}}\to \cat$ to
the 2-category $\cat$ of small categories, is the same thing as a lax diagram of bicategories in which
every bicategory $\f_i$, $i\in \mbox{Ob}I$, is a category
(i.e., a bicategory where all the 2-cells are identities) since this
condition forces all the modifications $\omega$, $\delta$, and $\gamma$ to be identities.

For any given category $I$, the lax diagrams of bicategories $\f:I^{\mathrm{op}}\to \bicat$ are the
objects of a  tricategory, denoted as \begin{equation}\label{bicatiop}\bicat^{ I^{\mathrm{op}}},\end{equation}
whose 1-cells, called here {\em lax $I$-homomorphisms}, are lax transformations
{\em all of whose coherence $3$-cells are invertible}, whose 2-cells, called {\em pseudo $I$-transformations}, are trimodifications,  and whose  3-cells, called {\em
$I$-modifications}, are perturbations, in the sense of \cite[3.3]{g-p-s}.
Then, the data for a lax $I$-homomorphism $$F=(F,\theta,\Pi,\Gamma):\f\to \f'$$
are comprised of: for $i$ an object of $I$, a homomorphism  $F_i:\f_i\to \f'_i$,
 for $a:j\to i$ a morphism of $I$, a pseudo transformation
 $$
 \xymatrix@R=14pt{\f_i\ar[d]_{\textstyle a^*}\ar[r]^{\textstyle F_i}\ar@{}[rd]|{\textstyle \overset{\textstyle \theta\!=\!\theta_a}\Rightarrow}&\f'_i\ar[d]^{\textstyle a^*} \\
 \f_j\ar[r]_{\textstyle F_j}&\f'_j\,,
 }
 $$
for  $k\overset{b}\to j\overset{a}\to i$ two composable arrows and  $j$ any object of
 $I$, the respective invertible modifications

$$\xymatrix@R=8pt@C=16pt{F_kb^*a^*\ar@2{->}[dd]_{\textstyle \theta a^*}\ar@{=>}[r]^{\textstyle F_k\chi} & F_k(ab)^*\ar@{=>}[rd]^{\textstyle \theta}&\\\ar@{}[rr]|(0.4){\overset{\textstyle \Pi\!=\!\Pi_{a,b}}\Rrightarrow} && (ab)^*F_i\\
b^*F_ja^*\ar@{=>}[r]_{\textstyle b^*\theta}&b^*a^*F_i\ar@{=>}[ur]_{\textstyle \chi'F_i}&
}\hspace{0.6cm}
\xymatrix@R=12pt@C=16pt{& F_j\ar@2{->}[rdd]^{\textstyle \iota'F_j}\ar@2{->}[ldd]_{\textstyle F_j\iota}\ar@{}[dd]|(.6){\textstyle{\overset{\Rrightarrow}{\Gamma=\Gamma_{j}}}} &\\ & & \\ F_j1_j^*
\ar@2{->}[rr]_{\textstyle \theta} && 1_j^*F_j\,, }
$$
and these are subject to the two coherence axioms {\bf (CC3)} and {\bf (CC4)}, as stated in Section 7.

When the pseudo-transformations  $\theta:F_ja^*\Rightarrow a^*F_i$ are
pseudo-natural equivalences, for
all arrows $a:j\to i$, then $F:\f\rightarrow \f'$ is termed a {\em pseudo
$I$-homomorphism}.

Given lax $I$-homomorphisms $F,F':\f\to \f'$, a {\em pseudo $I$-transformation}
 between them, $$m=(m,\mathrm{M}):F\Rightarrow F'$$ is merely a trimodification,
 so it consists of pseudo transformations
 $$\xymatrix@R=0pt@C=2pt{\f_i\ar@/^0.9pc/[rr]^{\textstyle F_i}
\ar@/_0.9pc/[rr]_{\textstyle F'_i}&\Downarrow\!m=m_i& \f'_i},$$ $i\in\mbox{Ob}I$, and  invertible modifications
$$\xymatrix{F_ja^*\ar@{=>}[r]^{\textstyle m a^*}\ar@{=>}[d]_{\textstyle \theta}\ar@{}[rd]|{\textstyle \underset{\Rrightarrow}{\textstyle \mathrm{M}\!=\!\mathrm{M}_a}}&F'_ja^*\ar@{=>}[d]^{\textstyle \theta'}\\
a^*F_i\ar@{=>}[r]_{\textstyle a^*m}&a^*F'_i,}
$$
\noindent one  for each arrow $a:j\to i$ of $I$, subject to the two coherence conditions {\bf (CC5)} and {\bf (CC6)}, as stated in Section 7. And, finally, say that if $m,m':F\Rightarrow F':\f\to \f'$ are pseudo $I$-transformations, then an {\em $I$-modification}
$\sigma:m\Rrightarrow m'$ is a family of modifications
$$\sigma_i:m_i\Rrightarrow m'_i:F_i\Rightarrow F'_i:\f_i\to \f'_i\,,$$
one for each object $i$ of $I$, subject to the coherence condition {\bf (CC7)}.

For lax $I$-diagrams of bicategories $\f$ and $\g$, compositions in $\bicat^{ I^{\mathrm{op}}}\!(\f,\g)$ are as follows:
2-cells $\sigma:m\Rrightarrow m'$ and $\sigma':m'\Rrightarrow m''$, where
$m,m',m'':F\Rightarrow F'$, $F,F':\f\to \g$, are vertically composed yielding the $I$-modification
$\sigma'\cdot\sigma:m\Rrightarrow m''$ such that, for any object $i$ of $I$,
$(\sigma'\cdot\sigma)_i=\sigma'_i\cdot \sigma_i:m_i\Rrightarrow m''_i$. The horizontal composition
${m'\circ m:F\Rightarrow F''}$ of 1-cells $F\overset{m}\Rightarrow
F'\overset{m'}\Rightarrow F'':\f\to \g$, is given by writing $(m'\circ m)_i=m'_i\circ m_i$ for
each $i\in\mbox{Ob}I$, while its component  at an arrow $a:j\to i$ is the modification obtained
by pasting the diagram
$$\xymatrix{F_ja^*\ar@{=>}[r]^{\textstyle ma^*}\ar@{=>}[d]_{\textstyle \theta}\ar@{}[rd]|{\textstyle \underset{\Rrightarrow}{\textstyle \mathrm{M}}}&F'_ja^*\ar@{=>}[d]_{ \theta'}\ar@{=>}[r]^{\textstyle m'a^*}\ar@{}[rd]|{\textstyle \underset{\Rrightarrow}{\textstyle \mathrm{M'}}}&F''_ja^*\ar@{=>}[d]^{\textstyle \theta''}\\
a^*F_i\ar@{=>}[r]_{\textstyle a^*m}&a^*F'_i\ar@{=>}[r]_{\textstyle a^*m'}&a^*F''_i\,.}
$$
Two $I$-modifications $\sigma:m\Rrightarrow n:F\Rightarrow F'$ and
${\sigma':m'\Rrightarrow n':F'\Rightarrow F''}$ compose horizontally  giving the $I$-modification
${\sigma'\circ\sigma:m'\circ m\Rrightarrow n'\circ n}$ such that,
 for any object $i$ of $I$,
$(\sigma'\circ \sigma)_i=\sigma'_i\circ\sigma_i:m'_i\circ m_i\Rrightarrow n'_i\circ n_i$. All the structure constraints in the bicategory $\bicat^{ I^{\mathrm{op}}}\!(\f,\g)$ are provided by using
 the
corresponding structure constraints of the tricategory $\mathbf{Bicat}$ in a pointwise fashion. Thus, for example, For pseudo $I$-transformations $F\overset{m}\Rightarrow
F'\overset{m'}\Rightarrow F''\overset{m''}\Rightarrow F''':\f \to \g$, the $I$-modification
$\boldsymbol{a}:m''\circ (m'\circ m)\Rrightarrow (m''\circ m')\circ m$ is that defined by the family of associativity modifications
$\boldsymbol{a}:m''_i\circ(m'_i\circ m_i)\Rrightarrow (m''_i\circ m'_i)\circ m_i$ of
$\bicat(\f_i,\g_i)$, $i\in \mbox{Ob}I$.

For lax $I$-diagrams of bicategories $\f$, $\g$, and $\h$, the composition homomorphism
\begin{equation}\label{e1.2.9}
\bicat^{ I^{\mathrm{op}}}\!(\g,\h)\times \bicat^{ I^{\mathrm{op}}}\!(\f,\g)\to \bicat^{ I^{\mathrm{op}}}\!(\f,\h)
\end{equation}
carries  lax $I$-homomorphisms ${\f\overset{F}\to \g\overset{G}\to \h}$ to the lax $I$-homomorphism
${GF:\f\to \h}$, whose component at an object $i$ of $I$ is the composite homomorphism ${G_iF_i:\f_i\to \h_i}$, its
component at an arrow $a:j\to i$ is the composed pseudo transformation
$\xymatrix@C=10pt{G_jF_ja^*\ar@{=>}[r]^{G_j\theta}&G_ja^*F_i\ar@{=>}[r]^{\theta F_i}&a^*G_iF_i}$,
its
component at a pair of composable arrows $k\overset{b}\to j$ and $j\overset{a}\to i$ is the modification obtained, from those
 of $F$ and $G$, by pasting the diagram
$$
\xymatrix@C=10pt@R=5pt{G_kF_kb^*a^* \ar@2{->}[rr]^{\textstyle G_k\theta a^*} \ar@2{->}[ddd]_{\textstyle G_kF_k\chi} && G_kb^*F_ja^* \ar@2{->}[rr]^{\textstyle \theta F_j a^*} \ar@2{->}[ddd]_(0.4){G_kb^*\theta} && b^*G_jF_ja^* \ar@2{->}[ddd]^{\textstyle b^*G_j\theta}\\
& \ar@3{->}[d]_(0.4){\textstyle G_k\Pi}&& \overset{(\ref{4})}\cong &\\
&&&&\\
G_kF_k(ab)^* \ar@2{->}[ddddrr]_{\textstyle G_k\theta} && G_kb^*a^*F_i \ar@2{->}[rr]^{\theta a^*F_i}\ar@2{->}[dddd]_(0.4){G_k\chi F_i} && b^*G_ja^*F_i \ar@2{->}[dd]^{\textstyle b^*\theta F_i}\\ &&&\ar@3{->}[d]_{\textstyle ~~~\Pi F_i}&
\\
&&&& b^*a^*G_iF_i\ar@2{->}[dd]^{\textstyle \chi G_iF_i}\\
&&&&\\
&& G_k(ab)^*F_i \ar@2{->}[rr]_{\textstyle \theta F_i} && (ab)^*G_iF_i }$$
and, finally, its component $\Gamma$ at an object $j$ of $I$ is the modification obtained from those
of $F$ and $G$, by pasting the diagram

$$
\xymatrix@C=37pt@R=9pt{&G_jF_j1_j^*\ar@2{->}[r]^{\textstyle G_j \theta}\ar@3{->}[d]^(0.7){\textstyle ~G_j\Gamma}&G_j1_j^*F_j\ar@2{->}[rddd]^{\textstyle \theta F_j}&\\
&&\ar@3{->}[d]_{\textstyle \Gamma F_j~}&\\
&&&\\
 G_jF_j\ar@2{->}[uuurr]_{G_j\iota F_j}\ar@2{->}@<3pt>[uuur]^{\textstyle G_jF_j \iota}\ar@2{->}@<3pt>[rrr]_{\textstyle \iota G_jF_j}&&&
  1_j^*G_jF_j.}$$
     If $m:F\Rightarrow F':\f\to \g$ and
$n:G\Rightarrow G':\g\to \h$ are pseudo $I$-transformations, then their composition
is $nm:GF\Rightarrow G'F':\f\to \h$, whose component at an object $i$ is $n_im_i:G_iF_i\Rightarrow
G'_iF'_i:\f_i\to \h_i$, and whose component at an arrow $a:j\to i$ is the modification obtained by pasting the
diagram

$$
\xymatrix@C=4pt@R=1pt{G_jF_ja^*\ar@2{->}[rr]^{\textstyle G_jm_ja^*}\ar@2{->}[dd]_{\textstyle G_j\theta} && G_jF'_ja^*\ar@2{->}[rr]^{\textstyle n_jF'_ja^*}\ar@2{->}[dd]|-{G_j\theta'} && G'_jF'_ja^* \ar@2{->}[dd]^{\textstyle G'_j\theta}\\
& \underset{\Rrightarrow}{G_j\mathrm{M}} && \overset{(\ref{4})}\cong &\\
G_ja^*F_i\ar@2{->}[dd]_{\textstyle \theta F_i}\ar@2{->}[rr]^{G_ja^*m_i} && G_ja^*F'_i \ar@2{->}[dd]|-{\theta F'_i}\ar@2{->}[rr]^{n_ja^*F'_i} && G'_ja^*F'_i \ar@2{->}[dd]^{\textstyle \theta' F'_i} \\
& \overset{(\ref{4})}\cong && \underset{\Rrightarrow}{NF'_i}&\\
a^*G_iF_i \ar@2{->}[rr]_{\textstyle a^*G_im_i} && a^*G_iF'_i \ar@2{->}[rr]_{\textstyle a^*n_iF'_i} && a^*G'_iF'_i\\  }
$$

And the composition of $I$-modifications $\sigma:m\Rrightarrow m':F\Rightarrow
F':\f\to \g$ and
$\tau:n\Rrightarrow n':G\Rightarrow G':\g\to \h$ is $\sigma\tau:nm\Rrightarrow n'm'$, with
$(\sigma\tau)_i=\sigma_i\tau_i:n_im_i\Rrightarrow n'_im'_i$ for every $i\in \mbox{Ob}I$. Moreover, given lax $I$-homomorphisms $\f\overset{F}\to \g\overset{G}\to\h$ and
 pseudo
$I$-transformations $F\overset{m}\Rightarrow F'\overset{m'}\Rightarrow F'':\f \to \g$ and
$G\overset{n}\Rightarrow G'\overset{n'}\Rightarrow G'':\g\to\h$, the structure constraint  of the composition homomorphism (\ref{e1.2.9}), $1_{GF}\Rrightarrow 1_{G}1_{F}$ and $ n'm'\circ nm\Rrightarrow
(n'\circ n)(m'\circ m)$, are
provided by the family of modifications (\ref{e1.5}) $1_{G_iF_i}\Rrightarrow 1_{G_i}1_{F_i}$ and $ n'_im'_i\circ n_im_i\Rrightarrow
(n'_i\circ n_i)(m'_i\circ m_i)$,
 $i\in \mbox{Ob}I$, respectively.

The associativity and unit pseudo natural equivalences
$$
\xymatrix@R=0pt@C=-4pt{\bicat^{ I^{\mathrm{op}}}\!\!\!(\h,\T)\!\!\times\! \bicat^{ I^{\mathrm{op}}}\!\!(\g,\h)\!\!\times\!
\bicat^{ I^{\mathrm{op}}}\!\!(\f,\g)\ar[rr]\ar[dd]&&\bicat^{ I^{\mathrm{op}}}\!\!(\g,\T)\!\!\times\!\bicat^{ I^{\mathrm{op}}}\!\!(\f,\g)
\ar[dd]\\&\boldsymbol{a}\Rightarrow& \\
\bicat^{ I^{\mathrm{op}}}\!(\h,\T)\times \bicat^{ I^{\mathrm{op}}}\!(\f,\h)\ar[rr]&&\bicat^{ I^{\mathrm{op}}}\!(\f,\T)
}
$$
$$
\xymatrix@R=0pt@C=-8pt{\bicat^{ I^{\mathrm{op}}}\!\!(\g,\!\g)\!\!\times\!\!\bicat^{ I^{\mathrm{op}}}\!\!(\f\!,\!\g)\ar[ddrr]&&
\bicat^{ I^{\mathrm{op}}}\!\!(\f\!,\!\g) \ar[ll]_-{1_\g\times
1}\ar[rr]^-{1\times 1_\f}\ar[dd]^1&&\bicat^{ I^{\mathrm{op}}}\!\!(\f\!,\!\g)\!\!\times\!\! \bicat^{ I^{\mathrm{op}}}\!\!(\f\!,\!\f)\ar[ddll]\\
&\hspace{1.5cm}\overset{\textstyle \boldsymbol{l}}\Rightarrow \hspace{-0.75cm}&
&\hspace{-0.5cm}\overset{\textstyle \boldsymbol{r}}\Leftarrow\hspace{1cm} &\\&&\bicat^{ I^{\mathrm{op}}}\!\!(\f\!,\!\g)&&}
$$
are as follows: For any lax $I$-homomorphisms $\f\overset{F}\to \g \overset{G}\to\h
\overset{H}\to \mathcal{K}$, $H(GF)\overset{\textstyle \boldsymbol{a}}\Rightarrow (HG)F$ is the pseudo $I$-equivalence whose component at any object $i$ of $I$ is the identity on the homomorphism $H_iG_iF_i$, and whose component at an arrow $a:j\to i$ is the modification obtained by pasting
$$
\xymatrix{H_jG_jF_ja^*\ar@{=>}[r]^{\textstyle 1}\ar@{=>}[dd]_{\textstyle H_j(\theta F_i\circ G_j\theta)}&H_jG_jF_ja^*\ar@{=>}[r]^{\textstyle H_jG_j\theta}\ar@{=>}@<-0.5ex>[rd]_(0.35){ H_j(\theta F_i\circ G_j\theta)\ \ \ \ }\ar@{}@<2.2ex>[rd]|(0.6){\textstyle \ \cong\! \widehat{H}_j}&H_jG_ja^*F_i\ar@{=>}[d]^{\textstyle H_j\theta F_i}\\&&H_ja^*G_iF_i\ar@{=>}[lld]_{\textstyle 1}\ar@{=>}[d]^{\textstyle \theta G_iF_i}\\ \ \ H_ja^*G_iF_i\ar@{=>}[r]_{\textstyle \theta G_iF_i}\ar@{}[ruu]|{\textstyle \cong}&a^*H_iG_iF_i\ar@{=>}[r]_{\textstyle a^*1}\ar@{}[ru]|(0.4){\textstyle \cong}&a^*H_iG_iF_i.}
$$
Besides, for any pseudo $I$-transformations $(m,n,t):(H,G,F)\Rightarrow
(H',G',F')$, the corresponding $I$-modification
$\widehat{\boldsymbol{a}}:\boldsymbol{a}_{_{H'\!,G'\!,F'}}\circ m(nt)\Rrightarrow (mn)t\circ \boldsymbol{a}_{_{H,G,F}}$, is given by the family of modifications (\ref{asso}) $\widehat{\boldsymbol{a}}:1\circ m_i(n_it_i)\Rrightarrow
(m_in_i)t_i\circ 1$, $i\in \mbox{Ob}I$.
For any $I$-homomorphism $F:\f\to \g$,  $\boldsymbol{\l}:1_{\g}F\Rightarrow F$ and $\boldsymbol{r}:F1_{\f}\Rightarrow F$,
are the pseudo $I$-equivalences whose components at any object $i\in \mbox{Ob}I$
are both the identity on the homomorphism $F_i$  and, at an arrow $a:j\to i$, are the
canonical isomorphism $a^*1_{F_i}\circ \theta_a\cong \theta_a\circ 1_{F_ja^*}$. Besides,  for $m:F\Rightarrow F'$ any pseudo $I$-transformation, the corresponding
 $I$-modifications $\widehat{\boldsymbol{l}}:\boldsymbol{l}_{F'}\circ 1_{1_\g}m
 \Rrightarrow m\circ \boldsymbol{l}_F$ and $\widehat{\boldsymbol{r}}:\boldsymbol{r}_{F'}
 \circ m 1_{1_\f}  \Rrightarrow m\circ \boldsymbol{r}_F$ are respectively given by the family of
 modifications (\ref{unibicat}), $\widehat{\boldsymbol{l}}:1_{F'_i}\circ
 1_{1_{\g_i}}m_i\Rrightarrow m_i\circ 1_{F_i}$ and $\widehat{\boldsymbol{r}}:1_{F'_i}
 \circ m_i
1_{1_{\f_i}}\Rrightarrow m_i\circ 1_{F_i}$.

In $\bicat^{ I^{\mathrm{op}}}$, the structure invertible $I$-modifications $\pi$ and $\mu$,
as in the definition of a tricategory, for any lax $I$-homomorphisms
$\f\overset{F}\to \g \overset{G}\to\h\overset{H}\to \mathcal{K}\overset{K}\to \mathcal{T}$,
$$\begin{array}{rcl}
(\boldsymbol{a}_{K,H,G}1_F\circ \boldsymbol{a}_{K,HG,F})\circ 1_K\boldsymbol{a}_{H,G,F}&
\overset{\textstyle \pi}\Rrightarrow &\boldsymbol{a}_{KH,G,F}\circ \boldsymbol{a}_{K,H,GF},
\\[5pt] \boldsymbol{r}_G1_F\circ \boldsymbol{a}_{G,1_\g,F}&\overset{\textstyle \mu} \Rrightarrow &
1_G\boldsymbol{l}_F,\end{array}
$$
are given by the family of modifications $(\ref{pimu})$
$\pi_{_{K_i,H_i,G_i,F_i}}$ and $\mu_{_{G_i,F_i}}$,
$i\in \mbox{Ob}I$.

Finally, note that considering lax $I$-diagrams of categories, that is,  lax functors ${\f:I^{\mathrm{op}}\to \cat}$ to
the 2-category $\cat$ of small categories, then $$\cat^{ I^{\mathrm{op}}}\subseteq \bicat^{ I^{\mathrm{op}}}$$ is a full subtricategory of $\bicat^{ I^{\mathrm{op}}}$. But note that $\cat^{ I^{\mathrm{op}}}$ is actually a 2-category, since all its 3-cells are identities.

\section{The bicategorical Grothendieck construction}
\subsection{The Grothendieck construction on lax diagrams of bicategories}Let $I$ be a small category. The well-known {\em Grothendieck construction} on a lax diagram of categories $I^\mathrm{op}\to \cat$ \cite{grothendieck, g-j, jardine, thomason} admits an extension to a lax diagram of bicategories
$$\f=(\f,\chi,\iota, \omega,\gamma,\delta): I^\mathrm{op}\to \bicat,$$
 assembling it into a large bicategory
\begin{equation}\label{gc}\xymatrix{\int_I\f}\end{equation}
which is a lax colimit  of the bicategories $\f_i$, $i\in \mbox{Ob}I$ and, as we shall detail later, it can be thought as its  homotopy colimit. This
bicategory is defined as follows (cf. \cite{bakovik, ccg}):

\underline{The objects} of  $\xymatrix{\int_I\f}$ are pairs $(x,i)$, where $i$ is an object of $I$ and $x$ one of the bicategory $\f_i$, so that
$$\xymatrix{\mbox{ Ob}\!\int_I\f=\bigsqcup\limits_{i\in \mbox{\scriptsize{ Ob}}I}\hspace{-5pt}\mbox{ Ob}\f_i.}$$

\underline{The hom-categories}  are
$$\xymatrix{\int_I\f\big((y,j),(x,i)\big)=\bigsqcup\limits_{j\overset{a}\to i}\f_j(y,a^*x),}$$
where the disjoint union is over all arrows $a:j\to i$ in $I$. Then, a morphism $(u,a):(y,j)\to (x,i)$ in $\int_I\f$ is a pair of morphisms where $a:j\to i$ is  in $I$ and $u:y\to a^*x$ is in $\f_j$; and given two morphisms $(u,a),(u',a'):(y,j)\to(x,i)$, the existence of a $2$-cell  $(u,a)\Rightarrow (u',a')$ requires that $a=a'$,  and then, such a $2$-cell
$\xymatrix@C=4pt{(y,j)  \ar@/^1.1pc/[rr]^{\textstyle (u,a)} \ar@/_1.1pc/[rr]_{\textstyle (u',a)} & {}_{\textstyle \Downarrow\!(\alpha,a)} &(x,i) }  $ consists of a $2$-cell  $\xymatrix @C=0.5pc {y  \ar@/^1pc/[rr]^{\textstyle\  u} \ar@/_1pc/[rr]_{\textstyle \ u'} & {}_{\textstyle\Downarrow\alpha} &a^*x }$ in $\f_j$.

\underline{The horizontal composition} functor
$$\bigsqcup\limits_{j\overset{a}\to i}\f_j(y,a^*x) \times \bigsqcup\limits_{k\overset{b}\to j}\f_k(z,b^*y) \overset{\circ}\longrightarrow
\bigsqcup\limits_{k\overset{c}\to i}\f_k(z,c^*i),$$ for each triplet of objects $(z,k)$, $(y,j)$, and $(x,i)$ of $\int_I\f$,  maps the component at two morphisms $a:j\to i$ and $b:k\to j$ of $I$ into the component at the composite $ab:k\to i$ via the composition
$$\xymatrix{\f_j(y,a^*\!x)\!\!\times\! \f_k(z,b^*\!y) \ar[r]^-{\textstyle b^*\times 1} &
\f_k(b^*y,b^*a^*x)\!\times\! \f_k(z,b^*y)\ar[d]^-{\textstyle \circ}&\\
&\f_k(z,b^*a^*x)\ar[r]^{\textstyle \chi_{_*}}&\f_k(z,(ab)^*x), }$$
where $\chi=\chi_{_{a,b}}x:b^*a^*x\to (ab)^*x$, that is,
$$\xymatrix@C=1pt{(z,k)  \ar@/^1.3pc/[rr]^{\textstyle (v,b)} \ar@/_1.3pc/[rr]_{\textstyle (v',b)} & {}_{\textstyle \Downarrow\!(\beta,b)} &(y,j)  \ar@/^1.3pc/[rr]^{\textstyle (u,a)} \ar@/_1.3pc/[rr]_{\textstyle (u',a)} & {}_{\textstyle \Downarrow\!(\alpha,a)} &(x,i)&\ar@{|->}[rrr]^-{\textstyle \circ} && &&(z,k) \ar@/^1.5pc/[rr]^{\textstyle (\chi\!\circ\!(b^*u\!\circ\! v),ab)} \ar@/_1.5pc/[rr]_{\textstyle (\chi\!\circ\!(b^*u'\!\circ\! v'),ab)}&{}_{\textstyle \Downarrow\!\big(1_{\chi}\!\circ\!(b^*\alpha\!\circ\!\beta),ab\big)}& (x,i).  }  $$

\underline{The structure associativity}  isomorphism
$$(u,a)\circ ((v,b)\circ(w,c)) \cong ((u,a)\circ (v,b))\circ (w,c),$$
 for any three composable morphisms $(t,\ell)\overset{\textstyle (w,c)}\longrightarrow (z,k)\overset{\textstyle (v,b)}\longrightarrow  (y,j)\overset{\textstyle (u,a)}\longrightarrow  (x,i)$ in
$\int_I\c$,
is provided by pasting, in the bicategory $\f_\ell$, the diagram
$$\xymatrix@C=6pt@R=1.5pt{
&&&&& (bc)^*y \ar@{}[ddr]|{\textstyle \overset{\textstyle \widehat{\chi}_u}\cong }\ar[rr]|{\,(bc)^*\!\!u} && (bc)^*a^*x \ar[rdd]^{\textstyle \chi} &\\
&&&\cong &&& &&\\
t \ar[rr]^w \ar@/^2.5pc/[rrrrrrruu]^(0.5){\textstyle (bc)^*\!u\!\circ\! (\chi \circ (c^*\!v\circ
w))}\ar@/_1.5pc/[rrrrrrrdd]_(.45){\textstyle c^*(\chi\circ(b^*\!u\circ v))\circ w} && c^*z\ar[rr]^{c^*\!v} && c^*b^*y  \ar[rr]^{c^*\!b^*\!u}
\ar[ruu]^{\chi} && c^*b^*a^*x \ar[ruu]^{\chi
a^*}\ar[rdd]_{c^*\chi} & \overset{\textstyle \omega}\cong & (abc)^*x\,. \\
&&&&& \cong &&&\\
&&&&&&& c^*(ab)^*x \ar[ruu]_{\textstyle \chi} }
$$

\underline{The identity morphism}, for each object $(x,i)$ in $\int_I\f$, is provided by the pseudo-transformation $\iota:1_{\f_i}\Rightarrow 1_i^*$, by $$ 1_{(x,i)}=(\iota x,1_i):(x,i)\to (x,i).$$

\underline{The left and right identity constraints}
$$\begin{array}{c}
1_{(x,i)}\circ (u,a)=(\chi \circ (a^*\iota x\circ u),a)\cong (u,a),\\[6pt]
(u,a)\circ 1_{(y,j)}=(\chi \circ(1_j^*u\circ \iota y),a)\cong (u,a),
\end{array}
$$
for each morphism
$(u,a):(y,j)\to (x,i)$,
are respectively given by pasting the diagrams
$$ \xymatrix@C=2pt@R=-1pt{ y\ar[rrr]^{\textstyle u} \ar@/_/[rrrddd]_{\textstyle u} &&& a^*x\ar[rr]^{\textstyle a^*\iota}\ar[ddd]_{1}&& a^*1_i^*x \ar@/^/[llddd]^{\textstyle \chi}
&&& y
\ar[rrrr]^{\textstyle \iota} \ar[rrrrdd]^{ u} \ar@/_1pc/[rrrrdddd]_{\textstyle u} &&&& 1_j^*y \ar[rrr]^{\textstyle 1_j^*u}\ar@{}[dd]|(.4){\overset{\textstyle \widehat{\iota}_u}\cong} &&& 1_j^*a^*x \ar@/^1pc/[llldddd]^{\textstyle \chi}\\
& &\cong  && \overset{\textstyle \delta}\cong &&&&&&&&  &&& \\
&&&&&&&&&&&\cong& a^*x \ar[rrruu]^-{\iota a^*} \ar[dd]_{1} &\overset{\textstyle \gamma}\cong&&\\
&&& a^*x &&&&&& && &&  &&\\
&&&&&&&&&&&& a^*x &&&\\ }$$

The coherence pentagon for associativity in $\int_I\f$ holds thanks to the coherence condition {\bf (CC1)} in Section 7, and the coherence triangles for unit constraints in $\int_I\f$ follows from {\bf (CC2)}. Hence $\int_I\f$ is actually a bicategory.

\vspace{0.2cm}
\subsection{The Grothendieck construction trihomomorphism}
The assignment $$\xymatrix{\f\mapsto \int_I\f}$$ is the function on objects of a trihomomorphism of tricategories \begin{equation}\xymatrix{\int_I :\bicat^{ I^{\mathrm{op}}}\rightarrow \bicat,}\end{equation}
described below.

\underline{The homomorphism} of bicategories
\begin{equation}\label{int2}\xymatrix{\int_I :\bicat^{ I^{\mathrm{op}}}(\f,\g)\rightarrow \bicat(\int_I\f,\int_I\g)},\end{equation}
for any two lax $I$-diagrams $\f,\g:I^\mathrm{op}\to \bicat,$ carries a lax $I$-homomorphism  $F=(F,\theta,\Pi,\Gamma):\f\to \g$ to the  homomorphism
\begin{equation}\label{G.7} \xymatrix{\int_I\!F:\int_I\f\rightarrow \int_I\g }\end{equation}
defined on objects by  $\xymatrix{\int_I\!F(x,i)=(F_ix,i)}$, and, for each pair of objects $(y,j)$ and $(x,i)$ of $\int_I\f$, the functor
$$\xymatrix@C=16pt{\int_I\!F:\bigsqcup\limits_{j\overset{a}\to i}\f_j(y,a^*x) \ar@<1.5ex>[r]& \bigsqcup\limits_{j\overset{a}\to i}\g_j(F_jy,a^*F_ix)}$$
is defined on the components at each morphism $j\overset{a}\to i$  by the composition of functors
$\f_j(y,a^*x) \overset{ F_j}\to \g_j(F_jy,F_ja^*x) \overset{ \theta_{_*}}\to\g_j(F_jy,a^*F_ix)$,
where $\theta=\theta_ax:F_ja^*x\to a^*F_ix$.

 If $(z,k)\overset{ (v,b)}\longrightarrow (y,j)\overset{ (u,a)}\longrightarrow (x,i)$ are any two composible morphisms of
$\int_I\f$, then the invertible structure 2-cell
$$\xymatrix{\int_I\!F(u,a)\circ\int_I\!F(v,b)\cong \int_I\!F\big((u,a)\circ(v,b)\big)}$$
is provided by pasting the diagram in $\f_k$
$$\xymatrix@C=6pt@R=-1pt{&&&&&&&&&&\\
&& & &{ b^*F_jy} \ar[rrddd]^{b^*\!F_ju}\ar@{}[ru]|>{\textstyle \cong} &&&& b^*a^*F_ix \ar[rrddd]^{\textstyle \chi F_i} && \\
&&&&&&&&&&\\
&&&&&&&&&&\\
F_kz \ar[rr]^{F_kv} \ar@/^3pc/[rrrrrrrruuu]^{\textstyle b^*(\theta\circ F_ju)\circ (\theta\circ F_kv)}
\ar@/_3pc/[rrrrrrrrddd]_{\textstyle F_k(\chi\circ(b^*\!u\circ v))} && F_kb^*y \ar[rruuu]^{\theta} \ar[rrddd]_{F_kb^*\!u} && \cong
\widehat{\theta}_u && b^*F_ja^*x \ar[rruuu]^{b^*\!\theta}&& \cong \Pi && (ab)^*F_ix \\
&&&&&&&&&&\\
&&&&&&&&&&\\
&&  & & F_kb^*\!a^*x \ar[rrrr]_{F_k\chi}\ar[rruuu]^{\theta a^*} \ar@{}[d]|(1.4){\textstyle \cong}&&&& F_k(ab)^*x \ar[rruuu]_{\textstyle \theta} && \\
&&&&&&&&&&}$$
where $\Pi=\Pi_{_{a,b}}x$. And, for each
object $(x,i)$ of $\int_I\f$, the isomorphism
$$\xymatrix{1_{\int_IF(x,i)} \cong \int_I\!F(1_{(x,i)})}$$
is provided by the invertible deformation $\Gamma_i$.

Since the commutativity coherence conditions follow from {\bf (CC3)} and {\bf (CC4)}, $\int_I\!F:\int_I\f\rightarrow \int_I\g$ is actually a homomorphism of bicategories. This describes the
function on objects of (\ref{int2}), which acts as follows on the hom-categories:
Any pseudo $I$-transformation, $m:F\Rightarrow G:\f\rightarrow \g$, gives rise to a pseudo-transformation
\begin{equation}\label{pstint}\xymatrix{ \int_I\!m:\int_I\!F\Rightarrow \int_I\!G:\int_I\f\to \int_I\g,}\end{equation}
whose component at an object $(x,i)$ of $\int_I\f$ is
$$\xymatrix{\int_I\!m(x,i)=(\iota_iG_ix\circ m_ix,1_i):(F_ix,i)\to (G_ix,i)},$$
and whose component at a morphism $(u,a):(y,j)\to (x,i)$
$$\xymatrix{\widehat{\int_I\!m}_{(u,a)}:\int_I\!m(x,i)\circ \int_I\!F(u,a)\cong \int_I\!G(u,a)\circ \int_I\!m(y,j)}$$
is given by pasting
$$\xymatrix@C=35pt@R=35pt{F_jy\ar[r]^{\textstyle F_ju}\ar[d]_{\textstyle m}&F_ja^*\!x\ar[r]^{\textstyle \theta}\ar[d]|{ ma^*}&a^*\!F_ix\ar[d]|{a^*\!m}\ar[rd]^{\textstyle a^*\!(\iota G_i\circ m)}&\\
\ar@{}[ru]|{\overset{\textstyle \widehat{m}_u}\cong} G_jy\ar[r]|{\, G_ju}\ar[d]_{\textstyle \iota G_j}&\ar@{}[ru]|{\overset{\textstyle \mathrm{M}}{\textstyle \cong}} G_ja^*\!x\ar[r]|{\, \theta}\ar[d]|{ \iota G_ja^*}&a^*\!G_ix\ar[r]|(0.42){\, a^*\!\iota G_i}\ar[d]|{\, \iota a^*\!G_i}\ar[rd]|{ 1}\ar@{}[ru]|(0.3){\textstyle \cong}&a^*\!1_i^*G_ix\ar[d]^{\textstyle \chi G_i}\\ \ar@{}[ru]|{\overset{\textstyle \widehat{\iota}_{G_ju}}{\textstyle \cong}}
1_j^*G_jy\ar[r]|(0.45){1_j^*\!G_ju}\ar@/_1.6pc/[rr]_{\textstyle 1_j^*\!(\theta \circ G_ju)}^{^{\textstyle \cong}}&\ar@{}[ru]|{\overset{(\ref{4})}{\textstyle \cong}}1_j^*G_ja^*\!x \ar[r]|{ 1_j^*\theta}&1_j^*a^*\!G_ix \ar[r]_{\textstyle \chi G_i}\ar@{}[ru]|(0.32){\overset{\textstyle \gamma G_i}{\textstyle \cong}}|(0.7){\overset{\textstyle \delta G_i}{\textstyle \cong}} &a^*\!G_ix\,.
}
$$
So defined, $\int_I\!m$ is indeed a pseudo transformation thanks to the coherence conditions {\bf (CC5)} and {\bf (CC6)}.

If $m,m':F\Rightarrow G:\f\to\g$ are two pseudo $I$-transformations, it follows from the commutativity
of the squares in {\bf (CC7)} that every $I$-modification $\sigma:m \Rrightarrow m'$ defines a
modification

$$\xymatrix{ \int_I\!\sigma: \int_I\!m \Rrightarrow \int_I\!m',}$$
by writing
$\xymatrix{\int_I\!\sigma(x,i)=(1_{\iota G_ix}\circ \sigma_ix, 1_i):(\iota G_ix\circ m_ix,1_i)\Rightarrow
(\iota G_ix\circ m'_ix,1_i)}$.

For $I$-modifications $\sigma:m\Rrightarrow m'$ and $\sigma':m'\Rrightarrow m''$, where
$m, m',m'':F\Rightarrow G$, the equality $\int_I\!(\sigma'\circ\sigma)=\int_I\!\sigma' \circ
\int_I\!\sigma$ is easily verified. Moreover, for the horizontal composition $n\circ m:F\Rightarrow H$ of
pseudo $I$-transformations $F\overset{m}\Rightarrow G \overset{n}\Rightarrow H:\f\to \g$, the invertible structure
 modification
\begin{equation}\label{gc1}\xymatrix{ \int_I\!n\circ \int_I\!m ~~\Rrightarrow ~~\int_I\!(n\circ m):\int_I\!F\Rightarrow \int_I\!H}\end{equation}
is given by pasting
$$
\xymatrix@R=30pt{&&\ar@{}[d]|(.75){\textstyle \cong}&&\\&1_i^*G_ix\ar[r]^{1_i^*n_i} \ar@/^2pc/[rr]^{\textstyle 1_i^*(\iota H_i\circ n_i)} &1_i^*H_ix\ar[r]^{1_i^*\iota H_i}\ar[rd]|{\ 1\ }\ar@{}@<1ex>[rd]^(.55){\textstyle \cong\!\delta}\ar@{}@<-1ex>[rd]_(.5){\textstyle\cong}&1_i^*1_i^*H_ix \ar[d]^{\textstyle \chi H_i}\\
F_ix\ar[r]_{\textstyle m_i}\ar[ru]^{\textstyle \iota G_i\circ m_i}\ar@{}@<-1ex>[ru]_{\textstyle =} &G_ix\ar[r]_{\textstyle n_i}\ar[u]|{\iota G_i}\ar@{}[ur]|{\textstyle \overset{(\ref{4})}\cong}&H_ix\ar[r]_{\textstyle \iota  H_i}\ar[u]|{\iota H_i}&1_i^*H_ix.
}
$$
If $m:F\Rightarrow G:\f\to\g$ is any pseudo
$I$-transformation, then $\int_I\!1_{m}=1_{\int_I\!m}$ and, for any lax $I$-homomorphism $F:\f\to
\g$, the invertible structure constraint
\begin{equation}\label{gc2}
\xymatrix{1_{\int_I\!F} ~~ \Rrightarrow ~~ \int_I\!1_F
}\end{equation}
is provided by the canonical isomorphisms $\boldsymbol{r}:\iota F_i\circ 1_{F_i} \cong \iota F_i$, $i\in \mbox{Ob}I$.

 \underline{The pseudo-natural equivalence}
$$\xymatrix@C=12pt@R=0pt{
\bicat^{ I^{\mathrm{op}}}(\g,\h)\times \bicat^{ I^{\mathrm{op}}}(\f,\g) \ar[rr]^-{\textstyle \int_I~\times~\int_I} \ar[dd] &&
\bicat(\int_I\!\g,\int_I\!\h)\times \bicat(\int_I\!\f,\int_I\!\g) \ar[dd]\\
&\Downarrow \Sigma &\\
\bicat^{ I^{\mathrm{op}}}(\f,\h) \ar[rr]^(0.6){\textstyle \int_I} && \bicat(\int_I\!\f,\int_I\!\h),\\ }$$ for any three lax $I$-diagrams $\f,\g$ and $\h$,  is that whose component at any pair of lax
$I$-homomorphisms, $\f\overset{F}\to \g \overset{G}\to \h$, is the pseudo natural equivalence
\begin{equation}\label{gc3}\xymatrix{
\Sigma=\Sigma_{G,F}:\int_I\!G\int_I\!F \Rightarrow \int_I\!GF,
}\end{equation}
which is the identity on objects, that is,
$$\Sigma(x,i)=1_{(G_i\!F_ix,i)}= (\iota_iG_iF_ix,1_i):(G_iF_ix,i)\longrightarrow (G_iF_ix,i)
, $$
and its component at a morphism $(u,a):(y,j)\to (x,i)$ is the 2-cell
$$\xymatrix{
\widehat{\Sigma}_{(u,a)}:1_{(G_iF_ix,i)}\circ \int_I\!G\!\int_I\!F(u,a)\Rightarrow \int_I\!GF(u,a)\circ 1_{(G_jF_jx,j)}},
$$
canonically obtained from the 2-cell $\int_I\!G\!\int_I\!F(u,a)\Rightarrow \int_I\!GF(u,a)$ given by the composition in $\h_j$
$$ \theta_aF_ix\circ G_j(\theta_ax\circ F_ju)\overset{1\circ \widehat{G}_j^{-1}}\cong  \theta_aF_ix\circ (G_j\theta_ax\circ G_jF_ju)\overset{\boldsymbol{a}}\cong  (\theta_aF_ix\circ G_j\theta_ax)\circ G_j F_ju .$$

For $(n, m):(G,F)\Rightarrow (G',F')$, the component of $\Sigma$ at $(n,m)$ is the invertible modification
$
\xymatrix{\Sigma_{G',F'}\circ \int_I\!n \!\int_I\!m \Rrightarrow \int_I\!nm \circ \Sigma_{G,F}},
$ canonically obtained from the modification $\int_I\!n \!\int_I\!m \Rrightarrow \int_I\!nm$, which assigns to each object $(x,i)$ of $\int_I\!\f$ the 2-cell of $\int_I\!\h$ provided by pasting in $\h_i$

$$\xymatrix@C=45pt@R=25pt{&G_i1_i^*F'_ix \ar@{}[rd]|(.4){\textstyle \underset{\textstyle \cong}{\Gamma F'_i}}\ar@{}[rdd]_(.65){ \underset{\textstyle \cong}{(\ref{4})}}\ar[r]^{\textstyle \theta F'_i }& 1_i^*G_iF'_ix\ar[rd]^{\textstyle\  \ \  1_i^*\!(\iota G'_iF'_i\!\circ\!n_iF'_i)} \ar[d]|{1_i^*n_iF'_i}& \\
G_iF_ix\ar@{}[r]|{\textstyle \cong}\ar[ru]^{\textstyle G_i(\iota F'_i\!\circ\!m_i)\ }\ar[rd]_{\textstyle G_im_i}&&1_i^*G'_iF'_ix\ar[r]|-{1_i^*\!\iota G'_iF'_i}\ar[rd]|{\ 1\ }\ar@{}[ru]|(.35){\textstyle \cong}&1_i^*1_i^*G'_iF'_ix \ar[d]^{\textstyle \chi G'_iF'_i}\\
&G_iF'_ix\ar[r]_{\textstyle nF'_i}\ar[uu]|{G_i\iota F'_i}\ar[ruu]|{\iota G_iF'_i}&G'_iF'_ix\ar[r]_{\textstyle \iota G'_iF'_i}\ar[u]|{\iota G'_iF'_i}\ar@{}[ru]|(.3){\textstyle \cong}\ar@{}[ru]|(.65){\textstyle \cong\!\delta}&1_i^*G'_iF'_ix.
}$$

\underline{The pseudo-natural equivalence}
\begin{equation}\label{gc4}\xymatrix{
\Sigma_\f:\int_I\! 1_\f \Rightarrow 1_{\int_I\!\f},}
\end{equation}
for any lax
$I$-diagram $\f:I^{\mathrm{op}}\to \bicat$, is the identity on objects and its component at any 1-cell, $(u,a):(y,j)\to (x,i)$ of $\int_I\!\f$, is the 2-cell $$\xymatrix{1_{(x,i)}\circ \int_I\! 1_\f(u,a)\Rightarrow (u,a)\circ 1_{(y,j)}}$$ obtained by pasting
$$
\xymatrix@C=40pt{(y,j)\ar@<-2pt>@{}[r]_{\textstyle \Downarrow\!(\boldsymbol{l}_u,a)}\ar[d]_{\textstyle 1}\ar[r]^{(\textstyle 1_{a^*\!x}\circ u,a)}\ar@/_1.8pc/[r]|{(u,a)}&(x,i)\ar[d]^{\textstyle 1}\\(y,j)\ar[r]_{\textstyle (u,a)}&(x,i).\ar@{}@<-1.5pt>[l]_{_{\textstyle \cong}}
}
$$

\underline{The structure invertible modifications} $\omega$, $\delta$ and $\gamma$,  as in the definition of a trihomomorphism, for  $\int_I$,
\begin{equation}\label{odg}
 \xymatrix@C=-17pt@R=15pt{ &&\int_I\!HG\int_I\!F \ar@2{->}[rrd]^{\textstyle \Sigma}&&\\
 \int_I\!H\int_I\!G\int_I\!F \ar@{}[rrr]|(0.7){\textstyle \overset{\textstyle \omega}\cong}\ar@2{->}[urr]^{\textstyle \Sigma \,1} \ar@2{->}[dr]_{\textstyle 1\,\Sigma}&&&& \int_I\!(HG)F \\
 &\int_I\!H\int_I\!GF \ar@2{->}[rr]_{\textstyle \Sigma}&&\int_I\!H(GF) \ar@{=>}[ur]_{\textstyle \int_I\!\boldsymbol{a}}&}
  \hspace{-0.7cm}
 \xymatrix@C=6pt@R=15pt{\int_I\!F\int_I\!1_\f \ar@{=>}[ddrr]_{\textstyle 1\Sigma_\f}\ar@{=>}[rr]^{\textstyle \Sigma_{F,1_\f}}&&
 \int_I\!F1_\f \ar@{}[lldd]|(0.3){\textstyle \overset{\textstyle \delta}\cong}\ar@{=>}[dd]^{\textstyle \int_I\!\boldsymbol{r} }\\&&\\&&\int_I\!F}
  \hspace{-0.1cm}
   \xymatrix@C=6pt@R=15pt{\int_I\!1_\g\int_I\!F \ar@{=>}[ddrr]_{\textstyle \Sigma_\g 1}\ar@{=>}[rr]^(0.55){\textstyle \Sigma_{1_\g,F}}&&
 \int_I\!1_\g F\ar@{}[ddll]|(0.3){\textstyle \overset{\textstyle \gamma}\cong}\ar@{=>}[dd]^{\textstyle \int_I\!\boldsymbol{l} }\\&&\\ &&\int_I\!F}
 \end{equation}
 for any lax I-homomorphisms $\f\overset{F} \to \g \overset{G}\to \h \overset{H}\to \mathcal{T}$,
 are, respectively, the unique coherence 2-cells, $(x,i)\in \mbox{Ob}\!\int_I\!\f$,
 $$
\xymatrix@C=10pt{&(H_iG_iF_ix,i)\ar[r]^{\textstyle 1}&(H_iG_iF_ix,i)\ar[r]^{\textstyle1}&(H_iG_iF_ix,i)\\
(H_iG_iF_ix,i)\ar[ur]^{\textstyle \int_I\!H\int_I\!G\,1_{(F_ix,i)}\,\,\,\,\,\,} \ar[dr]_{\textstyle \int_I\!H\,1_{(G_iF_ix,i)}}&\ar@{}[r]|{\textstyle\cong}&&\\
&(H_iG_iF_ix,i)\ar[r]^{\textstyle 1}&(H_iG_iF_ix,i)\ar[r]^{\textstyle 1}&(H_iG_iF_ix,i)\ar[uu]_{\textstyle (\iota_iH_iG_iF_ix\!\circ\! 1_{H_iG_iF_ix},1_i)}\ar@/^3pc/[uu]^{\textstyle 1}_(.6){\textstyle\overset{\textstyle (\boldsymbol{r},1_i)}\cong},
}
 $$
 $$\xymatrix@C=20pt@R=35pt{
  (F_ix,i)\ar[d]_{\textstyle \int_i\!F\,1_{(x,i)}}\ar[rr]^{\textstyle 1}&&(F_ix,i)\ar[d]^{\textstyle (\iota_iF_ix\circ 1_{F_ix},1_i)}\ar@/_2.5pc/[d]_{\textstyle 1}^{\textstyle\overset{\textstyle (\boldsymbol{r},1_i)}\cong}\\(F_ix,i)\ar[rr]^{\textstyle 1}&\ar@{}[ul]_(0.4){\textstyle \cong} &(F_ix,i), }
  $$
and
$$\xymatrix@C=40pt@R=45pt{
  (F_ix,i) \ar@/^2.5pc/[d]^{\textstyle 1}_{\textstyle \overset{\textstyle (\boldsymbol{l},1_i)}\cong}
  \ar[d]_{\textstyle (1_{1_i^*F_ix}\circ \iota_iF_ix,1_i)}\ar[rr]^{\textstyle 1}&&(F_ix,i)
  \ar[d]^{\textstyle (\iota_iF_ix\circ 1_{F_ix},1_i)}\ar@/_2.5pc/[d]_{\textstyle 1}^{\textstyle \overset{\textstyle (\boldsymbol{r},1_i)}\cong}\\(F_ix,i)\ar[rr]^{\textstyle 1}&\ar@{}[u]|{\textstyle \cong} &(F_ix,i).
  }
  $$

\vspace{0.2cm}

\subsection{The bicategorical Grothendieck construction as a tricolimit}
In \cite[\S \, 8]{gray}, Gray proved that the functor $\int_I:\cat^{ I^{\mathrm{op}}}\to \cat$ carries any lax $I$-diagram of categories to its lax colimit (or 2-colimit) in the 2-category of small categories. Next we shall prove the parallel fact for lax diagrams of bicategories. To do that, fix $I$ any small category and let $$\ct :\bicat\to \bicat^{ I^{\mathrm{op}}}$$
denote the diagonal trihomomorphism mapping any bicategory $\b$ to the constant lax $I$-diagram $\ct(\b):I^{\mathrm{op}}\to \bicat$ that $\b$ canonically defines. Then, we have the following theorem, whose proof this subsection is dedicated to.

\begin{theorem}\label{tricoli} The trihomomorphism $\int_I:\bicat^{ I^{\mathrm{op}}}\to \bicat$ is left triadjoint to the trihomomorphism $\ct :\bicat\to \bicat^{ I^{\mathrm{op}}}$\!.
\end{theorem}
\begin{proof} Remark first that a trihomomorphism $L:\T\to\T'$, where $\T$ and $\T'$ are tricategories, is called a {\em left triadjoint} for a trihomomorphism $R:\T'\to \T$, and $R$ is called a {\em right triadjoint} for $L$, if there is a biequivalence \cite[Definition 3.5]{g-p-s}
$\T'(L(-),-)\Rightarrow \T(-,R(-))$ in the tricategory  of trihomomorphisms with domain $\T^{\mathrm{op}}\times \T'$ and codomain $\bicat$, $\mathbf{Tricat}(\T^{\mathrm{op}}\times \T',\bicat)$, whose 1-cells are  tritransformations, whose 2-cells are trimodifications, and whose 3-cells are perturbations \cite[3.3]{g-p-s}.

Hence, we must prove that there is a tritransformation
$$\xymatrix{\bicat(\int_I(-),-)\Rightarrow \bicat^{ I^{\mathrm{op}}}\!(-,\ct(-))}$$
such that, for any lax diagram of bicategories $\f:I^{\mathrm{op}}\to \bicat$ and bicategory $\b$, the associated homomorphism
$$\xymatrix{\bicat(\int_I\f,\b)\rightarrow \bicat^{ I^{\mathrm{op}}}\!(\f,\ct(\b))}$$
is a biequivalence of bicategories. In more elementary terms, we shall prove the existence of tritransformations (the {\em unit} and {\em counit})
$$\xymatrix{\eta:1_{\bicat^{ I^{\mathrm{op}}}}\Longrightarrow \ct\int_I},\hspace{0.3cm}
\xymatrix{\epsilon: \int_I\ct \Longrightarrow 1_{\bicat}},$$
and equivalences (the {\em triangulators})
\begin{equation}\label{trian}
\xymatrix@C=35pt@R=35pt{\int_I\ar@{=>}[r]^{\textstyle \int_I\eta}\ar@{=>}[rd]_{\textstyle 1}&\int_I\ct\int_I\ar@{=>}[d]^{\textstyle \epsilon \int_I}\\
\ar@{}[ru]|(0.7){\overset{\textstyle T} \Lleftarrow}&\int_I\,,} \hspace{0.6cm}
\xymatrix@C=35pt@R=35pt{\ct\ar@{=>}[r]^{\textstyle \eta\,\ct}\ar@{=>}[rd]_{\textstyle 1}&\ct\int_I\ct\ar@{=>}[d]^{\textstyle \ct\,\epsilon }\\
\ar@{}[ru]|(0.7){\overset{\textstyle S} \Rrightarrow}&\ct\,,}
\end{equation}
that is, trimodifications $T$ and $S$ as above, such that, for any  $\f:I^{\mathrm{op}}\to \bicat$, the pseudo transformation $$\xymatrix{T\f:\epsilon\!\int_I\!\f\,\int_I\!\eta\f\Rightarrow 1_{\int_I\!\f}}$$
is a pseudo equivalence, and, for any bicategory $\b$, the corresponding 
$$\xymatrix{S\b:1_{\ct(\b)}\Rightarrow \ct(\epsilon\b)\,\eta\ct(\b)}$$
is a pseudo $I$-equivalence.

The proof is then divided into three parts.

\noindent \underline{Part I}. Here we exhibit the unit tritransformation $\xymatrix{\eta:1_{\bicat^{ I^{\mathrm{op}}}}\Rightarrow \ct\int_I}.$

At any lax $I$-diagram of bicategories $\f:I^{\mathrm{op}}\to \bicat$, the lax $I$-homomorphism $$\xymatrix{\eta=\eta\f:\f\to\ct(\int_I\!\f)}$$ works as follows:
For any object $i$ of $I$, $\eta_i:\f_i\to \int_I\!\f$ is the embedding homomorphism defined by
\begin{equation}\label{eeta}\xymatrix@C=3pt{y \ar@/^0.7pc/[rr]^{\textstyle u} \ar@/_0.7pc/[rr]_{\textstyle u'} & {}_{\textstyle \Downarrow\phi} &x}\ \overset{\textstyle \eta_i}\mapsto  \
\xymatrix@C=3pt{(y,i)  \ar@/^1.2pc/[rr]^{\textstyle (\iota_ix\!\circ\! u,1_i)} \ar@/_1.2pc/[rr]_{\textstyle (\iota_ix\!\circ\! u',1_i)} & {\textstyle \Downarrow\!(1_{\iota_ix}\!\circ\! \phi,1_i)} &(x,i) ,}\end{equation}
where, for the horizontal composition of
$1$-cells $z\overset{v}\to y \overset{u}\to x$ in $\f_i$, the invertible  structure $2$-cell $\eta_{i}(u)\circ \eta_{i}(v)\cong \eta_{i}(u\circ v)$ is provided by pasting in $\f_i$ the diagram
\begin{equation}\label{ji}
\xymatrix@R=30pt{&&\ar@{}[d]|(.7){\textstyle \cong}&&\\&1_i^*y\ar[r]^{1_i^*u} \ar@/^2pc/[rr]^{\textstyle 1_i^*(\iota\circ u)} &1_i^*x\ar[r]^{1_i^*\iota }\ar[rd]|{\ 1\ }\ar@{}@<1ex>[rd]^(.6){\textstyle \overset{\textstyle \delta}\cong}\ar@{}@<-1ex>[rd]_(.5){\textstyle\cong}&1_i^*1_i^*x \ar[d]^{\textstyle \chi }\\
z\ar[r]_{\textstyle v}\ar[ru]^{\textstyle \iota \circ v}\ar@{}@<-1ex>[ru]_{\textstyle =} &y\ar[r]_{\textstyle u}\ar[u]_{\textstyle\iota }\ar@{}[ur]|{\textstyle \overset{\textstyle \widehat{\iota}_u}\cong}&x\ar[r]_{\textstyle \iota }\ar[u]_{\textstyle \iota }&1_i^*x,
}
\end{equation}
and, for any object $x$ in $\f_i$, the structure isomorphism $1_{\eta_{i}x}\cong \eta_{i}(1_x)$ is the one given  by the canonical isomorphisms $\iota_ix\cong \iota_ix \circ 1_x$.
If $a:j\to i$ is any morphism in $I$, then the component at an object $x\in \f_i$ of the attached pseudo transformation
$$\xymatrix{\theta:\eta_j a^*\Rightarrow \eta_i:\f_i\to\int_I\!\f}$$
is the morphism $$\theta x=(1_{a^*x},a):(a^*x,j)\to (x,i),$$
 and, for each morphism $u:y\to x$ in $\f_i$, the invertible 2-cell $$\widehat{\theta}_u:\theta x\circ \eta_ja^*\!u\cong \eta_iu\circ \theta y$$ is that obtained by pasting the diagram
\begin{equation}\label{205}\xymatrix@R=15pt{
&&1_j^*a^*x\ar@{}[rd]|{\textstyle \underset{\textstyle \cong}{\widehat{\iota}}}\ar[rr]^{\textstyle 1_j^*(1_{a^*x})}&
&1_j^*a^*x\ar[rd]^{\textstyle \chi}\ar@{}[d]|(.6){\textstyle \underset{\textstyle \cong}{\gamma}}&\\
&a^*x\ar@{}[rd]|{\textstyle \cong}\ar[ru]^{\textstyle \iota a^*}\ar[rr]^{1_{a^*x}}&& a^*x\ar[rd]|{a^*\!\iota
}\ar@{}[d]|(.6){\textstyle  \cong}
\ar[ru]^{\iota a^*}\ar[rr]|{1_{a^*\!x}}&\ar@{}[d]|(.45){\textstyle \underset{\textstyle \cong}{\delta}}&a^*x\\
a^*y\ar[ru]^{\textstyle a^*u}\ar[rr]_{\textstyle 1_{a^*y}}&&a^*y\ar[ru]^{a^*u}\ar[rr]_{\textstyle a^*(\iota \circ u)}&
&a^*1_j^*x\,.\ar[ru]_{\textstyle \chi}&}
\end{equation}
For $k\overset{b}\to j\overset{a}\to i$, two composable morphisms of $I$, and any object $i$,  the invertible modifications
$$\xymatrix{\eta_kb^*a^*\ar@{=>}[r]^{\textstyle \eta_k\chi}\ar@{=>}[d]_{\textstyle \theta_{\!b}a^*}\ar@{}[rd]|{\textstyle \overset{\textstyle\Pi}\cong}&\eta_k(ab)^*\ar@{=>}[d]^{\textstyle \theta_{\!ab}}\\\eta_ja^*\ar@{=>}[r]_{\textstyle \theta_{\!a}}&\eta_i}
\hspace{0.6cm}
\xymatrix@R=8pt@C=14pt{& \eta_i\ar@2{->}[rdd]^{\textstyle 1_{\eta_i}}\ar@2{->}[ldd]_{\textstyle \eta_i\iota}\ar@{}[dd]|(.6){\textstyle{\overset{\textstyle \Gamma}{\textstyle\cong}}} &\\ & & \\ \eta_i1_i^*
\ar@2{->}[rr]_{\textstyle \theta_{1_i}} && \eta_i\,, }
$$
are respectively provided, at each object $x$ of $\f_i$, by pasting the diagrams
\begin{equation}\label{206}
\xymatrix@C=20pt@R=8pt{b^*a^*x\ar[r]^{\textstyle b^*1_{a^*}}&b^*a^*x\ar[r]^{\textstyle \chi}&(ab)^*x&\\
\ar@{}[rr]|(.35){\textstyle \cong}&&\ar@{}[r]|(.3){\textstyle \cong}&1_k^*(ab)^*x\ar[lu]_{\textstyle \chi}
\\
b^*a^*x\ar[uu]_{\textstyle 1_{b^*\!a^*}}\ar[r]_{\textstyle \chi}&(ab)^*x\ar[ruu]^{1}\ar@<-2ex>@{}[ruu]|(0.6){\textstyle \overset{\textstyle \gamma}\cong}\ar[r]_{\textstyle \iota(ab)^*}&1_k^*(ab)^*x\ar[ru]_{\textstyle \ 1_k^*\!(1_{(ab)^*x})\,,}\ar[uu]_\chi&}
\xymatrix@C=35pt@R=12pt{ 1_i^*x\ar@<0.5ex>@{}[rrdd]^(.45){\textstyle \cong\!\! \gamma}\ar[rrdd]_{1}\ar[r]^{\textstyle \iota 1_i^*}&1_i^*1_i^*x\ar[r]^{\textstyle 1_i^*(1_{1_i^*x})}\ar[rdd]^(.3){\chi}&1_i^*1_i^*x\ar[dd]^{\textstyle \chi}\\
&\ar@{}[ru]_(.6){\textstyle \cong}&\\
x\ar@{}[ru]|(.6){\textstyle \cong}\ar[uu]_{\textstyle \iota}\ar[rr]_{\textstyle \iota}&&1_i^*x.
}
\end{equation}

 If $F:\f\to\g$ is any lax $I$-homomorphism,  then the attached pseudo $I$-equivalence $$\xymatrix{\widehat{\eta}=\widehat{\eta}_F:\ct(\int_I\!F)\,\,\eta\f\Rightarrow \eta\g\,\,F,}$$ is, at any object  $i$ of $I$, the identity on objects pseudo equivalence $$\xymatrix{\widehat{\eta}_i: \int_I\!F\,\,\eta_i\Rightarrow \eta_i\,\,F_i:\f_i\to\int_I\!\g,}$$  that is, with $\widehat{\eta}_ix=1_{(F_ix,i)}$, and whose component at a morphism $u:y\to x$ of the bicategory $\f_i$ is canonically obtained from pasting
\begin{equation}\label{207}
\xymatrix{&F_i1_i^*x \ar[rd]^{\textstyle \theta}&\\F_iy\ar@<1ex>[ru]^{\textstyle F_i(\iota\circ u)}\ar[r]_{\textstyle F_iu}&F_ix\ar[u]|{F_i\iota}\ar[r]_{\textstyle \iota F_i}\ar@{}[ru]|(.35){\textstyle \cong\!\Gamma}\ar@{}[lu]|(.35){\textstyle \cong}&1_i^*F_ix;}
\end{equation}
and, for $a:j\to i$, the corresponding invertible modification
$$\xymatrix@R=8pt@C=25pt{\int_I\!F\,\eta_j\,a^*\ar@{=>}[dd]_{\textstyle \int_I\!F\,\theta}\ar@{=>}[r]^{\textstyle \widehat{\eta}_ja^*}&\eta_j\,F_j\,a^*\ar@{=>}[dr]^{\textstyle \eta_j\theta} &\\
\ar@{}[rr]|(0.4){\textstyle \cong} && \eta_j\,a^*\,F_i\ar@{=>}[ld]^{\textstyle \theta F_i}\\
\int_I\!F\,\eta_i\ar@2{->}[r]_{\textstyle \widehat{\eta}_i} & \eta_i\,F_i&}
$$
is, at any object $x$ of $\f_i$, that canonically obtained from pasting in $\g_j$
\begin{equation}\label{208}
\xymatrix@R=40pt@C=55pt{
F_ja^*x\ar@<-1ex>@{}[rrd]|(.3){\textstyle\cong}\ar[r]^{\textstyle F_j1_a^*}\ar[d]_{\textstyle \theta}&F_ja^*x\ar@<-1ex>@{}[rd]^(.7){\textstyle\cong}\ar[r]^{\textstyle \theta}&a^*F_ix\\
a^*F_ix\ar@<-1ex>@{}[rru]_{\textstyle\cong\!\gamma}\ar[rru]^{1}\ar[r]_{\textstyle \iota a^*F_i}&1_j^*a^*F_ix\ar[ru]_\chi\ar[r]_{\textstyle 1_j^*(1_{a^*F_ix})}&1_j^*a^*F_ix\ar[u]_{\textstyle \chi}.
}
\end{equation}

The tritransformation $\eta$ takes any pseudo $I$-transformation $m:F\Rightarrow F':\f\to\g$ to the invertible $I$-modification
$$\xymatrix{\ct(\int_I\!m)\,\eta\f \circ \widehat{\eta}_F\cong \widehat{\eta}_{F'}\circ \eta\g\, m,}$$
whose component at an object $x$ of $\f_i$, for any $i$ of $I$, is the canonical isomorphism $1_{\eta_ix}\circ \eta_i(m_ix)\cong \eta_i(m_ix)\circ 1_{\eta_ix}$ of the bicategory   $\int_I\!\g$.

If $\f\overset{F}\to\g\overset{G}\to\h$ are any two composable lax $I$-homomorphisms, then the estructure invertible $I$-modification for the tritransformation $\eta$
$$
\xymatrix@C=50pt@R=35pt{\ar@{}[rd]|{\textstyle \cong}\ct(\int_I\!G)\,\ct(\int_I\!F)\,\eta
\ar@{=>}[r]^{\textstyle \ct(\int_I\!G)\widehat{\eta}}
\ar@{=>}[d]_{\textstyle \ct(\Sigma)\eta}&\ct(\int_I\!G)\,\eta\,F\ar@{=>}[d]^{\textstyle \widehat{\eta}F}\\
\ct(\int_I\!(GF))\,\eta\ar@{=>}[r]^{\textstyle \widehat{\eta}}&\eta\,\, G\,F}
$$
is, for any objects $i$ of $I$ and $x$ of $\f_i$, the canonical isomorphism in the bicategory   $\int_I\!\h$ $$\xymatrix{1_{(G_iF_ix,i)}\circ \int_I\!G(1_{(F_ix,i)})\cong 1_{(G_iF_ix,i)}\circ 1_{(G_iF_ix,i)}}.$$
And, finally, say that for any lax $I$-diagram of bicategories $\f$, the equality $$\xymatrix{\ct(\Sigma_\f)\,\eta\f=\widehat{\eta}_{1_\f} :\ct(\int_I\!1_\f)\,\eta\f\Rightarrow \eta\f}$$
holds, and the corresponding $I$-modification of $\eta$ attached to $\f$ between them is the identity one. This makes complete the description of the tritransformation $\eta$.

\vspace{0.2cm}
\noindent \underline{Part II}. Here we shall describe the counit tritransformation $\xymatrix{\epsilon: \int_I\ct \Rightarrow 1_{\bicat}},$
which is easier to describe than  the unit, since the composite $\int_I\ct $ can be identified with the trihomomorphism $(- )\times I:\bicat\to\bicat$,  and then $\epsilon$ with the projection on the first factor. More precisely, for any bicategory $\b$, $$\xymatrix{\epsilon=\epsilon\b: \int_I\!\ct(\b) \to \b}$$ is the normal homomorphism
$$
 \xymatrix@C=3pt@R=3pt{ (y,j) \ar@/^1pc/[rr]^{\textstyle (u,a)} \ar@/_1pc/[rr]_{\textstyle (u',a)} &
\Downarrow(\alpha,a) & (x,i)}~~~ \overset{\textstyle \epsilon}\mapsto ~~~\xymatrix@C=2pt@R=2pt{
y\ar@/^0.7pc/[rr]^{\textstyle u} \ar@/_0.7pc/[rr]_{\textstyle u'} & \Downarrow
\alpha& x,}
$$
whose structure constraints for horizontal compositions of 1-cells are given by the left identity constraints of the
bicategory $\b$. For any  two bicategories $\b$ and $\c$, the diagram
$$
\xymatrix@C=10pt@R=12pt{\bicat(\b,\c)\ar[rr]^{\textstyle \int_I\ct}\ar[dr]_{\textstyle \epsilon^*}&&\bicat(\int_I\ct(\b),\int_I\ct(\c))\ar[ld]^{\textstyle \epsilon_*}
\\ &\hspace{0.3cm}\bicat(\int_I\ct(\b),\c)&
}$$
commutes, and  the corresponding pseudo natural equivalence $$\xymatrix{\widehat{\epsilon}:\epsilon_*\int_I\ct \Rightarrow \epsilon^*}$$ is the identity.

For any two homomorphism $\b\overset{F}\to\c\overset{G}\to\d$, the invertible modification
$$
\xymatrix@C=50pt{\ar@{}[rd]|{\textstyle \cong}\epsilon\int_I\ct(G)\,\int_I\ct(F)\ar@{=>}[r]^-{\textstyle \widehat{\epsilon}\,\int_I\ct(F)}\ar@{=>}[d]_{\textstyle \epsilon\Sigma}&G\,\epsilon\int_I\ct(F)\ar@{=>}[d]^{\textstyle G\widehat{\epsilon}}\\ \epsilon\int_I\ct(GF)\ar@{=>}[r]^{\textstyle \widehat{\epsilon}}&GF\epsilon
}
$$
is, at any object $x$ of $\b$, the canonical isomorphism $G1_{Fx}\circ 1_{GFx}\cong 1_{GFx}\circ 1_{GFx}$ in $\d$, and, for any $\b$, we have $\epsilon\b\,\Sigma_{\ct(\b)}=\widehat{\epsilon}_{1_\b}$, and the corresponding invertible modification for $\epsilon$ at $\b$ is the identity.

\vspace{0.2cm}
\noindent \underline{Part III}. We conclude here the proof by showing the triangulators  $T$ and $S$ in (\ref{trian}).

   The component of $T$ at any lax $I$-diagram $\f:I^{\mathrm{op}}\to \bicat$, is the pseudo equivalence $T\f:\epsilon\!\int_I\!\f\, \int_I\!\eta\f\Rightarrow 1_{\int_I\!\f}$ with $T\f(x,i)=1_{(x,i)}$ for any object $(x,i)$ of $\int_I\!\f$, and whose component at a morphism $(u,a):(x,i)\to (y,j)$ is canonically provided by the 2-cell in $\f_i$ pasted of
$$\xymatrix{a^*x\ar[rd]|{\ \iota a^*\ }\ar[r]^{\textstyle \iota a^*}\ar[d]_{\textstyle
1_{a^*}}&1_i^*a^*x\ar[d]^{\textstyle 1_i^*(1_{a^*x})}\\a^*x\ar@{}[ru]|(.3){\textstyle \cong\!\gamma}|(.7){\textstyle
\cong}&1_i^*a^*x. \ar[l]^{\textstyle \chi}}
$$
For any lax $I$-homomorphism $F:\f\to\g$, the structure invertible modification
$$
\xymatrix@C=70pt{\ar@{}[rd]|{\textstyle \cong}\epsilon\!\int_I\!\!\g\,\,\ct(\int_IF)\int_I\!\eta\f\ar@{=>}[r]^{\textstyle \epsilon\!\int_I\!\g\,\int_I\!\widehat{\eta}_\f}\ar@{=>}[d]_{\textstyle\widehat{\epsilon}_F\int_I\!\eta\f}&\epsilon\!\int_I\!\g\,\int_I\!\eta\g\,\int_I\!F\ar@{=>}[d]^{\textstyle T\g\,\int_I\!F}\\\int_I\!F\,\,\epsilon\!\int_I\!\f\,\int_I\!\eta\f\ar@{=>}[r]_{\textstyle \int_I\!F\,\,T\f}&\int_I\!F}
$$
is, at any object $(x,i)$ of $\int_I\!\f$, the canonical isomorphism $\xymatrix{1_{\int_I\!F(x,i)} \cong \int_I\!F(1_{(x,i)})}$.

And when it comes to S, say that,
for any bicategory $\b$, the component of ${S\b\!:\!1_{\ct(\b)}\Longrightarrow \ct(\epsilon\b)\,\eta\ct(\b)}$ at an object $i$ of $I$ is the pseudo equivalence which is the identity on objects of $\b$, and whose component at a morphism $u:y\to x$ is the canonical isomorphism $1_x\circ (1_x\circ u)\cong u\circ 1_y$. For any homomorphism $F:\b\to\c$, the structure invertible modification
$$
\xymatrix@C=70pt{\ar@{=>}[d]_{\textstyle S\c\,\,\ct(F)}\ct(F)\ar@{}[rd]|-{\textstyle\cong}\ar@{=>}[r]^-{\textstyle \ct(F)\,\,S\b}&\ct(F)\,\,\ct(\epsilon\b)\,\,\eta\ct(\b)\ar@{=>}[d]^{\textstyle \ct(\widehat{\epsilon}_\b)\,\,\eta\ct(\b)}\\\ct(\epsilon\c)\,\,\eta\ct(\c)\,\,\ct(F)\ar@{=>}[r]_-{\textstyle \ct(\epsilon\c)\,\,\widehat{\eta}_{\ct(F)}}&\ct(\epsilon\c)\,\,\ct(\int_I\!\ct(F))\,\,\eta\ct(\b),
}
$$
at any $i\in\mbox{Ob}I$ and $x\in\mbox{Ob}\b$, is the canonical 2-cell   $1_{Fx}\circ F(1_x)\cong 1_{Fx}\circ 1_{Fx}$ in the bicategory $\c$.
\end{proof}

\section{Rectification}\label{rectification} Following Giraud \cite{giraud}, Street \cite{street72}, Thomanson \cite{thomason}, and May \cite{may}, we shall show here how any lax $I$-diagram of bicategories
$\f=(\f,\chi,\iota, \omega,\gamma,\delta): I^\mathrm{op}\to \bicat$ has, naturally associated to it, a genuine $I$-diagram of
bicategories, that is, a functor \begin{equation}\label{recti}\f^{^\mathrm{r}}:I^\mathrm{op}\to \mathbf{Hom}\subset  \bicat \end{equation} that, as we will prove later, represents
the same homotopy type as $\f$. This $I$-diagram of bicategories $\f^\mathrm{r}$ is built as follows.
For each object $i$ of $ I$, let $i/I$ be the comma category whose objects are the arrows in $I$ of the form $b:i\to k$
and whose morphisms are the appropriate commutative triangles. By composing $\f$ with the obvious forgetful functor
$i/I\overset{\pi_i}\to I$, we obtain  the  lax $(i/I)$-diagram
$$\f\pi_i:(i/I)^\mathrm{op}\to \bicat,$$ and then, by the Grothendieck construction, a new bicategory
$$\xymatrix{\f^{^\mathrm{r}}_i=\int_{i/I}\f\pi_i, }$$
whose set of objects is $\bigsqcup\limits_{i\overset{b}\to k}\!\!\mbox{Ob}\f_k$  and hom-categories
$$\f^{^\mathrm{r}}_i\big((y,i\overset{c}\to l), (x,i\overset{ b}\to k)\big)=\bigsqcup_{\scriptsize \begin{array}{c}l\overset{\textstyle d}\to k\\ dc=b\end{array}}\!\! \f_l(y,d^*x).$$

An arrow $a:j\to i$ in $I$ induces a functor $a^*:i/I\to j/I$ with $\pi_j a^*=\pi_i$ and hence a strict 2-functor
$
a^*:\f^{^\mathrm{r}}_i\to \f^{^\mathrm{r}}_j,
$
$$ \xymatrix@C=7pt@R=7pt{ (y,i\overset{c}\to l) \ar@/^1.1pc/[rr]^{\textstyle (u,d)} \ar@/_1.1pc/[rr]_{\textstyle (u',d)} &
\Downarrow(\alpha,d) & (x,i\overset{b}\to k)}~~~ \overset{\textstyle a^*}\mapsto ~~~\xymatrix@C=7pt@R=7pt{
(y,j\overset{ca}\to l) \ar@/^1.1pc/[rr]^{\textstyle (u,d)} \ar@/_1.1pc/[rr]_{\textstyle (u',d)} & \Downarrow(\alpha,d) &
(x,j\overset{ba}\to k).}
$$

For $i$ any object of $I$, we have $1_i^*=1_{\f^{^\mathrm{r}}_i}$, and for $k\overset{b} \to j \overset{a} \to i$, any
two composible arrows of $I$, the equality $b^*a^*=(ab)^*:\f^{^\mathrm{r}}_i \to \f^{^\mathrm{r}}_k$ holds. Therefore,
we have defined a genuine $I$-diagram of bicategories and strict functors
$\f^{^\mathrm{r}}:I^{\mathrm{op}}\to \mathbf{Hom}\subset \bicat$, which we refer to as the {\em rectification} of $\f$.

\begin{proposition} The assignment $\f\mapsto \f^{^\mathrm{r}}$ is the function on objects of a triendomorphism $(\
)^{^\mathrm{r}} :\bicat^{ I^{\mathrm{op}}}\to \bicat^{ I^{\mathrm{op}}}$, which we call {\em rectification}.
\end{proposition}
\begin{proof}
If $F:\f\to \g$ is any given lax $I$-homomorphism between lax $I$-diagrams
$\f,\g:I^{\mathrm{op}}\to \bicat$, then, for each object $i$ of $I$, the composite $F\pi_i:\f\pi_i\to \g\pi_i$ is a lax
$(i/I)$-homomorphism inducing a homomorphism
$$
\xymatrix{F^{^\mathrm{r}}_i=\int_{i/I}\!F\pi_i:\f^{^\mathrm{r}}_i \to \g^{^\mathrm{r}}_i.}
$$
The assignment $i\mapsto F^{^\mathrm{r}}_i$ completely determines an $I$-homomorphism
$
F^{^\mathrm{r}}:\f^{^\mathrm{r}}\to \g^{^\mathrm{r}},
$
that is, a lax $I$-homomorphism such that, for any arrow $a:j\to i$ in $I$, the equality
$
F^{^\mathrm{r}}_j a^* = a^* F^{^\mathrm{r}}_i
$
holds, the pseudo-transformations $\theta$ for $F^{^\mathrm{r}}$ are identities, and the invertible
modifications $\Pi$ and $\Gamma$ are given by the unit constraints. Call
$F^{^\mathrm{r}}$ {\em the rectification of $F$}.

Similarly, for $m:F\Rightarrow G:\f\to \g$ a pseudo $I$-transformation, we define its {\em
rectification}
$
m^{{_\mathrm{r}}}:F^{^\mathrm{r}}\Rightarrow G^{^\mathrm{r}}:\f^{^\mathrm{r}}\to \g^{^\mathrm{r}} $ to be the $I$-transformation given by
writing
$$\xymatrix{ m^{{_\mathrm{r}}}_i=\int_{i/I}\! m\pi_i:F^{^\mathrm{r}}_i\Rightarrow G^{^\mathrm{r}}_i,} $$
for each object $i$ of $I$. For any arrow $a:j\to i$ in $I$, the equality
$ a^*m^{{_\mathrm{r}}}_i =m^{{_\mathrm{r}}}_j a^* $ holds, so that $m^{{_\mathrm{r}}}$ is a genuine $I$-transformation in
the sense that the corresponding invertible modification $\mathrm{M}$ is that given by the unit
constraints. And finally, for $\sigma:m \Rrightarrow n$ an $I$-modification, we take
$
\sigma^{{_\mathrm{r}}}:m^{{_\mathrm{r}}}\Rrightarrow n^{{_\mathrm{r}}}
$
to be the $I$-modification defined at any object $i$ of $I$ by
$$\xymatrix{  \sigma^{{_\mathrm{r}}}_i=\int_{i/I}\!\sigma\pi_i: m^{{_\mathrm{r}}}_i\Rrightarrow n^{{_\mathrm{r}}}_i.}$$

The rectification constructions above actually lead to a triendomorphism of the tricategory of lax $I$-diagrams, simply thanks to the Grothendieck construction
$\xymatrix{\int_I \!:\!\bicat^{ I^{\mathrm{op}}}\!\rightarrow \bicat}$ being a trihomomorphism. Thus, the structure isomorphisms of the rectification homomorphism $$(\ )^{^\mathrm{r}}\!:\!\xymatrix{\bicat^{ I^{\mathrm{op}}}\!(\f,\g)\to \bicat^{ I^{\mathrm{op}}}\!(\f^{^\mathrm{r}},\g^{^\mathrm{r}})  }$$
are as follows: for any  pseudo $I$-transformations $F\overset{m}\Rightarrow G \overset{n}\Rightarrow H:\f\to \g$, the structure
invertible $I$-modification $n^{{_\mathrm{r}}}\circ m^{{_\mathrm{r}}}\Rrightarrow (n\circ m)^{{_\mathrm{r}}}$ at an object $i$ of $I$ is
$$n^{{_\mathrm{r}}}_i\circ m^{{_\mathrm{r}}}_i=\xymatrix{ \int_{i/I}\!n\pi_i\circ \int_{i/I}\!m\pi_i ~~\overset{(\ref{gc1})}\Rrightarrow ~~\int_{i/I}\!(n\circ m)\pi_i}=(n\circ m)^{{_\mathrm{r}}}_i,$$ while the invertible $I$-modification $1_{F^{^\mathrm{r}}}\Rrightarrow 1_F^{^\mathrm{r}}$ at an object $i$ is
$$(1_{F^{^\mathrm{r}}})_i=\xymatrix{1_{\int_{i/I}\!F\pi_i} ~~ \overset{(\ref{gc2})}\Rrightarrow ~~ \int_{i/I}\!1_F\pi_i}=(1_F^{^\mathrm{r}})_i.$$
Furthermore, for any pair of lax
$I$-homomorphisms $\f\overset{F}\to \g \overset{G}\to \h$, the components, at an object $i\in \mbox{Ob}I$, of  the pseudo-natural $I$-equivalences
$\xymatrix{
\Sigma^{^\mathrm{r}}_{G,F}:G^{^\mathrm{r}}F^{^\mathrm{r}} \Rightarrow (GF)^{^\mathrm{r}}
}$ and $\Sigma^{^\mathrm{r}}_\f:1_\f^{^\mathrm{r}}\Rightarrow 1_{\f^{^\mathrm{r}}}$ are, respectively,

$$\begin{array}{rcl}
G^{^\mathrm{r}}_iF^{^\mathrm{r}}_i=\xymatrix{\int_{i/I}\!G\pi_i\,\int_{i/I}\!F\pi_i}&
\overset{(\ref{gc3})}\Rightarrow& \xymatrix{\int_{i/I}GF\pi_i}=(GF)^{^\mathrm{r}}_i,\\[5pt]
(1^{^\mathrm{r}}_\f)_i=\xymatrix{\int_{i/I}\!1_{\f\pi_i}}&
\overset{(\ref{gc4})}\Rightarrow&\xymatrix{1_{\int_{i/I}\!\f\pi_i}}=(1_{\f^{^\mathrm{r}}})_i,\end{array}
$$
that is, $(\Sigma^{^\mathrm{r}}_{G,F})_i=\Sigma_{G\pi_i,F\pi_i}$ and $(\Sigma^{^\mathrm{r}}_\f)_i=\Sigma_{\f\pi_i}$. For any morphism $a:j\to i$, the equalities $\Sigma^{^\mathrm{r}}_ja^*=a^*\Sigma^{^\mathrm{r}}_i$ hold, and the components $\mathrm{M}_a$, both for $\Sigma^{^\mathrm{r}}_{G,F}$ and $\Sigma^{^\mathrm{r}}_\f$,  are the canonical modifications given by the identity constraints.

\noindent Given lax $I$-homomorphisms $\f\overset{F} \to \g \overset{G}\to \h \overset{H}\to \mathcal{K}$, the structure invertible $I$-modifications $\omega^{^\mathrm{r}}$, $\delta^{^\mathrm{r}}$ and $\gamma^{^\mathrm{r}}$ for $(~)^{^\mathrm{r}}$, as in the definition of a trihomomorphism,

$$
 \xymatrix@C=-17pt@R=15pt{ &&(HG)^{^\mathrm{r}}F^{^\mathrm{r}} \ar@2{->}[rrd]^{\textstyle \Sigma^{^\mathrm{r}}}&&\\
 H^{^\mathrm{r}}G^{^\mathrm{r}}F^{^\mathrm{r}} \ar@{}[rrr]|(0.7){\textstyle \overset{\textstyle \omega^{^\mathrm{r}}}\cong}\ar@2{->}[urr]^{\textstyle \Sigma^{^\mathrm{r}} \,1} \ar@2{->}[dr]_{\textstyle 1\,\Sigma^{^\mathrm{r}}}&&&& ((HG)F)^{^\mathrm{r}} \\
 &H^{^\mathrm{r}}(GF)^{^\mathrm{r}} \ar@2{->}[rr]_{\textstyle \Sigma^{^\mathrm{r}}}&&(H(GF))^{^\mathrm{r}} \ar@{=>}[ur]_(0.6){\textstyle \boldsymbol{a}^{^\mathrm{r}}}&}
  \hspace{-0.2cm}
 \xymatrix@C=8pt@R=15pt{F^{^\mathrm{r}}1_\f^{^\mathrm{r}} \ar@{=>}[ddrr]_{\textstyle 1\Sigma_\f^{^\mathrm{r}}}\ar@{=>}[rr]^{\textstyle \Sigma_{F,1_\f}^{^\mathrm{r}}}&&
 F^{^\mathrm{r}}1_\f^{^\mathrm{r}} \ar@{}[lldd]|(0.3){\textstyle \overset{\textstyle \delta^{^\mathrm{r}}}\cong}\ar@{=>}[dd]^{\textstyle \boldsymbol{r}^{^\mathrm{r}} }\\&&\\&&F^{^\mathrm{r}}}
  \hspace{0.1cm}
   \xymatrix@C=8pt@R=15pt{1_\g^{^\mathrm{r}}F^{^\mathrm{r}} \ar@{=>}[ddrr]_{\textstyle \Sigma_\g^{^\mathrm{r}} 1}\ar@{=>}[rr]^(0.45){\textstyle \Sigma_{1_\g,F}^{^\mathrm{r}}}&&
 (1_\g F)^{^\mathrm{r}}\ar@{}[ddll]|(0.3){\textstyle \overset{\textstyle \gamma^{^\mathrm{r}}}\cong}\ar@{=>}[dd]^{\textstyle \boldsymbol{l}^{^\mathrm{r}} }\\&&\\ &&F^{^\mathrm{r}}}
 $$
are, respectively, given by the families of modifications (\ref{odg}), $\omega_{_{H\pi_i, G\pi_i,F\pi_i}}$, $\delta_{_{F\pi_i}}$ and $\gamma_{_{F\pi_i}}$,
$i\in \mbox{Ob}I$.
\end{proof}

Every lax diagram of bicategories is related to its rectification by a canonical lax homomorphism, which we describe as follows:

\begin{lemma}\label{lemj}
Given $\f=(\f,\chi,\iota, \omega,\gamma,\delta):I^\mathrm{op}\to \bicat$, any lax $I$-diagram of bicategories, there is a lax $I$-homomorphism
$$J=(J,\theta,\Pi,\Gamma):\f\to \f^{^\mathrm{r}},$$
whose component at an object $i$ of $I$ is the homomorphism $J_i:\f_i\to \f^{^\mathrm{r}}_i$ acting by
\begin{equation}\label{defji}
\xymatrix@C=2pt{y \ar@/^0.6pc/[rr]^{\textstyle u} \ar@/_0.7pc/[rr]_{\textstyle u'} & {\textstyle \Downarrow\!\phi} &x}\ \overset{\textstyle J_i}\mapsto  \
\xymatrix@C=3pt{(y,1_i)  \ar@/^1.2pc/[rr]^{\textstyle (\iota_ix\!\circ\! u,1_i)} \ar@/_1.2pc/[rr]_{\textstyle (\iota_ix\!\circ\! u',1_i)} & {\textstyle \Downarrow\!(1_{\iota_ix}\!\circ\! \phi,1_i)} &(x,1_i) }.
\end{equation}
For the horizontal composition of
$1$-cells $z\overset{v}\to y \overset{u}\to x$ in $\f_i$, the structure invertible $2$-cell $J_i(u)\circ J_i(v)\cong J_i(u\circ v)$ is provided by pasting the diagram $(\ref{ji})$ in $\f_i$
and, for any object $x$ in $\f_i$, the structure isomorphism $1_{J_ix}\cong J_i(1_x)$ is that given  by the canonical isomorphisms $\iota_ix\cong \iota_ix \circ 1_x$.

If $F=(F,\theta,\Pi,\Gamma):\f\to \g$ is any lax $I$-homomorphism, then there is a pseudo $I$-equivalence
\begin{equation}\label{mj}
\xymatrix@R=15pt{\f\ar[r]^{\textstyle F}\ar[d]_{\textstyle J}&\g\ar[d]^{\textstyle J}\\ \f^{^\mathrm{r}}\ar[r]_{\textstyle F^{^\mathrm{r}}}\ar@{}[ru]|{\textstyle \overset{\textstyle m}\Rightarrow }& \g^{^\mathrm{r}}}
\end{equation}
\end{lemma}
\begin{proof}
Given a morphism $a:j\to i$ in the category $I$, the component of the pseudo-transformation \begin{equation}\label{ctheta}\theta:J_ja^*\Rightarrow a^*J_i\end{equation}
 at an object $x\in \text{Ob}\f_i$ is the morphism, in $\f^{^\mathrm{r}}_j$, $\theta x=(1_{a^*x}, a):(a^*x,1_j)\to (x,a)$. Moreover, for each morphism $u:y\to x$ in $\f_i$, the invertible $2$-cell $$\widehat{\theta}_u:\theta x\circ J_ja^*\!u\cong a^*\!J_iu\circ \theta y$$ is that obtained by pasting the diagram (\ref{205}).

For $k\overset{b}\to j\overset{a}\to i$, two composable morphisms of $I$, and any object $i$,  the invertible modifications

$$\xymatrix@R=8pt@C=25pt{J_kb^*a^*\ar@{=>}[dd]_{\textstyle \theta a^*}\ar@{=>}[r]^{\textstyle J_k\chi}&J_k(ab)^*\ar@{=>}[dr]^{\textstyle \theta} &\\
\ar@{}[rr]|(0.4){\textstyle\overset{\textstyle \Pi}\cong} && (ab)^*J_i\\
b^*J_ja^*\ar@2{->}[r]_{\textstyle b^*\theta} & b^*a^*J_i\ar@{=>}[ru]_{\textstyle 1 J_i}&}
\hspace{0.6cm}
\xymatrix@R=12pt@C=16pt{& J_i\ar@2{->}[rdd]^{\textstyle 1_{J_i}}\ar@2{->}[ldd]_{\textstyle J_i\iota}\ar@{}[dd]|(.6){\textstyle{\overset{\textstyle \Gamma}{\textstyle\cong}}} &\\ & & \\ J_i1_i^*
\ar@2{->}[rr]_{\textstyle \theta} && 1_i^*J_i\,, }
$$
are respectively provided, at each object $x$ of $\f_i$, by pasting the diagrams (\ref{206}).

Given $F:\f\to\g$, a lax $I$-homomorphism, for each object $i$ of $I$, the pseudo-natural equivalence (\ref{mj}) at $i$, $m_i: F^{^\mathrm{r}}_iJ_i\Rightarrow J_iF_i:\f_i\to\g^{^\mathrm{r}}_i$, is the identity on objects, that is, $m_ix=1_{(F_ix,1_i)}$, while its component at a 1-cell $u:y\to x$ of the bicategory $\f_i$ is canonically obtained from pasting (\ref{207}).

Finally, for $a:j\to i$ a morphism of $I$, the corresponding invertible modification
$$\xymatrix{F_j^{^\mathrm{r}}J_ja^*\ar@{=>}[r]^{\textstyle m_j a^*}\ar@{=>}[d]_{\textstyle \theta}\ar@{}[rd]|{\textstyle \underset{\cong}{\textstyle \mathrm{M}}}&J_jF_ja^*\ar@{=>}[d]^{\textstyle \theta}\\
a^*F_i^{^\mathrm{r}}J_i\ar@{=>}[r]_{\textstyle a^*m_i}&a^*J_iF_i}
$$
is that obtained from (\ref{208}).
\end{proof}

We should comment that the data in the previous lemma describes the components at objects and morphisms for a tritransformation $J:1_{\bicat^{ I^{\mathrm{op}}}}\Rightarrow (\ )^{^\mathrm{r}}$, whose full description is left to the reader. Furthermore, although for any given lax diagram $\f$, the lax $I$-homomorphism $J:\f\to \f^{^\mathrm{r}}$  does not have any right biadjoint (in the tricategory $\bicat^{ I^{\mathrm{op}}}$), we have the following:

\begin{lemma}\label{lemj2}  Let $\f=(\f,\chi,\iota, \omega,\gamma,\delta):I^\mathrm{op}\to \bicat$ be a lax $I$-diagram
of bicategories. For any object $i$ of the category $I$, the homomorphism in $(\ref{defji})$, $J_i:\f_i\to
\f^{^\mathrm{r}}_i$, has a right biadjoint.
\end{lemma}
\begin{proof}

 The right biadjoint to $J_i$ is the homomorphism $R_i: \f^{^\mathrm{r}}_i\to\f_i$ such that
$$
 \xymatrix@C=7pt@R=7pt{ (y,i\overset{b}\to k) \ar@/^1.1pc/[rr]^{\textstyle (u,d)} \ar@/_1.1pc/[rr]_{\textstyle (u',d)} &
\Downarrow(\alpha,d) & (x,i\overset{a}\to j)}~~~ \overset{\textstyle R_i}\mapsto ~~~\xymatrix@C=7pt@R=7pt{
b^*y\ar@/^1.1pc/[rr]^{\textstyle \chi\circ b^*u} \ar@/_1.1pc/[rr]_{\textstyle \chi \circ b^*u'} & \Downarrow
1_{\chi}\!\circ b^*\alpha & a^*x.}
$$
If $(z,i\overset{c}\to l)\overset{(v,e)}\longrightarrow (y,i\overset{b}\to k) \overset{
(u,d)}\longrightarrow (x,i\overset{a}\to j)$ are any two composible 1-cells of $ \f^{^\mathrm{r}}_i$, then the
structure invertible 2-cell $\xymatrix{R_i(u,d)\circ R_i(v,e)\cong R_i((u,a)\circ(v,b)\big)}$
is provided by pasting the diagram in $\f_i$
$$\xymatrix@C=65pt{c^*z\ar@{}[rd]|{\textstyle \cong}\ar[d]_{\textstyle c^*\!v}\ar[r]^(0.45){\textstyle c^*(\chi\circ(e^*\!u\circ v))}&c^*(de)^*x\ar[dr]^{\textstyle \chi} &\\
c^*e^*y\ar[d]_{\textstyle \chi}\ar[r]_{c^*\!e^*\!u}&c^*e^*d^*x\ar@{}[r]|
(0.4){\textstyle \overset{\textstyle \omega}\cong}\ar[d]^{\chi
d^*}\ar[u]_{c^*\chi}&a^*x \\
b^*y\ar@{}[ur]|{\textstyle \overset{\textstyle \widehat{\chi}_u}\cong }\ar[r]_{\textstyle b^*\!u}&b^*d^*x\ar[ru]_{\textstyle \chi}&
}$$
and, for each object $(x,i\overset{a}\to j)$, the  identity structure constraint $1_{R_i(x,a)}\cong R_i(1_{(x,a)})$
is $\delta_a:\chi\circ a^*\!\iota \Rightarrow 1_{a^*}$.

The unit of the biadjunction is the pseudo-transformation $\eta:1_{\f_i}\Rightarrow R_iJ_i$, with
$\eta x=\iota_ix:x\to 1_i^*x$, for each object $x$ of $\f_i$, and whose component at a 1-cell $u:y\to x$ is the
invertible 2-cell obtained by pasting
$$
\xymatrix@R=13pt{y\ar[r]^{\textstyle u}\ar[dd]_{\textstyle \iota}&x\ar@{}[rd]|(.3){\textstyle \cong}\ar[r]^{\textstyle \iota}\ar[d]_{\iota}&1_i^*x\\
\ar@{}[r]^{\textstyle\underset{\textstyle \cong}{\widehat{\iota}}}&1_i^*x\ar@{}[d]|(.6){\textstyle
\cong}\ar[rd]\ar@{}@<-1ex>[rd]^{1_i^*\iota}\ar[ru]|{\ 1\ } \ar@{}[r]|(.6){\textstyle \cong\!\delta}&
\\1_i^*y\ar[rr]_{\textstyle \,1_i^*(\iota\circ u)}\ar[ru]\ar@{}@<-1ex>[ru]^{1_i^*u}&&1_i^*1_i^*x,\ar[uu]_{\textstyle
\chi}}
$$
and the counit of the biadjunction is the pseudo-transformation $\epsilon: J_iR_i\Rightarrow 1_{\f^{^\mathrm{r}}_i}$,
with
$\epsilon(x,i\overset{a}\to j)=(1_{a^*x},a):(a^*x,1_i)\longrightarrow (x,a)$,
and whose component at a morphism $(u,c):(y,i\overset{b}\to k)\to (x,i\overset{a}\to j)$ in $\f^{^\mathrm{r}}_i$ is the
invertible deformation provided from pasting in the bicategory $\f_i$
$$
\xymatrix@R=35pt{ b^*y\ar[r]^{\textstyle b^*\!u}\ar[d]_{\textstyle 1_{b^*}}&b^*c^*x\ar@{}[d]|{\textstyle
\cong}\ar[r]^{\textstyle \chi}&a^*x\ar[r]^{\textstyle \iota a^*}
\ar[rd]|{\ \iota a^*\ }\ar[d]_{1_{a^*}}&1_i^*a^*x\ar[d]^{\textstyle 1_i^*(1_{a^*x})}\\
b^*y\ar[r]_{\textstyle b^*\!u}&b^*c^*x\ar[r]_{\textstyle \chi}&a^*x\ar@{}[ru]|(.3){\textstyle
\cong\!\gamma}|(.7){\textstyle \cong}&1_i^*a^*x.\ar[l]^{\textstyle \chi}}$$

The invertible modification triangulators $1_{R_i}\Rrightarrow R_i\epsilon\circ \eta R_i$ and $\epsilon J_i\circ J_i\eta \Rrightarrow 1_{J_i}$
are, at objects $(x,i\overset{a}\to j)$ of $\f^{^\mathrm{r}}_i$ and $x$ of $\f_i$, respectively obtained from  pasting
 the diagrams below in $\f_i$.
$$
\xymatrix{a^*x\ar[rd]|{\ \iota a^*\ }\ar[r]^{\textstyle \iota a^*}\ar[d]_{\textstyle
1_{a^*}}&1_i^*a^*x\ar[d]^{\textstyle 1_i^*(1_{a^*x})}\\a^*x\ar@{}[ru]|(.3){\textstyle \cong\!\gamma}|(.7){\textstyle
\cong}&1_i^*a^*x \ar[l]^{\textstyle \chi}}\hspace{0.6cm} \xymatrix{x\ar[r]^{\textstyle \iota}\ar[dr]_{\textstyle
\iota}&1_i^*x\ar@{}[rd]|(.3){\textstyle \cong\!\gamma}|(.7){\textstyle
\cong}\ar[d]\ar@<-0.5ex>@{}[d]|{1}\ar[r]^{\textstyle \iota 1_i^*}&\ar[ld]\ar@<1ex>@{}[ld]|(.25){ \chi
}\ar[d]^{\textstyle 1_i^*(1_{1_i^*x})}1_i^*1_i^*x\\\ar@{}[ru]|(.7){\textstyle \cong
}&1_i^*x&1_i^*1_i^*x\ar[l]^{\textstyle \chi}}
$$
\end{proof}
For any  lax diagram of bicategories $\f:I^\mathrm{op}\to \bicat$, the lax $I$-homomorphism $J:\f\to \f^{^\mathrm{r}}$
induces a corresponding homomorphism on the Grothendieck constructions
$\xymatrix{\int_I\!J:\int_I\f\to\int_I\f^{^\mathrm{r}}}$. Up to a pseudo-natural equivalence, this homomorphism $\int_I\!J$ can easier be
described in terms of the following normal homomorphism

\begin{equation}\label{inj}\xymatrix{\textbf{j}:\int_I\f\to\int_I\f^{^\mathrm{r}},}
\end{equation}
$$ \xymatrix@C=7pt@R=7pt{ (y,j) \ar@/^1.1pc/[rr]^{\textstyle (u,a)} \ar@/_1.1pc/[rr]_{\textstyle (u',a)} &
\Downarrow(\alpha,a) & (x,i)}~~~ \overset{\textstyle \textbf{j}}\mapsto ~~~\xymatrix@C=7pt@R=7pt{ ((y,1_j),j)
\ar@/^1.1pc/[rr]^{\textstyle ((u,a),a)} \ar@/_1.1pc/[rr]_{\textstyle ((u',a),a)} & \Downarrow\! ((\alpha,a),a) &
((x,1_i),i)),}
$$
whose structure constraints for horizontal compositions of 1-cells are given by the left identity constraints of the
bicategories $\f^{^\mathrm{r}}_i$.
\begin{lemma}\label{simj} There is a pseudo-natural equivalence
$\int_I\!J\Rightarrow\mathbf{j}:\int_I\f\to\int_I\f^{^\mathrm{r}}$.
\end{lemma}
\begin{proof} The claimed pseudo-natural equivalence is the identity transformation on objects and, at  each 1-cell
$(u,a):(y,j)\to(x,i)$ of the bicategory $\int_I\!J$, its component is the composite 2-cell
$$\xymatrix@C=45pt{1_{((x,1_i),i) }\circ
\int_I\!J(u,a)\overset{\boldsymbol{l}}{\cong}\int_I\!J(u,a)\ar@{=>}@<-2pt>[r]^{\textstyle ((\alpha,a),a)}&
\mathbf{j}(u,a)\overset{\boldsymbol{r}^{-1}}{\cong} \mathbf{j}(u,a)\circ 1_{((y,1_j),j)}, }$$ where the invertible
2-cell $\alpha$ is that obtained by pasting the diagram in $\f_j$
$$\xymatrix@R=12pt{y\ar[r]^{\textstyle u}\ar[rrddd]_{\textstyle u}& a^*x\ar@{}[ddr]_{\textstyle \cong}\ar[rr]^{\textstyle \iota a^*}\ar[dr]^{1_{a^*x}} && 1_j^*a^*x \ar[r]^{\textstyle 1_j^*(1_{a^*x})} & 1_ja^*x \ar[dddll]^{\textstyle \chi}\\
&& a^*x\ar@{}[rd]|(.45){\textstyle \underset{\textstyle \cong}{\gamma}}\ar[rru]_{\iota a^*} \ar[dd]^{1_{a^*x}} \ar@{}[ru]^(0.3){\textstyle  \cong \!{\widehat{\iota}}}&&\\
&&&&\\
&& a^*x && }$$
\end{proof}
\begin{proposition}\label{proad}For any  lax diagram of bicategories $\f:I^\mathrm{op}\to \bicat$, the homomorphism
$\xymatrix{\mathbf{j}:\int_I\!\f\to\int_I\!\f^{^\mathrm{r}}}$ has a left biadjoint.
\end{proposition}
\begin{proof} The left biadjoint to $\mathbf{j}$ is the normal homomorphism $\mathbf{p}: \xymatrix{\int_I\!\f^{^\mathrm{r}}\to\int_I\!\f}$ defined by
$$
 \xymatrix@C=7pt@R=7pt{ ((y,j\overset{d}\to l),j) \ar@/^1.1pc/[rr]^{\textstyle ((u,b),a)} \ar@/_1.1pc/[rr]_{\textstyle ((u',b),a)} &
\Downarrow((\alpha,b),a) & ((x,i\overset{c}\to k),i)}~~~ \overset{\textstyle \mathbf{p}}\mapsto ~~~\xymatrix@C=7pt@R=7pt{
(y,l)\ar@/^1pc/[rr]^{\textstyle (u,b)} \ar@/_1pc/[rr]_{\textstyle (u',b)} & \Downarrow
(\alpha,b)& (x,k),}
$$
whose structure constraints for horizontal compositions of 1-cells are given by the left identity constraints of the
bicategories $\f^{^\mathrm{r}}_i$.

The unit of the biadjunction is the pseudo-transformation $\eta:1 \Rightarrow \mathbf{j p}$, with
$$\eta ((x,i\overset{c}\to k),i)=((\iota_kx,1_k),c):((x,i\overset{c}\to k),i)\to ((x,k\overset{1_k}\to k),k),$$ for each object $((x,i\overset{c}\to k),i)$ of $\int_I\!\f^{^\mathrm{r}}$, and whose component   $$\widehat{\eta}:((\iota_kx,1_k),c)\circ ((u,b),a)\cong ((u,b),b)\circ (\iota_ly,1_l),d),$$ at a 1-cell ${((u,b),a):((y,j\overset{d}\to l),j)\to ((x,i\overset{c}\to k),i)}$, is provided by the 2-cell obtained by pasting in the bicategory $\f_l$ $$
\xymatrix{y\ar@{}[rd]|{\textstyle \overset{\textstyle\widehat{\iota}}{\cong}}\ar[d]_{\textstyle \iota}\ar[r]^{\textstyle u}&b^*x\ar[d]\ar@{}@<1ex>[d]_{\iota b^*}\ar[rd]\ar@{}@<-1ex>[rd]^{1_{b^*}}\ar[r]^{\textstyle b^*\iota}&b^*1_k^*x\ar[d]^{\textstyle \chi}\\1_l^*y\ar[r]_{\textstyle 1_l^*u}&1_l^*b^*x\ar[r]_{\textstyle \chi}\ar@{}[ru]|(.3){\textstyle
\cong\!\gamma}|(.7){\textstyle \cong\!\delta}&b^*x.
}
$$
 One easily sees the equalities  $\mathbf{p j}=1_{\int_I\!\f}$, $\eta\mathbf{j}=1_{\mathbf{j}}$, and $\mathbf{p}\eta=1_{\mathbf{p}}$, showing that $\mathbf{p}\dashv \mathbf{j}$ is a biadjunction.
\end{proof}

\section{classifying spaces}

For the general background on simplicial sets,  we mainly refer  to \cite{g-j}. The {\em simplicial category} is denoted by $\Delta$, and its objects, that is, the ordered sets ${[n]=\{0,1,\dots,n\}}$, are usually considered as categories with only one morphism  $j\rightarrow i$ when $0\leq i\leq j\leq n$. Then, a non-decreasing map $[n]\rightarrow [m]$ is the same as a functor, so that we see $\Delta $, the simplicial
category of finite ordinal numbers, as a full subcategory of $\cat$, the category (actually the 2-category)
of small categories. Recall that the category $\Delta$ is generated by the injections
$d^i:[n-1]\to[n]$ (cofaces), $0\leq i\leq n$, which omit  the $i$th element and the
surjections $s^i:[n+1]\to [n]$ (codegeneracies), $0\leq i\leq n$, which repeat the $i$th
element, subject to the well-known  {\em cosimplicial identities}:
$d^j d^i=d^i d^{j-1}$ if $i<j$, etc.

Given a bicategory $\c$, let \begin{equation}\label{clasi}\class\c\end{equation} denote its {\em classifying space}. We shall briefly recall from \cite{ccg} that $\class\c$ can be defined through several,  but always homotopy-equivalent, constructions. For instance, $\class \c$ may be thought of as the realization of the  normal pseudo-simplicial category, called the {\em pseudo-simplicial nerve} of the bicategory,
\begin{equation}\label{ps1.1}\ner\c=(\ner\c,\chi,1): \Delta^{\!^{\mathrm{op}}}\to \cat\,,
\end{equation}
whose category of $p$-simplices is
$$
\ner_{\!p}\c = \bigsqcup_{(x_0,\ldots,x_p)\in \mbox{\scriptsize Ob}\c^{p+1}}\hspace{-0.3cm}
\c(x_1,x_0)\times\c(x_2,x_1)\times\cdots\times\c(x_p,x_{p-1}),
$$
where a typical arrow is a string of 2-cells in $\c$
$$\xymatrix @C=8pt {x_0  & {\Downarrow\, \alpha_1} & x_1 \ar@/^1pc/[ll]^{\textstyle v_1}
\ar@/_1pc/[ll]_-{\textstyle u_1} & {\Downarrow\, \alpha_2} & x_2  \ar@/^1pc/[ll]^{\textstyle v_2}
\ar@/_1pc/[ll]_-{\textstyle u_2}&\hspace{-15pt}\cdots&\hspace{-10pt}
 x_{n-1}
  & {\Downarrow\, \alpha_n} & x_n
\ar@/^1pc/[ll]^{\textstyle v_n} \ar@/_1pc/[ll]_-{\textstyle u_n},}$$
and $\ner\c_0=\mbox{Ob}\,\c$, as a discrete category. The face  and degeneracy functors are defined in the standard way by using the horizontal composition of adjacent cells and the identity morphisms of the bicategory:
\begin{equation}\label{p.1.2}\begin{array}{l} d_i(\alpha_1,\dots,\alpha_p)=\left\{
\begin{array}{lcl}
\hspace{-0.2cm}  (\alpha_2,\dots,\alpha_p) &\text{ if }& i=0 ,  \\[4pt]
\hspace{-0.2cm}  (\alpha_1,\dots,\alpha_i\circ\alpha_{i\text{+}1},\dots,\alpha_p) & \text{ if } & 0<i<p , \\[4pt]
\hspace{-0.2cm}  (\alpha_1,\dots,\alpha_{p-1}) & \text{ if }& i=p,
\end{array}\right.\\[20pt]
s_i(\alpha_1,\dots,\alpha_p)=(\alpha_1,\dots,\alpha_i,1_{x_i},\alpha_{i\text{+}1},\dots, \alpha_p).
\end{array}
\end{equation}
If $a:[q]\to[p]$ is any non-identity map in $\Delta$, then we write
$a$ in the (unique) form (see \cite{may67}, for example) $a=d^{i_1}\cdots d^{i_s}s^{j_1}\cdots s^{j_t}$, where
$0\leq i_s<\cdots< i_1 \leq p$, $0\leq j_1<\cdots< j_t \leq q$ and $q+s=p+t$, and the induced functor
$a^*:\ner_{\!p}\c\to\ner_q\c$ is defined by $a^*=s_{j_t}\cdots
s_{j_1}d_{i_s}\cdots d_{i_1}$. Note that $d_jd_i=d_id_{j+1}$ for $i\leq j$, unless $i=j$ and $1\leq i\leq p-2$, in
which case the associativity constraint of $\c$ gives a  canonical natural isomorphism
\begin{equation}\label{ps1.2} d_id_i \overset{\chi}\cong d_id_{i+1}. \end{equation}
Similarly, all the equalities $d_0s_0=1$, $d_{p+1}s_p=1$, $d_is_j=s_{j-1}d_i$ if $i<j$ and $d_is_j=s_jd_{i-1}$ if $i>j+1$,
hold, and the unit constraints of $\c$ give canonical isomorphisms
\begin{equation}\label{ps1.3} d_is_i\overset{\chi}\cong 1, \hspace{0.6cm} d_{i+1}s_i\overset{\chi}\cong 1. \end{equation}

Then it is a fact that this family of natural isomorphisms (\ref{ps1.2}) and (\ref{ps1.3}), uniquely
determines a whole system of natural isomorphisms
$\chi_{a,b}: b^*a^*\cong (ab)^*$,
one for each pair of composible maps in $\Delta$, $[n]\overset{b}\to [q]\overset{a}\to [p]$, such that the
assignments $a\mapsto a^*$, $1_{[p]}\mapsto 1_{\ner\c_p}$, together with these isomorphisms
$b^*a^*\cong (ab)^*$, give the data for the pseudo-simplicial category
(\ref{ps1.1}), $\ner\c:\Delta^{\!^{\mathrm{op}}}\to\cat$. This fact can be easily proven by using
Jardine's supercoherence theorem \cite[Corollary 1.6]{jardine} since the commutativity of the seventeen diagrams of supercoherence, (1.4.1)-(1.4.17) in \cite{jardine},
 easily follows from the pentagon and triangle coherence diagrams in the bicategory $\c$.

When a category $\c$ is considered as a discrete bicategory, that is, where the deformations are all identities, then $\ner \c$ is the usual Grothendieck's nerve of the category.

 Since the horizontal composition involved is in general neither strictly associative nor unitary, $\ner\c$ is not a simplicial category (with a well understood simple geometric realization), which forces one to deal with defining the geometric realization of what is not simplicial but only `simplicial up to (coherent) isomorphisms'. Indeed, this has been done by Segal, Street, and Thomason using methods of Grothendieck, so that the classifying space of the bicategory is
$$\class\!\xymatrix{ \int_{\!\Delta}\!\ner\c,}$$
the ordinary classifying space of the category obtained as the Grothendieck construction on the pseudo-simplicial nerve of the bicategory $\ner\c$.

   A second possibility is to recall that the  {\em unitary geometric nerve of a bicategory} $\c$ \cite{street, duskin, gurski2, ccg} is the simplicial set
\begin{equation}\label{ngn} \Delta^{\hspace{-2pt}^\mathrm{u}}\!\c:\Delta^{\!^{\mathrm{op}}}\ \to \ \set, \hspace{0.6cm}[p]\mapsto \mathrm{NorLaxFunc}([p],\c),\end{equation}
whose $p$-simplices are the normal lax functors ${\xi:[p]\to \c}$. If $a:[q]\to [p]$ is any map in $\Delta$, that is, a functor,  the induced $a^*:\Delta^{\hspace{-2pt}^\mathrm{u}}\!\c_p\to\Delta^{\hspace{-2pt}^\mathrm{u}}\!\c_q$ carries ${\xi:[p]\to \c}$ to the composite ${\xi a:[q]\to \c}$, of $\xi$ with $a$.
This nerve $\Delta^{\hspace{-2pt}^\mathrm{u}}\!\c$ is   a simplicial set which is coskeletal in dimensions greater than 3, whose vertices are the objects $\xi 0$ of $\c$, the 1-simplices are the 1-cells $\xi_{0,1}:\xi 1\to \xi 0$  and, for $p\geq 2$, a $p$-simplex of $\Delta^{\hspace{-2pt}^\mathrm{u}}\c$
is geometrically represented by a diagram in $\c$ with the shape of the 2-skeleton of an
orientated standard $p$-simplex, whose faces are triangles
$$
\xymatrix{\ar@{}[drr]|(.6)*+{\!\Downarrow \xi_{{i,j,k}}}
                & \xi_j \ar[dr]^{\textstyle \xi_{i,j}}             \\
                     \xi_k   \ar[ur]^{\textstyle \xi_{j,k}} \ar[rr]_{\textstyle \xi_{i,k}} && \xi_i     }
$$
with objects $\xi i$ placed on the vertices, $0\leq i\leq p$, 1-cells $\xi_{i,j}:\xi j\rightarrow \xi i$ on the edges, $0\leq i<j\leq p$, and 2-cells $\xi_{i,j,k}:\xi_{i,j}\circ \xi_{j,k}\Rightarrow \xi_{i,k}$, for $0\leq i<j<k\leq p$. These data are required to satisfy the condition that, for $0\leq i<j<k<l\leq p$,  each tetrahedron is commutative in the sense that
$$
\xymatrix@C=30pt@R=30pt{\xi_l\ar[r]\ar[rd]\ar[d]& \xi_i  \ar@{}[ld]|(0.3){\textstyle \Uparrow}|(0.7){\textstyle \Rightarrow}\ar@{}@<7ex>[d]|{=}\\
\xi_k\ar[r]&\xi_j\ar[u]}\hspace{0.6cm}\xymatrix@C=30pt@R=30pt{\xi_l\ar[r]\ar[d]\ar@{}[rd]|(0.3){\textstyle \Uparrow}|(0.7){\textstyle \Leftarrow}& \xi_i \\
\xi_k\ar[r]\ar[ru]&\xi_j\ar[u] }
$$

The {\em geometric nerve of a bicategory} $\c$ is the simplicial set
 \begin{equation}\label{geb}\Delta\c:\Delta^{\!^{\mathrm{op}}}\ \to \ \set,\hspace{0.6cm}
[p]\mapsto \lfunc([p],\c),\end{equation}
that is, the simplicial set whose $p$-simplices are all lax functors ${\xi:[p]\to \c}$.
  Hence, the unitary geometric nerve  $\Delta^{\hspace{-2pt}^\mathrm{u}}\c$ becomes a simplicial subset of $\Delta\c$.  The $p$-simplices of the  geometric nerve $\Delta\c$ are described similarly to those of the normalized one, but now they include  2-cells
$\xi_i:1_{\xi i}\Rightarrow \xi_{i,i}$, $0\leq i\leq p$, with the
requirement that the diagrams below commute.
$$
\xymatrix@C=9pt@R=9pt{& \xi_{i,j}\circ 1_{\xi j} \ar@2{->}[dl]\ar@{}@<-3pt>[dl]_(0.6){\textstyle 1\circ \xi_j} \ar@2{->}[dr]^{\textstyle \boldsymbol{r}} \\
\xi_{i,j}\circ \xi_{j,j} \ar@2{->}[rr]_{\textstyle \xi_{i,j,j}} && \xi_{i,j} } \hspace{1cm} \xymatrix@C=9pt@R=9pt{& 1_{\xi i}\circ \xi_{i,j} \ar@2{->}[dl]_{\textstyle \xi_i \circ 1} \ar@2{->}[dr]^{\textstyle \boldsymbol{l}} \\
\xi_{i,i}\circ \xi_{i,j} \ar@2{->}[rr]_{\textstyle \xi_{i,i,j}} && \xi_{i,j} }
$$

We shall list below a number of required results from  \cite{ccg}:

\begin{fact}\label{f1}\cite[Theorem 6.1]{ccg} For any bicategory $\c$, there are natural homotopy equivalences
\begin{equation}\label{ef1}\class \c \simeq |\Delta^{\hspace{-2pt}^\mathrm{u}}\c|\simeq |\gner\c|.\end{equation}
\end{fact}

\begin{fact}\label{f2}\cite[(30) and Theorem 7.1]{ccg} $(i)$ Any homomorphism between bicategories $F:\b\to \c$ induces a continuous cellular map $\class F:\class\b\to\class\c$. Thus, the classifying space construction, $\c\mapsto \class\c$, defines a functor from the category of bicategories and homomorphisms to CW-complexes.

\vspace{0.2cm}
\noindent $(ii)$ If $F,F':\b\to\c$ are two homomorphisms between bicategories, then any lax (or oplax)
transformation, $F\Rightarrow F'$, canonically defines a homotopy between the induced maps on classifying spaces, $ \class F\simeq \class F':\class\b\to\class\c$.

\vspace{0.2cm}
\noindent (iii) If a homomorphism of bicategories has a left or right biadjoint, the map induced on classifying spaces is a homotopy
equivalence. In particular, any biequivalence of bicategories induces a homotopy equivalence on classifying spaces.
\end{fact}

\begin{fact}\label{f3}\cite[Theorem 7.3]{ccg} Suppose a category $I$ is given. For every functor $\f:I^{{\mathrm{op}}}\to\mathbf{Hom}\subset \bicat$, there exists a natural weak homotopy equivalence of simplicial sets
$$\xymatrix{ \hoco_I\gner\f\overset{\sim}\longrightarrow
\gner\!\int_I\!\f,}$$  where $\hoco_I\gner\f$ is the homotopy colimit construction by Bousfield and Kan \cite[\S XII]{bousfield-kan} of the diagram of simplicial sets $\gner\f:I^{{\mathrm{op}}}\to\sset$, obtained by composing $\f$ with the
geometric nerve functor $\gner:\mathbf{Hom}\to \sset$, and $\int_I\!\f$ is the bicategory obtained by the Grothendieck construction on $\f$.
\end{fact}

In \cite{segal68}, Segal extended Milnor's geometric realization process, $S\mapsto |S|$, to simplicial (compactly generated topological)
spaces, which provides, for instance, the notion of {\em classifying spaces for  simplicial bicategories} $\f:\Delta^{\!{^\mathrm{op}}}\to \mathbf{Hom}\,$. By replacing each bicategory $\f_p$, $p\geq 0$, by its classifying space
$\class\f_p$, one obtains a simplicial space, $[p]\mapsto \class\f_p$,  whose
Segal realization is, by definition, the classifying space  of the simplicial bicategory. But
note, as a consequence of Fact \ref{f1} and \cite[Lemma p.86]{quillen}, that there are  homotopy equivalences
\begin{equation}\label{fdia} |[p]\mapsto \class\f_p| \simeq |[p]\mapsto |\gner\f_p|| \simeq |\diag \Delta \f |,\end{equation}
where $\diag \Delta \f $ is the simplicial set diagonal of the bisimplicial set
obtained by composing the geometric nerve functor ${\Delta:\mathbf{Hom} \to\sset} $ with $\f$, that is, $$\Delta\f:([p],[q])\mapsto \lfunc([q],\f_p).$$

The above construction, for simplicial bicategories, leads to the more
general notion of {\em classifying space for diagrams of bicategories}: If ${\f: I^{{\mathrm{op}}}\to \mathbf{Hom}}$ is a functor,  where $I$ is any category, then
one applies the so-called {\em Borel construction}, obtaining the simplicial bicategory
$$
 E_I\f:   \Delta^{\!{^\mathrm{op}}}\to \mathbf{Hom}, \hspace{0.6cm}
 [p] \mapsto \bigsqcup\limits_{[p]\overset{\beta}\to I}\! \f_{\beta 0},
$$
where the disjoint union is over all functors $\beta:[p]\to I$ (i.e., the $p$-simplices of the nerve $\ner I=\gner I$). The induced homomorphism by a map $a:[q]\to [p]$, in $\Delta$, applies the bicategory component at $\beta:[p]\to I$ into the component  at the composite $\beta a:[q]\to I$, just by the homomorphism of bicategories $$\beta_{0,{\textstyle a} 0}^*:\f_{\beta 0}\to \f_{\beta{\textstyle a} 0}$$ attached
in diagram $\f:I^{{\mathrm{op}}}\to\mathbf{Hom} $ at the morphism $\beta_{0, {\textstyle a} 0}:\beta a0\to \beta0$ of $I$. Then,  the classifying space of the diagram of bicategories $\f: I^{{\mathrm{op}}}\to \mathbf{Hom}$ is
the classifying space, in the above sense,  of the simplicial bicategory  $E_I\f$.
But note that $$\diag\gner E_I\f=\hoco_I\gner\f,$$ that is, the simplicial set
$$
 [p] \mapsto \bigsqcup\limits_{[p]\overset{\beta}\to I} \lfunc([p],\f_{\beta 0}),
$$
and therefore, by (\ref{fdia}), the classifying space of $\f$ is homotopy equivalent to
$$|\hoco_I\Delta\f|,$$
the geometric realization of the homotopy colimit \cite{bousfield-kan} of the simplicial set diagram  $\Delta \f:I^{{^\mathrm{op}}}\to \sset$,   obtained by composing ${\f}$ with the geometric nerve functor $\gner:\mathbf{Hom}\to \sset$.
Since, for any simplicial bicategory $\f$, we have a natural weak homotopy
equivalence of simplicial sets \cite[XII, 4.3]{bousfield-kan} $\mbox{hocolim}_\Delta\Delta\f\overset{\thicksim}\to \mbox{diag}\Delta\f$, it follows that
both constructions above for the classifying space  of a simplicial
bicategory $\f:\Delta^{\!{^\mathrm{op}}}\to \mathbf{Hom}$ coincide up to a natural homotopy equivalence.

Furthermore,  the classifying space of any diagram $\f: I^{{\mathrm{op}}}\to \mathbf{Hom}$ is homotopy equivalent to the one of the
bicategory obtained by the Grothendieck construction on it, $\int_I\!\f$, thanks to
  the existence of the natural homotopy equivalences
\begin{equation}\label{cg}|\hoco_I\Delta\f|\simeq |\gner\!\!\xymatrix{\int_I\!\f}|\simeq \class\!\!\xymatrix{\int_I\!\f}, \end{equation}
by Facts \ref{f3} and \ref{f1} respectively. This suggests the following general definition for lax diagrams of bicategories:
\begin{definition} \label{defclass}The classifying space of a lax $I$-diagram  $\f: I^\mathrm{op}\to \bicat$, denoted $\class \f$,  is defined to be the classifying space of the bicategory obtained by the Grothendieck construction on $\f$.
\end{definition}

So defined, the assignment $\f\mapsto \class\f$ has the following basic properties:

\begin{proposition} \label{propclass}$(i)$ If $\f,\g:I^\mathrm{op}\to \bicat$ are lax $I$-diagrams, then each
lax $I$-homomorphism $F:\f\to \g$ induces a continuous map $\class F:\class\f\to\class\g$.

\vspace{0.2cm} \noindent $(ii)$ Any pseudo $I$-transformation,
$F\Rightarrow G:\f\rightarrow \g$ induces a homotopy  $\class
F\simeq \class G$.

\vspace{0.2cm}\noindent $(iii)$ For any lax $I$-diagram $\f$,  there is a homotopy $\class1_\f\simeq 1_{\class\f}$. For  any pair of composible lax $I$-homomorphisms $\f\overset{F}\to \g \overset{G}\to \h$, there is a homotopy  $$\class G\, \class F\simeq \class (GF).$$
\end{proposition}
\begin{proof}
$(i)$ As in (\ref{G.7}), the lax $I$-homomorphism $F:\f\to\g$
defines the homomorphism of bicategories
$\xymatrix{\int_I\!F:\int_I\f\rightarrow \int_I\g}$ which, by Fact
\ref{f2} (i), determines the claimed cellular map $\class
F:\class\f\rightarrow \class\g$.

\vspace{0.2cm}\noindent $(ii)$ As in (\ref{pstint}), any pseudo $I$-transformation  $m:F\Rightarrow G$ gives rise
to a pseudo-transformation
$\xymatrix{ \int_I\!m:\int_I\!F\Rightarrow \int_I\!G}$,
 which, by Fact \ref{f2} (ii), determines a homotopy $\class F\simeq\class G$.

\vspace{0.2cm}\noindent $(iii)$ The announced homotopies are respectively induced,  from Fact \ref{f2} (i), (ii),  by the pseudo-natural equivalences  (\ref{gc4}) and (\ref{gc3}).
\end{proof}

We have seen that for a diagram, that is, a functor,  $\f: I^\mathrm{op}\to \mathbf{Hom}$, both Borel and Grothendieck constructions lead to the same space $\class\f$, up to a natural homotopy equivalence. Next,  we show that the classifying space construction for lax diagrams of bicategories is consistent with
the so-called {\em rectification} process, $\f \mapsto \f^{^\mathrm{r}}$, developed in Section \ref{rectification}. Recall that this process associates
to any lax diagram  $\f:I^{{\mathrm{op}}}\to \bicat$ a  genuine diagram
$\f^{^\mathrm{r}}:I^{{\mathrm{op}}}\to \mathbf{Hom}\subset\bicat$.

\begin{proposition}\label{hip} Given $\f:I^{\mathrm{op}}\to \bicat$ any lax $I$-diagram of bicategories, the
lax $I$-homomorphism $ J:\f\to \f^{^\mathrm{r}}$ in Lemma {\em \ref{lemj}} induces a homotopy equivalence
$$\class J:\class\f\overset{_\simeq}\longrightarrow \class\f^{^\mathrm{r}}.$$

If $F:\f\to \g$ is any lax $I$-homomorphism between lax $I$-diagrams, then the induced diagram below is homotopy commutative.
$$
\xymatrix{\ar@{}@<0.5ex>[r]_{\simeq}\class \f\ar[r]^{\textstyle \class J}\ar[d]_{\textstyle \class F}&\class \f^{^\mathrm{r}} \ar[d]^{\textstyle \class F^{^\mathrm{r}}}\\
\class\g \ar[r]_{\textstyle \class J}\ar@{}@<-0.5ex>[r]^{\simeq}&\class \g^{^\mathrm{r}}
}
$$
\end{proposition}
\begin{proof}
By Lemma \ref{simj} and Fact \ref{f2} (ii), the map $\class J$ is
homotopic to the induced map $\class \mathbf{j}:\class\f\to\class
\f^{^\mathrm{r}}$ by the homomorphism (\ref{inj}), $
\xymatrix{\textbf{j}:\int_I\f\to\int_I\f^{^\mathrm{r}}}$, which, by
Proposition \ref{proad}, has a left biadjoint and therefore induces
a homotopy equivalence on classifying spaces, by Fact \ref{f2}
(iii). Hence, $\class J$ is a homotopy equivalence.

For the square in the proposition, note that, by Proposition \ref{propclass} (iii), there are homotopies $\class F^{^\mathrm{r}}\class J\simeq \class(F^{^\mathrm{r}}J)$ and  $\class J\, \class F\simeq \class(J\,F)$. Since the pseudo $I$-equivalence (\ref{mj}), $m:F^{^\mathrm{r}}J\Rightarrow JF$, induces a homotopy $\class(F^{^\mathrm{r}}J)\simeq \class(JF)$, by Proposition \ref{propclass} (ii), the result follows. \end{proof}

The following main theorem extends to lax diagrams of bicategories a well-known result by Thomason \cite[Corollary 3.3.1]{thomason} for lax diagrams of categories:

\begin{theorem}\label{p.1.23} If $F:\f\to\g$ is a lax $I$-homomorphism between lax $I$-diagrams $\f,\g:  I^{{\mathrm{op}}}\to\bicat$, such that the induced maps $BF_i:\class \f_i\to \class\g_i$ are homotopy equivalences, for all objects $i$ of
$I$, then the induced map $\class F:\class\f \to \class\g$ is a homotopy equivalence.
\end{theorem}
\begin{proof}
By Proposition \ref{hip} above, it suffices to prove that the induced map after rectification   ${\class F^{^\mathrm{r}}:\class\f^{^\mathrm{r}} \to \class\g^{^\mathrm{r}}}$ is a homotopy equivalence. Let us recall from Section \ref{rectification} that both $\f^{^\mathrm{r}}$ and $\g^{^\mathrm{r}}$ are genuine diagrams of bicategories, that is, functors $I^{{\mathrm{op}}}\to\mathbf{Hom}$, and $F^{^\mathrm{r}}:\f^{^\mathrm{r}} \to \g^{^\mathrm{r}}$ is merely  a natural transformation.
Then, by the natural homotopy equivalences (\ref{cg}), it will be enough to prove that the natural transformation $\gner F^{^\mathrm{r}}$ between the functors $\gner\f^{^\mathrm{r}}, \gner\g^{^\mathrm{r}}:I^{{\mathrm{op}}}\to\sset$,
induces a weak homotopy equivalence on the corresponding homotopy colimits $$\hoco_I\gner F^{^\mathrm{r}}:\hoco_I\gner \f^{^\mathrm{r}} \overset{\sim}\to \hoco_I\gner\g^{^\mathrm{r}}.$$

For, let us observe that, for each object $i$ of the category $I$, the square
$$
\xymatrix{\ar@{}@<0ex>[d]^{\simeq}\ar@{}@<0.5ex>[r]_{\simeq}\class \f_i\ar[r]^{\textstyle \class J_i}\ar[d]_{\textstyle \class F_i}&\class \f^{^\mathrm{r}}_i \ar[d]^{\textstyle \class F^{^\mathrm{r}}_i}\\
\class\g_i \ar[r]_{\textstyle \class J_i}\ar@{}@<-0.5ex>[r]^{\simeq}&\class \g^{^\mathrm{r}}_i
}
$$
is homotopy commutative because of the pseudo-natural equivalence  (\ref{mj}) at $i$, ${m_i: F^{^\mathrm{r}}_iJ_i\Rightarrow J_iF_i:\f_i\to\g^{^\mathrm{r}}_i}$ (see Fact \ref{f2} (i), (ii)). Moreover,  both maps $\class J_i$ in the square are homotopy equivalences, since all the homomorphisms $J_i$ have a right biadjoint by Lemma \ref{lemj2} (see Fact \ref{f2}(iii)). Since, by hypothesis, the map $BF_i:\class \f_i\to \class\g_i$ is also a homotopy equivalence, it follows that the remaining map in the square has the same property, that is, the map $\class F^{^\mathrm{r}}_i:\class \f^{^\mathrm{r}}_i\to \class \g^{^\mathrm{r}}_i$ is a homotopy equivalence.
By taking into account Fact \ref{f1}, the above means that, for every object $i$ of $I$, the induced simplicial map on geometric nerves $\gner F^{^\mathrm{r}}_i:\gner \f^{^\mathrm{r}}_i\to \gner \g^{^\mathrm{r}}_i$ is a weak homotopy equivalence, whence, by the Homotopy Lemma \cite[XII, 4-2]{bousfield-kan}, the result follows, that is, the simplicial map $\hoco_I\gner F^{^\mathrm{r}}$ is a weak homotopy equivalence.
\end{proof}

For any lax diagram $\f: I^\mathrm{op}\to \bicat$, the bicategory $\int_I\!\f$ assembles all bicategories $\f_i$ in the following sense: There is a projection 2-functor $\xymatrix{\mathbf{q} :\int_I\!\f \to I}$,
$$
 \xymatrix@C=2pt{(y,j)  \ar@/^1pc/[rr]^{\textstyle (u,a)} \ar@/_1pc/[rr]_{\textstyle (u',a)} & {}_{\textstyle \Downarrow(\alpha,a)} &(x,i)&\overset{\textstyle \mathbf{q}}\mapsto\underset{~}~& \hspace{-0.3cm}j\ar[rr]^-{\textstyle a}&&i \,,}
$$
and, for each object $i$ of $I$, there is a commutative square
\begin{equation}\label{square}
\xymatrix@C=18pt@R=18pt{\f_i\ar[d]_{\textstyle \eta_i}\ar[r]&[0] \ar[d]^{\textstyle i}\\\int_I\!\f\ar[r]^{\textstyle \mathbf{q}} &I}
\end{equation}
where  $\eta_i:\f_i\to \int_I\!\f$ is the embedding homomorphism
described in (\ref{eeta}).

\begin{theorem}\label{hc2} Suppose that $\f:I^{^{\mathrm{op}}}\to \bicat$ is a lax $I$-diagram of bicategories such that the induced map
$\class a^*:\class \f_{i}\to \class \f_{j}$, for each morphism $a:j\to i$ in  $I$,  is a homotopy equivalence. Then, for every object $i$ of $I$, the square induced by $(\ref{square})$
\begin{equation}\label{square2}
\xymatrix@C=18pt@R=18pt{\class\f_i\ar[r]\ar[d]&\star \ar[d]^i\\ \class \f\ar[r]&\class I}
\end{equation}
is homotopy cartesian. Therefore, for each object $x\in \f_i$ , there is an induced long exact sequence on homotopy groups relative to the base points $x$ of $\class\f_i$, $(x,i)$ of $\class \f$, and $i$ of $\class I$,
$$\xymatrix{\cdots \to \pi_{n+1}\class I\to\pi_n\class\f_i\to\pi_n\class \f\to\pi_n\class I\to\cdots.}$$
\end{theorem}
\begin{proof}
The square (\ref{square}) is the composite of the squares
$$
\xymatrix@C=18pt@R=22pt{
\f_i\ar[r]^{\textstyle J_i} \ar[d]_{\textstyle \eta_i\f}\ar@{}[rd]|(.45){(a)}
&
\f^{^\mathrm{r}}_i \ar[d]|(0.4){\eta_i\f^{^\mathrm{r}}}\ar[r] \ar@{}[rd]|(.50){(b)}
&
[0]\ar[d]^{\textstyle i}
\\
\int_I\!\f \ar[r]^{\textstyle \mathbf{j}}
&
\int_I\!\f^{^\mathrm{r}} \ar[r]^{\textstyle \mathbf{q}} &I
}
$$
where, in $(a)$, both horizontal homomorphisms $J_i$ (\ref{defji}) and $\mathbf{j}$ (\ref{inj}) induce homotopy equivalences on classifying spaces, by Lemma \ref{lemj2}, Proposition \ref{proad}, and Fact \ref{f2} (iii). Therefore, the induced square (\ref{square2}) is homotopy cartesian if and only if the one induced by $(b)$ is as well. But, recall that the rectification $\f^{^\mathrm{r}}:I^{^{\mathrm{op}}}\to \mathbf{Hom}\subset \bicat$ is a diagram, that is, a functor,  and we have the natural homotopy equivalences (\ref{cg}). Therefore,  it will be enough to prove that the induced pullback square of spaces
$$
\xymatrix@C=18pt@R=18pt{\class\f^{^\mathrm{r}}_i\ar[r]\ar[d]&\star \ar[d]^i\\ \class E_I\f^{^\mathrm{r}} \ar[r]&\class I}
$$
is homotopy cartesian, which is a consequence of Quillen's Lemma \cite[p. 90]{quillen} (see also \cite[\S IV, Lemma 5.7]{g-j}). To verify the hypothesis, simply note that, for each arrow $a:j\to i$ in $I$, the square
$$
\xymatrix@R=18pt{\class\f_i\ar[r]^{\textstyle \class J_i}\ar[d]_{\textstyle \class a^*}&\class\f^{^\mathrm{r}}_i\ar[d]^{\textstyle \class a^*}\\
\class\f_j\ar[r]^{\textstyle \class J_j}&\class\f^{^\mathrm{r}}_j}
$$
is homotopy commutative thanks to the pseudo-transformation (\ref{ctheta}) and Fact \ref{f2} (ii). Since the horizontal induced maps $\class J_i$ and $\class J_j$ are both homotopy equivalences by Lemma \ref{lemj2} and Fact \ref{f2} (iii), as well as the map $\class a^*: \class \f_i\to \class \f_j$,by hypothesis, we conclude that the map $\class a^*:\class\f^{^\mathrm{r}}_i\to \class\f^{^\mathrm{r}}_j$ is also a homotopy equivalence.
\end{proof}

\section{Classifying spaces of Braided Monoidal Categories}
Lax diagrams of bicategories form the foundation for the classifying spaces theory of (small) tricategories: any tricategory $\mathcal{T}=(\mathcal{T},\otimes, \mathrm{I}, \boldsymbol{a},\boldsymbol{l},\boldsymbol{r}, \pi,\mu,\lambda,\rho)$, as in \cite{g-p-s}, has associated a pseudo-simplicial bicategory, called its {\em nerve},
\begin{equation}\label{psnt}\ner\mathcal{T}=(\ner\mathcal{T},\chi,\iota, \omega,\gamma,\delta): \Delta^{\!^{\mathrm{op}}}\to \bicat,
\end{equation}
and  the {\em classifying space of the tricategory} is  the classifying space of its bicategorical pseudo-simplicial nerve . Briefly, say that the bicategory of $p$-simplices  of $\ner\mathcal{T}$ is
$$
\ner_p\mathcal{T} = \bigsqcup_{(x_0,\ldots,x_p)\in \mbox{\scriptsize Ob}\mathcal{T}^{p+1}}\hspace{-0.3cm}
\mathcal{T}(x_1,x_0)\times\mathcal{T}(x_2,x_1)\times\cdots\times\mathcal{T}(x_p,x_{p-1}),
$$
whose face and degeneracy homomorphisms are induced, following the formulas (\ref{p.1.2}), by the composition ${\mathcal{T}(y,x)\times \mathcal{T}(z,y)\overset{\otimes}\to \mathcal{T}(z,x)}$ and unit $\mathrm{I}_x:1\to \mathcal{T}(x,x)$ homomorphisms, respectively. If $a:[q]\to[p]$ is any map in $\Delta$, then one writes $a=d^{i_1}\cdots d^{i_s}s^{j_1}\cdots s^{j_t}$, where
$0\leq i_s<\cdots< i_1 \leq p$, $0\leq j_1<\cdots< j_t \leq q$, and the induced homomorphism is $a^*=s_{j_t}\cdots
s_{j_1}d_{i_s}\cdots d_{i_1}:\ner\mathcal{T}_p\to\ner\mathcal{T}_q$. The pseudo equivalences $\chi$ and $\iota$ arise from the associativity and unit constraints of $\mathcal{T}$, while the invertible modifications $\omega$, $\gamma$  and $\delta$ come from the structure modifications $\pi$, $\mu$, $\lambda$ and $\rho$.
 However, to prove that $\ner\mathcal{T}$ is actually a pseudo-simplicial diagram of bicategories is far from obvious and beyond the scope of this paper since a `supercoherence theorem' is needed. Instead, it will be the subject of an upcoming separate publication specially dedicated to the study of classifying spaces of tricategories and monoidal bicategories.
Hence, we shall only treat here an interesting particular instance: the case of {\em braided monoidal categories} \cite{joyal},  which can be regarded as  one-object, one-arrow tricategories \cite[Corollary 8.7]{g-p-s}.

We shall start by reviewing the notion of classifying space for a monoidal category.
A {\em monoidal} (tensor) category $(\m,\otimes)=(\m,\otimes,\text{I},\boldsymbol{a},\boldsymbol{l},\boldsymbol{r})$, \cite{maclane},  can be viewed as a
bicategory \begin{equation}\label{ome}\Omega^{^{-1}}\hspace{-5pt}\m\end{equation}
with only one object, say $*$,  the objects $u$ of $\m$ as 1-cells $u:*\rightarrow *$ and the morphisms of $\m$
as 2-cells. Thus, $\Omega^{^{-1}}\hspace{-5pt}\m(*,*)=\m$, and it is the horizontal composition of morphisms and deformations given by the  tensor functor
${\otimes:\m\times\m\rightarrow\m}$. The identity at the object is $1_*=\text{I}$, the unit object of the monoidal category, and the associativity, left unit and right unit constraints for $\Omega^{^{-1}}\hspace{-5pt}\m$ are precisely those of the monoidal category, that is,  $\boldsymbol{a}$, $\boldsymbol{l}$ and $\boldsymbol{r}$, respectively.

 The pseudo-simplicial nerve (\ref{ps1.1})  of the bicategory
 $\Omega^{^{-1}}\hspace{-5pt}\m$, hereafter denoted
 \begin{equation}\label{psnm}\ner(\m,\otimes):\Delta^{\!^{\mathrm{op}}} \to\cat,\hspace{0.4cm} [p]\mapsto \m^p,\end{equation}
 is exactly the pseudo-simplicial category that the monoidal category  defines by
 the reduced bar construction \cite[Corollary 1.7]{jardine}, whose category of $p$-simplices is the $p$-fold power of the underlying category $\m$,  and whose face and degeneracy functors are induced by the tensor $\m\times \m\overset{\otimes}\to \m$ and unit $\mathrm{I}:1\to \m$  functors, respectively, following the familiar formulas (\ref{p.1.2})
in analogy with those of the nerve of a monoid. $\ner(\m,\otimes)$
is called the {\em pseudo-simplicial nerve of the monoidal category}
and its classifying space $\class\ner(\m,\otimes)$ is {\em the
classifying space of the monoidal category}  (see \cite[\S
3]{jardine}, \cite[Appendix]{hinich}, \cite{b-c2} or \cite{b-f-s-v},
for example), hereafter denoted by
\begin{equation}\label{clmc}\class(\m,\otimes).\end{equation}

Hence, the classifying space of a monoidal category $(\m,\otimes)$ is the same as the classifying space of $\Omega^{^{-1}}\hspace{-5pt}\m$, the one-object bicategory it defines.
The observation, due to Benabou \cite{benabou}, that monoidal categories are essentially the same as bicategories with just one object is known as the {\em delooping principle}, and the bicategory $\Omega^{^{-1}}\hspace{-5pt}\m$ is called the {\em delooping of the category} induced by its monoidal structure \cite[2.10]{k-v}. This term arises from the existence of a natural map \begin{equation}\label{hgc}\class \m\to \Omega\class(\m,\otimes),\end{equation}  where $\class\m$ is the classifying space of the underlying category  and $\Omega\class(\m,\otimes)$ the loop space based at the 0-cell of $\class(\m,\otimes)$, which is up to group completion a homotopy equivalence (see \cite[Propositions 3.5 and 3.8]{jardine} or \cite[Corollary 4]{b-c2}, for example).

 A monoidal functor $F:(\m,\otimes)\to (\m',\otimes)$ amounts precisely to a homomorphism
 $\Omega^{^{-1}}\hspace{-2pt}F:\Omega^{^{-1}}\hspace{-5pt}\m\to \Omega^{^{-1}}\hspace{-5pt}\m'$ between the
 corresponding delooping bicategories and therefore, by Fact \ref{f2} $(i)$,
 it induces a cellular map
 $$\class(F,\otimes):\class(\m,\otimes)\to
 \class(\m',\otimes).$$
 More precisely, $\class(F,\otimes)$ is the induced on classifying spaces by the
 pseudo-simplicial funtor
 $\ner\Omega^{^{-1}}\hspace{-2pt}F$, hereafter denoted by
 $$
 \ner(F,\otimes):\ner(\m,\otimes)\to\ner(\m',\otimes),\hspace{0.5cm}[p]\mapsto F^p:\m^p\to\m'^p,$$
 whose structure natural
 isomorphisms $s_i F_* \cong F_*s_i$
and $d_iF_*\cong F_* d_i$ are those canonically obtained from  the
invertible  structure constraints of the monoidal functor,
 $\widehat{F}:\text{I}\cong F\text{I}$ and
$\widehat{F}: F(\alpha_i)\otimes F(\alpha_{i+1})\cong
F(\alpha_i\otimes \alpha_{i+1})$ (the commutativity of the needed
six coherence diagrams in \cite{jardine} is clear).

Thus, the classifying space construction, $(\m,\otimes)
\mapsto\class(\m,\otimes)$, defines a functor from monoidal
categories to CW-complexes.

\vspace{0.2cm}

We now consider the braided case.   Recall from \cite[Corollary 8.7]{g-p-s} that a
{\em braided monoidal category} $(\m,\otimes,\boldsymbol{c})=(\m,\otimes,\text{I},\boldsymbol{a},\boldsymbol{l},\boldsymbol{r},\boldsymbol{c})$, \cite[Definition 2.1]{joyal},  defines a one-object, one-arrow tricategory. More precisely, following \cite[2.5]{berger}, \cite[4.2]{k-v} and the categorical delooping principle, let
\begin{equation}\label{omenosdos}\Omega^{^{-2}}\hspace{-5pt}\m\end{equation}
denote the tricategory with only one object, say $*$,  only one arrow $*=1_*:*\to *$, the objects $u$ of $\m$ as 2-cells $u:*\rightarrow *$ and the morphisms of $\m$
as 3-cells. Thus, $\Omega^{^{-2}}\hspace{-5pt}\m(*,*)=\Omega^{^{-1}}\hspace{-5pt}\m$, the delooping bicategory associated to the underlying monoidal category (\ref{ome}), the composition is also (as the horizontal one in $\Omega^{^{-1}}\hspace{-5pt}\m$) given by the  tensor functor
${\otimes:\m\times\m\rightarrow\m}$ and the interchange 3-cell between the two different composites of 2-cells is given by the braiding $\boldsymbol{c}:u\otimes v\to v\otimes u$.

Call this tricategory $\Omega^{^{-2}}\hspace{-5pt}\m$ the {\em
double delooping} of the underlying category  $\m$ associated to the
given braided monoidal structure on it, and call its corresponding
bicategorical pseudo-simplicial nerve (\ref{psnt}) the {\em
pseudo-simplicial nerve of the braided monoidal category}, hereafter
denoted by $\ner(\m,\otimes,\boldsymbol{c})$. Thus, it is given by
\begin{equation}\label{psnbmc}\ner(\m,\otimes,\boldsymbol{c}):\Delta^{\!^{\mathrm{op}}}
\to\bicat, \hspace{0.4cm} [p]\mapsto
 (\Omega^{^{-1}}\hspace{-5pt}\m)^p,\end{equation}
and next we see that $\ner(\m,\otimes,\boldsymbol{c})$ is actually a
pseudo-simplicial bicategory.

Because of the  braiding, the
pseudo-simplicial nerve  of the monoidal category,
$\ner(\m,\otimes):[p]\mapsto \m^p$, is actually the underlying
pseudo-simplicial category of the
pseudo-simplicial monoidal category, $$[p]\mapsto (\m^p,\otimes)=
(\m,\otimes)^p.$$ Indeed, this follows because the functors $a^*:(\m^q,\otimes)\to(
\m^p,\otimes)$ and the structure natural isomorphisms
$\chi:b^*a^*\cong (ab)^*$ are monoidal  (it suffices to observe the
monoidal structure for the face and degeneracy functors
(\ref{p.1.2})  and also for the natural isomorphisms (\ref{ps1.2}) and (\ref{ps1.3}), which can be respectively deduced from Propositions 5.2 and 5.1 in   \cite{joyal}).

Then we have that $\ner(\m,\otimes,\boldsymbol{c})$ is just the pseudo-simplicial bicategory obtained
 as the composite
$$
\xymatrix@R=0pt{\hspace{-0.9cm}&\Delta^{\!^{\mathrm{op}}}\ar[rrr]^-{\textstyle
(\ner(\m,\otimes),\otimes)}&&&\mcat\ar[r]^-{\textstyle \Omega^{^{-1}}}&\bicat\,.\\
&[p]\ar@{|->}[rrr]&&&
(\m^p,\otimes)\ar@{|->}[r]&\Omega^{^{-1}}\hspace{-3pt}\m^p}
$$

  Hence,
$\ner(\m,\otimes,\boldsymbol{c})$ is actually a pseudo-simplicial
diagram of one-object bicategories (with the structure modifications
$\omega$, $\gamma$, and $\delta$ all being identities) and,
following the general Definition \ref{defclass}, we give the following:

\begin{definition}\label{dcbmc} The classifying space  of the braided monoidal category, denoted by
\begin{equation}\label{cbmc}\class(\m,\otimes,\boldsymbol{c}),\end{equation}
  is  defined to be the classifying space of its pseudo-simplicial nerve $(\ref{psnbmc})$.
\end{definition}
\begin{remark}  {\em By replacing each delooping bicategory $\Omega^{^{-1}}\hspace{-5pt}\m^p$ by  its
pseudo simplicial nerve (\ref{ps1.1}), that is, by the nerve
(\ref{psnm}) of the monoidal category  $(\m^p,\otimes)$,  the
pseudo-simplicial nerve of the braided monoidal category
(\ref{psnbmc})  determines a  pseudo bisimplicial category
$\Delta^{\!^{\mathrm{op}}}\times \Delta^{\!^{\mathrm{op}}}\to \cat$,
$([p],[q])\mapsto \m^{pq}$, which (for $(\m,\otimes)$ strict) is
taken in \cite{b-f-s-v} to construct the double delooping space of
$\class\m$.}
\end{remark}

The basic properties of the classifying space construction for braided monoidal categories can be stated as follows:

\begin{proposition} $(i)$ Any braided monoidal functor between braided monoidal categories,
$F:(\m,\otimes,\boldsymbol{c})\to (\m',\otimes,\boldsymbol{c})$,
induces a continuous map between the corresponding classifying
spaces, $$\class
(F,\otimes,\boldsymbol{c}):\class(\m,\otimes,\boldsymbol{c})\to\class(\m',\otimes,\boldsymbol{c}).$$
Therefore, the classifying space construction,
$(\m,\otimes,\boldsymbol{c})\mapsto
\class(\m,\otimes,\boldsymbol{c})$, defines a functor from the
category of braided monoidal categories to CW-complexes.

\vspace{0.2cm} \noindent $(ii)$ If two braided monoidal functors
$F,F':(\m,\otimes,\boldsymbol{c})\to (\m',\otimes,\boldsymbol{c})$
are  related by a monoidal transformation $F\Rightarrow F'$, then
the induced maps on classifying spaces, $\class
(F,\otimes,\boldsymbol{c})$ and $\class(F',\otimes,\boldsymbol{c})$,
are homotopic.

\vspace{0.2cm} \noindent (iii) If
$F:(\m,\otimes,\boldsymbol{c})\overset{\sim}\to
(\m',\otimes,\boldsymbol{c})$ is a braided  monoidal equivalence,
then the induced map on classifying spaces $\class
(F,\otimes,\boldsymbol{c}):\class(\m,\otimes,\boldsymbol{c})\overset{\sim}\to\class(\m',\otimes,\boldsymbol{c})$
is a homotopy equivalence.

\vspace{0.2cm}
\noindent (iv) If $F:(\m,\otimes,\boldsymbol{c})\to  (\m',\otimes,\boldsymbol{c})$ is a braided monoidal functor such that the underlying functor induces a homotopy equivalence $\class F:\class \m\overset{\sim}\to \class\m'$, then the induced map $\class (F,\otimes,\boldsymbol{c}):\class(\m,\otimes,\boldsymbol{c})\overset{\sim}\to\class(\m',\otimes,\boldsymbol{c})$ is a homotopy equivalence (as is also the induced map  between the  classifying spaces of the underlying monoidal categories, $\class (F,\otimes):\class(\m,\otimes)\overset{\sim}\to\class(\m',\otimes)$).
\end{proposition}

\begin{proof}
$(i)$ If
$F:(\m,\otimes,\boldsymbol{c})\to(\m',\otimes,\boldsymbol{c})$ is
any braided monoidal functor, then the pseudo-simplicial functor
$\ner (F,\otimes):\ner (\m,\otimes)\to\ner (\m',\otimes)$ underlies
a  pseudo-simplicial monoidal functor $$\ner (F,\otimes):(\ner
(\m,\otimes),\otimes)\to (\ner (\m',\otimes),\otimes);$$ that is,
every functor $F^p:(\m^p,\otimes)\to (\m'^p,\otimes)$ is monoidal
 and, moreover, every natural isomorphism ${ F^p a^*\cong
a^* F^q}$, for any  $a:[q]\to [p]$  in $\Delta$, is monoidal (it
suffices to prove this for the natural isomorphisms $s_i F^p \cong
F^{p+1}s_i$ and $d_iF^p\cong F^{p-1} d_i$,  which it is
straightforward).

Hence, we have a pseudo-simplicial homomorphism of pseudo-simplicial
bicategories $\Omega^{^{-1}}\hspace{-2pt}\ner (F,\otimes)$ (with the
structure modifications $\Pi$ and $\Gamma$ all being identities),
hereafter denoted by
$$\ner(F,\otimes,\boldsymbol{c}):\ner(\m,\otimes,\boldsymbol{c})\to
\ner(\m',\otimes,\boldsymbol{c}),$$ which, by
Proposition \ref{propclass} $(i)$, gives the claimed cellular map
$$\class
(F,\otimes,\boldsymbol{c}):\class(\m,\otimes,\boldsymbol{c})\to\class(\m',\otimes,\boldsymbol{c}).
$$

 Following now the proof of part $(iii)$ in Proposition \ref{propclass}, we see that the classifying  space
 construction defines a functor from the category of braided monoidal categories to the category of spaces:
 for $(\m,\otimes,\boldsymbol{c})\overset{F}\to
 (\m',\otimes,\boldsymbol{c})\overset{G}\to (\m'',\otimes,\boldsymbol{c})$, any two composable
 braided monoidal functors, the equality $$\xymatrix{\int_\Delta\!\! \ner(G,\otimes,\boldsymbol{c})\,
\int_\Delta\!\! \ner(F,\otimes,\boldsymbol{c})=\int_\Delta
\!\!\ner(GF,\otimes,\boldsymbol{c})}$$ holds (and the corresponding
pseudo-natural equivalence (\ref{gc3}) is an identity), whence the
equality $\class(G,\otimes,\boldsymbol{c})
\class(F,\otimes,\boldsymbol{c})=\class(GF,\otimes,\boldsymbol{c})$
follows from Fact \ref{f2} $(i)$. Analogously, the equality $\class
1_{(\m,\otimes,\boldsymbol{c})}=
1_{\class(\m,\otimes,\boldsymbol{c})}$ holds since the
pseudo-natural equivalence (\ref{gc4}) at any
$\ner(\m,\otimes,\boldsymbol{c})$ is an identity.

\vspace{0.1cm} $(ii)$  Any monoidal transformation,  $m:F\Rightarrow
F'$, between monoidal functors $F,F':(\m,\otimes)\to (\m',\otimes)$,
gives rise to a pseudo-simplicial  transformation $$\ner(m,\otimes):\ner
(F,\otimes)\Rightarrow \ner(F',\otimes): \ner (\m,\otimes)\to
\ner(\m',\otimes),$$ where $$\ner_p(m,\otimes)=m^p:F^p\Rightarrow F'^p:\m^p\to
\m'^p$$  (see \cite[p.125]{jardine}). When both $F$ and $F'$ are
braided between braided monoidal categories
$(\m,\otimes,\boldsymbol{c})$ and $(\m',\otimes,\boldsymbol{c})$,
then every $m^p:F^p\Rightarrow F'^p: (\m^p,\otimes)\to
(\m'^p,\otimes)$ is monoidal and $\ner(m,\otimes)$ becomes a pseudo-simplicial
monoidal transformation giving rise to a pseudo-transformation of
pseudo-simplicial homomorphisms of bicategories
$$\Omega^{^{-1}}\hspace{-2pt}\ner(m,\otimes):\ner(F,\otimes,\boldsymbol{c})\Rightarrow\ner(F',\otimes,\boldsymbol{c}):\ner(\m,\otimes,\boldsymbol{c})
\to\ner(\m',\otimes,\boldsymbol{c}),$$ whence the result follows
from part $(ii)$ of Proposition \ref{propclass}.

\vspace{0.1cm}
$(iii)$ It is a consequence of parts $(i)$ and $(ii)$ (and also of part $(iv)$) .

\vspace{0.1cm} $(iv)$ Since, for any $p\geq 0$, the induced map
$\class F^p:\class \m^p\to\class\m'^p$ is a homotopy equivalence,
Thomason's theorem \cite[Corollary 3.3.1]{thomason} means that the
pseudo-simplicial functor $\ner(F,\otimes):\ner(\m,\otimes)\to
\ner(\m',\otimes)$ induces a homotopy equivalence on classifying
spaces,
$\class(F,\otimes):\class(\m,\otimes)\overset{\sim}\to\class(\m',\otimes)$.
Then, each map
${\class(F^p,\otimes):\class(\m^p,\otimes)\overset{\sim}\to\class(\m'^p,\otimes)}$
is also a homotopy equivalence whence, by Theorem \ref{p.1.23}, the
pseudo-simplicial homomorphism of bicategories
$\ner(F,\otimes,\boldsymbol{c})$ induces a homotopy equivalence
$\class(F,\otimes,\boldsymbol{c}):
\class(\m,\otimes,\boldsymbol{c})\overset{\scriptsize
\sim}\to\class(\m',\otimes,\boldsymbol{c})$, as claimed.
\end{proof}

Returning to the monoidal case, if $(\m,\otimes)$ is any given
monoidal category, then the delooping bicategory
$\Omega^{^{-1}}\hspace{-4pt}\m$ has a corresponding unitary
geometric nerve (\ref{ngn}),
$\Delta^{\hspace{-2pt}^\mathrm{u}}\!\Omega^{^{-1}}\hspace{-4pt}\m$.
But, hereafter, we shall follow the terminology of \cite[\S 4]{cegarra3} and \cite[Definition 4.1]{cg},  where  a 2-{\em cocycle} of a (small) category $I$ in
the monoidal category $(\m,\otimes)$ is  defined as a normal lax
functor $I\to \Omega^{^{-1}}\hspace{-5pt}\m$. Therefore, such a
2-cocycle is a system of data $$\xi:I\to (\m,\otimes) $$ consisting
of an object $\xi_\sigma\in \m$ for each arrow $\sigma:j\to i$ in
$I$ and of a morphism $\xi_{\sigma,\tau}:\xi_\sigma\otimes
\xi_\tau\to\xi_{\sigma\tau}$ for each pair of composible arrows in
$I$, $k\overset{\tau}\to j\overset{\sigma}\rightarrow i$, such that,
for any three composable arrows in $I$, $l\overset{\gamma}\to k
\overset{\tau}\to j\overset{\sigma}\to i$, the diagram in $\m$
$$\xymatrix{\xi_\sigma\otimes(\xi_\tau\otimes \xi_\gamma)\ar[rr]^-{\textstyle \boldsymbol{a}}\ar[d]_-{\textstyle 1\otimes
\xi_{\tau,\gamma}}&&(\xi_\sigma\otimes \xi_\tau)\otimes \xi_\gamma\ar[d]^-{\textstyle \xi_{\sigma,\tau}\otimes 1}\\
\xi_\sigma\otimes \xi_{\tau\gamma}\ar[r]^-{\textstyle \xi_{\sigma,\tau\gamma}}&\xi_{\sigma\tau\gamma}&
\xi_{\sigma\tau}\otimes \xi_{\gamma}\ar[l]_-{\textstyle \xi_{\sigma\tau,\gamma}}
 }
$$
is commutative,  $\xi_{1} =I$,  $\xi_{1,\sigma}=\boldsymbol{l}:I\otimes \xi_{\sigma}\to \xi_\sigma$, and $\xi_{\sigma,1}=
\boldsymbol{r}:\xi_\sigma\otimes I\to \xi_\sigma$.

\vspace{0.2cm}

These 2-cocycles of $I$ in $(\m,\otimes)$ form the set, denoted by
$$Z^2(I,(\m,\otimes)),$$ and they are the objects of a category
\begin{equation}\label{c2co}\cdco(I,(\m,\otimes)),\end{equation} where a morphism $f :\xi\to
\xi'$ consists of a family of morphisms $f_\sigma:\xi_\sigma\to
\xi'_\sigma$ in $\m$, one for each arrow $\sigma:j\to i$ in $I$,
such that $f_{1}=1_\mathrm{I}$ and for any two arrows
$k\overset{\tau}\to j\overset{\sigma}\rightarrow i$ the following
square commutes:
$$
\xymatrix{\xi_\sigma\otimes \xi_\tau\ar[r]^{\textstyle \xi_{\sigma,\tau}} \ar[d]_{\textstyle
f_{\sigma}\otimes f_\tau}&\xi_{\sigma\tau}\ar[d]^{\textstyle f_{\sigma\tau}}\\
\xi'_\sigma\otimes \xi'_\tau\ar[r]^{\textstyle \xi'_{\sigma,\tau}}&\xi'_{\sigma\tau}.}
$$

We should note that the category $\cdco(I,(\m,\otimes))$ is a
subbicategory of the bicategory ${\bf
Lax}(I,\Omega^{^{-1}}\hspace{-3pt}\m)$, defined in \cite[p.
569]{street}. Namely, that subbicategory  given by the normal lax functors and
those lax transformations and modifications whose components at any
objects are identities.

The {\em geometric nerve of the monoidal category} $(\m,\otimes)$, \cite{b-c2}, is
then the simplicial set ($\cong
\Delta^{\hspace{-2pt}^\mathrm{u}}\!\Omega^{^{-1}}\hspace{-4pt}\m$)
 \begin{equation}\label{gnmc}Z^2(\m,\otimes):
 \Delta^{\!^{\mathrm{op}}}\ \to \ \set, \hspace{0.5cm}[p]
 \mapsto Z^2([p],(\m,\otimes)).\end{equation}
And this is the simplicial set of objects of the {\em categorical
geometric nerve of the monoidal category}, that is, the simplicial
category
\begin{equation}\label{cgnm}\cdco(\m,\otimes):\Delta^{\!^{\mathrm{op}}}\ \to \ \cat,
\hspace{0.5cm}[p]\mapsto \cdco([p],(\m,\otimes)).\end{equation}

 This geometric nerve $Z^2(\m,\otimes)$ is a $3$-coskeletal reduced (1-vertex) simplicial  set
  whose simplices have the following
simplified interpretation: the $1$-simplices  are the objects
$\xi_{0,1}$ of $\m$ and, for $p\geq 2$, the  $p$-simplices are
families of morphisms $$\xi_{i,j,k}:\xi_{i,j}\otimes \xi_{j,k}\to
\xi_{i,k} ,$$ $0\leq i<j<k\leq p$, making commutative the diagrams
$$
\xymatrix{ (\xi_{i,j}\otimes \xi_{j,k})\otimes
\xi_{k,l}\ar[rr]^{\textstyle  \boldsymbol{a}}\ar[d]_{\textstyle
\xi_{i,j,k}\otimes 1}&&\xi_{i,j}\otimes (\xi_{j,k}\otimes
\xi_{k,l})\ar[d]^{\textstyle 1\otimes \xi_{j,k,l}}\\\xi_{i,k}\otimes
\xi_{k,l}\ar[r]^-{\textstyle \xi_{i,k,l}}&\xi_{i,l}&\xi_{i,j}\otimes
\xi_{j,l}\ar[l]_{\textstyle \xi_{i,j,l}} }
$$
for $0\leq i<j<k<l\leq p$.

There is a pseudo-simplicial functor \cite[p. 325]{b-c2}
\begin{equation}\label{E} (~)^\mathrm{e}:\ner(\m,\otimes)\to\cdco(\m,\otimes), \end{equation}
taking an object  $X=(X_1,\dots,X_p)\in \ner_p(\m,\otimes)= \m^p$ to
the 2-cocycle \begin{equation}\label{ef}X^\mathrm{e}:[p]\to
(\m,\otimes),\end{equation}
 with $X_{i,i+1}^\mathrm{e}=X_{i+1}$ and, inductively, $X_{i,j+1}^\mathrm{e}=X_{i,j}^\mathrm{e}\otimes X_{j+1}$. The morphisms $X_{i,j,j+1}^\mathrm{e}:X_{i,j}^\mathrm{e}\otimes X_{j,j+1}^\mathrm{e}\to X_{i,j+1}^\mathrm{e}$ are all identities, and the remaining morphisms  $X_{i,j,k+1}^\mathrm{e}:X_{i,j}^\mathrm{e}\otimes X_{j,k+1}^\mathrm{e}\to X_{i,k+1}^\mathrm{e}$ are inductively  determined by the associativity constraints
of $\m$, through the commutative diagrams

$$
\xymatrix@C=45pt{ (X_{i,j}^\mathrm{e}\otimes X_{j,k}^\mathrm{e})\otimes X_{k,k+1}^\mathrm{e} \ar[d]_{\textstyle \boldsymbol{a}}\ar[r]^-{\textstyle X_{i,j,k}^\mathrm{e}\otimes 1} & X_{i,k}^\mathrm{e}\otimes X_{k,k+1}^\mathrm{e}\\X_{i,j}^\mathrm{e}\otimes X_{j,k+1}^\mathrm{e} \ar[r]^{\textstyle X_{i,j,k+1}^\mathrm{e}} & X_{i,k+1}^\mathrm{e}\,.\ar@{=}[u]}
$$
Further, the functor $(~)^\mathrm{e}$ on a morphism
$F=(F_1,\dots,F_p):X\to Y$ in $\m^p$ is the 2-cocycle morphism
$F^\mathrm{e}: X^\mathrm{e}\to Y^\mathrm{e}$,  inductively given by
$$F_{i,j+1}^\mathrm{e}=\left\{
\begin{array}{lll}
   F_{i+1} & \text{if} & j=i, \\
    &&\\
F_{i,j}^\mathrm{e}\otimes F_{j+1} & \text{if} & j>i. \\
\end{array}
\right.$$
For any map $a:[q]\to [p]$ in the simplicial category, the natural isomorphisms $ (a^*X)^\mathrm{e}\cong a^*(X^\mathrm{e})$ are canonically induced by the associativity and unit constraints $\boldsymbol{a}$, $\boldsymbol{l}$, and $\boldsymbol{r}$ of the monoidal category.

The main purpose in \cite{b-c2} (cf. Fact  \ref{f1})  was to prove the following:
\begin{fact}\label{fm}
For any monoidal category $(\m,\otimes)$, both
$(~)^\mathrm{e}:\ner(\m,\otimes)\to \cdco(\m,\otimes)$ and the
inclusion $Z^2(\m,\otimes)\to
\cdco(\m,\otimes)$ induce homotopy equivalences on classifying
spaces. In particular, there is a  homotopy equivalence
$$\class(\m,\otimes)\simeq  |Z^2(\m,\otimes)|.$$
\end{fact}

Going further  towards the  braided case, we shall start with the following observation:
\begin{lemma}\label{lmon2} Let $(\m,\otimes,\boldsymbol{c})$ be a braided monoidal category.

\noindent$(i)$ For any small category $I$, the category of
$2$-cocycles $\cdco(I,(\m,\otimes))$, $(\ref{c2co})$, has a natural
monoidal structure. The tensor product $\xi'\!\otimes \xi$ of
$2$-cocycles is given by putting  $(\xi'\!\otimes
\xi)_\sigma=\xi'_\sigma\otimes \xi_\sigma$, and $(\xi'\!\otimes
\xi)_{\sigma,\tau}$ is the composite dotted arrow in the diagram
$$ \xymatrix@C=-5pt@R=25pt{(\xi'_\sigma\otimes \xi_\sigma)\otimes (\xi'_\tau \otimes \xi_\tau)
\ar@{.>}[rr]^-{\textstyle (\xi'\otimes \xi)_{\sigma.\tau}}\ar[d]_{\textstyle
\cong }& &\xi'_{\sigma\tau}\otimes \xi_{\sigma\tau}\\
(\xi'_\sigma\otimes ( \xi_\sigma\otimes \xi'_\tau))\otimes \xi_\tau \ar[dr]
\ar@{}@<-3pt>[dr]_(0.4){\textstyle (1\otimes \boldsymbol{c})\otimes 1} & &(\xi'_\sigma\otimes \xi'_\tau)\otimes
( \xi_\sigma\otimes \xi_\tau) \ar[u]_{\textstyle \xi'_{\sigma,\tau}\otimes\xi_{\sigma,\tau}}\\
&(\xi'_\sigma\otimes(\xi'_\tau\otimes \xi_\sigma))\otimes
\xi_\tau\ar[ru]_{\textstyle \cong} &}$$ where the arrows labeled
with $\cong$ are (iterated) isomorphisms of associativity. The
tensor product  of morphisms $f'$ and $f$ is $f'\!\otimes f$, where
$(f'\otimes f)_\sigma= f'_\sigma\otimes f_\sigma$. The unit object
is the {\em trivial} $2$-cocycle,  denoted by $\mathrm{I}_0,$ which
is defined by the equalities $(\mathrm{I}_0)_{\sigma}=\mathrm{I}$
and
$({\mathrm{I}}_0)_{\sigma,\tau}=\boldsymbol{l}=\boldsymbol{r}:\mathrm{I}\otimes
\mathrm{I}\to \mathrm{I}$. The associativity and identity
constraints of $(\m,\otimes)$ yield associativity and identity
constraints in $\cdco(I,(\m,\otimes))$.

\vspace{0.2cm} \noindent$(ii)$ The categorical geometric nerve of
the underlying monoidal category  $(\ref{cgnm})$ underlies  the
simplicial monoidal category
$$(\cdco(\m,\otimes),\otimes):\Delta^{\!^{\mathrm{op}}}\ \to
\mathbf{MonCat},\hspace{0.4cm}[p]\mapsto
(\cdco([p],(\m,\otimes)),\otimes).$$

\vspace{0.2cm}
 \noindent$(iii)$ The pseudo-simplicial functor $(\ref{E})$, is actually a pseudo-simplicial
monoidal functor $$(~)^\mathrm{e}:(\ner(\m,\otimes),\otimes)\to
(\cdco(\m,\otimes),\otimes).$$If $Y=(Y_1,\dots,Y_p)$ and
$X=(X_1,\dots,X_p)$ are in $\m^p$, then the structure isomorphism
$\Phi:Y^\mathrm{e}\!\otimes X^\mathrm{e}\to (Y\!\otimes
X)^\mathrm{e}$ is as follows: $\Phi_{i,i+1}=1:Y_{i+1}\otimes
X_{i+1}\to Y_{i+1}\otimes X_{i+1}$ and, for $0\leq i<j<p$,
$\Phi_{i,j+1}$ is inductively defined as the composite dotted
arrow
$$
\xymatrix@C=-40pt@R=25pt{(Y_{\hspace{-2pt}i,j}^\mathrm{e}\otimes
Y_{\hspace{-1pt}j+1})\otimes (X_{i,j}^\mathrm{e}\otimes
X_{j+1})\ar@{.>}[rr]^-{\textstyle \Phi_{i,j+1}}\ar[d]_{\textstyle
\cong }& &(Y\!\otimes X)_{i,j}^\mathrm{e}\otimes (Y_{\hspace{-1pt}j+1}\otimes X_{j+1})\\
(Y_{\hspace{-1pt}i,j}^\mathrm{e}\otimes (Y_{\hspace{-1pt}j+1}\otimes
X_{i,j}^\mathrm{e}))\otimes X_{j+1}
     \ar[dr]\ar@{}@<-3pt>[dr]_(0.4){\textstyle (1\otimes \boldsymbol{c})\otimes 1} & &
(Y_{\hspace{-1pt}i,j}^\mathrm{e}\otimes X_{i,j}^\mathrm{e})\otimes
(Y_{\hspace{-1pt}j+1}\otimes X_{j+1})
 \ar[u]_{\textstyle \textstyle \Phi_{i,j}\otimes 1}\\
&(Y_{\hspace{-1pt}i,j}^\mathrm{e}\otimes( X_{i,j}^\mathrm{e}\otimes
Y_{j+1}))\otimes X_{j+1} \ar[ru]_{\textstyle \cong} &}$$ The
structure isomorphism $(\mathrm{I},\dots,\mathrm{I})^\mathrm{e}\to
\mathrm{I}_0$ is given by the canonical isomorphism in $\m$,
$(\cdots(\cdots\otimes \mathrm{I})\otimes \mathrm{I})\otimes
\mathrm{I}\cong \mathrm{I}$.

\vspace{0.2cm}\noindent $(iv)$ There is an induced pseudo-simplicial homomorphism
\begin{equation}\label{ipsh}
\Omega^{^{-1}}\hspace{-3pt}(~)^\mathrm{e}:\ner(\m,\otimes,\boldsymbol{c})
\to \Omega^{^{-1}}\hspace{-2pt}\cdco(\m,\otimes).
\end{equation}
\end{lemma}

Next, again following \cite[\S 4]{cegarra3},  where a 3-{\em cocycle} of a  category $I$ in a braided monoidal category $(\m,\otimes,\boldsymbol{c})$ is defined to be a normal lax functor $I\to \Omega^{^{-2}}\hspace{-5pt}\m$, \cite[Definition 3.1]{g-p-s}, we establish the following:

\begin{definition}\label{3-cocy}  Let $(\m,\otimes,\boldsymbol{c})$ be a braided monoidal category. For any given small category $I$, a {\em $3$-cocycle}  $$\lambda: I\to (\m,\otimes,\boldsymbol{c})$$ is  a system of data consisting  of:

- for each two composible arrows in $I$,  $k\overset{\tau}\to
j\overset{\sigma}\rightarrow i$, an object $\lambda_{\sigma,\tau} \in \m$,

- for each triplet of composible arrows in $I$, $l\overset{\gamma}\to k\overset{\tau}\to
j\overset{\sigma}\rightarrow i$, a morphism in $\m$
$$\xymatrix@C=35pt{\lambda_{\tau,\gamma}\otimes \lambda_{\sigma,\tau\gamma}\ar[r]^-{\textstyle \lambda_{\sigma,\tau,\gamma}}&\lambda_{\sigma,\tau}\otimes \lambda_{\sigma\tau,\gamma}},$$
such that, for any four composible arrows in $I$, $m\overset{\delta}\to l\overset{\gamma}\to k
\overset{\tau}\to j\overset{\sigma}\to i$, the following diagram in $\m$ commutes
$$\xymatrix@C=90pt@R=36pt{(\lambda_{\gamma,\delta}\otimes \lambda_{\tau,\gamma\delta})\otimes \lambda_{\sigma,\tau\gamma\delta}\ar[d]_{\textstyle \lambda_{\tau,\gamma,\delta}\otimes 1}\ar[r]^{\textstyle \boldsymbol{a}(1\otimes \lambda_{\sigma,\tau,\gamma\delta})\boldsymbol{a}^{-1}}&
(\lambda_{\gamma,\delta}\otimes \lambda_{\sigma,\tau})\otimes \lambda_{\sigma\tau,\gamma\delta}\ar[d]^{\textstyle \boldsymbol{c}\otimes 1}\\
(\lambda_{\tau,\gamma}\otimes \lambda_{\tau\gamma,\delta})\otimes \lambda_{\sigma,\tau\gamma\delta}\ar[d]_{\textstyle (1\otimes \lambda_{\sigma,\tau\gamma,\delta})\boldsymbol{a}}& (\lambda_{\sigma,\tau}\otimes \lambda_{\gamma,\delta})\otimes \lambda_{\sigma,\tau\gamma\delta}\ar[d]^{\textstyle (1\otimes \lambda_{\sigma\tau,\gamma,\delta})\boldsymbol{a}}\\
\lambda_{\tau,\gamma}\otimes (\lambda_{\sigma,\tau\gamma}\otimes \lambda_{\sigma\tau\gamma,\delta})\ar[r]^{\textstyle \boldsymbol{a}^{-1}(\lambda_{\sigma,\tau,\gamma}\otimes 1)\boldsymbol{a}}& \lambda_{\sigma,\tau}\otimes (\lambda_{\sigma\tau,\gamma}\otimes \lambda_{\sigma\tau\gamma,\delta})
 }
$$
and, moreover, the following equalities hold:  $\lambda_{1,\sigma} =I= \lambda_{\sigma,1}$,  $\lambda_{1,\sigma,\tau}=\boldsymbol{c}_{_{I,\lambda_{\sigma,\tau}}}$, $\lambda_{\sigma,1,\tau}=1$ and $\lambda_{\sigma,\tau,1}=\boldsymbol{c}_{_{\lambda_{\sigma,\tau},I}}$.
\end{definition}
 The  $3$-cocycles of $I$ in the braided monoidal
category $(\m,\otimes,\boldsymbol{c})$ form the set,   denoted by
$$ Z^3(I,(\m,\otimes,\boldsymbol{c})),$$
 which is the set of objects of
a bicategory
$$\ctco(I,(\m,\otimes,\boldsymbol{c})),$$
\underline{whose 1-cells}
$\xi:\lambda\to \lambda'$ consist  of pairs of maps assigning

\vspace{0.2cm}
- to each arrow $\sigma:j\to i$ in $I$, an object $\xi_\sigma\in \m$,

\vspace{0.2cm}
- to each pair of composible arrows in $I$, $k\overset{\tau}\to
j\overset{\sigma}\rightarrow i$, a morphism in $\m$
$$\xymatrix{(\xi_\sigma\otimes \xi_\tau)\otimes \lambda_{\sigma,\tau}\ar[r]^-{\textstyle \xi_{\sigma,\tau}}
&\lambda'_{\sigma,\tau}\otimes \xi_{\sigma\tau},}$$ such that, for any three composible arrows in $I$, $l\overset{\gamma}\to k
\overset{\tau}\to j\overset{\sigma}\to i$, the diagram below (where we have omitted the associativity constraints) is
commutative
$$\xymatrix@C=60pt@R=30pt{\xi_{\sigma}\otimes \xi_{\tau}\otimes \xi_{\gamma}\otimes \lambda_{\tau,\gamma}\otimes \lambda_{\sigma,\tau\gamma}
\ar[r]^{\textstyle 1\otimes \xi_{\tau,\gamma}\otimes 1}\ar[d]_-{\textstyle\  \  1\otimes
\lambda_{\sigma,\tau,\gamma}}&\xi_\sigma\otimes \lambda'_{\tau,\gamma}\otimes
\xi_{\tau\gamma}\otimes \lambda_{\sigma,\tau\gamma}\ar[d]^-{\textstyle \boldsymbol{c}\otimes 1}\\
\xi_{\sigma}\otimes \xi_{\tau}\otimes \xi _{\gamma}\otimes\otimes \lambda_{\sigma,\tau}\otimes
\lambda_{\sigma\tau,\gamma}\ar[d]_-{\textstyle 1\otimes \boldsymbol{c}\otimes 1} &\lambda'_{\tau,\gamma}\otimes \xi_\sigma\otimes
\xi_{\tau\gamma}\otimes \lambda_{\sigma,\tau\gamma}\ar[d]^-
{\textstyle \ 1\otimes \xi_{\sigma,\tau\gamma}}\\
\xi_\sigma\otimes \xi_\tau\otimes \lambda_{\sigma,\tau}\otimes \xi_\gamma\otimes \lambda_{\sigma\tau,\gamma}\ar[d]_-{\textstyle
\xi_{\sigma,\tau}\otimes 1} &\lambda'_{\tau,\gamma}\otimes
\lambda'_{\sigma,\tau\gamma}\otimes \xi_{\sigma\tau\gamma}\ar[d]^-{\textstyle \lambda'_{\sigma,\tau,\gamma}\otimes 1}\\
\lambda'_{\sigma,\tau}\otimes \xi_{\sigma\tau}\otimes \xi_\gamma\otimes \lambda_{\sigma\tau,\gamma}\ar[r]^-{\textstyle 1\otimes
\xi_{\sigma\tau,\gamma}}&\lambda'_{\sigma,\tau}\otimes \lambda'_{\sigma\tau,\gamma}\otimes \xi_{\sigma\tau\gamma},
 }
$$
\noindent moreover,  $\xi_{1_k}=I$ and, for every arrow $\tau:k\to l$, the  squares below commute.
$$
\xymatrix@C=30pt{(I\otimes \xi_\tau)\otimes I\ar[r]^(.6){\textstyle \xi_{1,\tau}}\ar[d]_{\textstyle
\boldsymbol{l}\otimes 1} & I\otimes \xi_\tau\ar[d]_{\textstyle \boldsymbol{l}} & \ar[l]_(.5){\textstyle
\xi_{\tau,1}}(\xi_\tau\otimes I)\otimes I\ar[d]^{\textstyle \boldsymbol{r}\otimes 1}\\ \xi_\tau\otimes
I\ar[r]^{\textstyle \boldsymbol{r}}&\xi_\tau& \xi_\tau\otimes I\ar[l]_{\textstyle \boldsymbol{r}}}
$$

\underline{A $2$-cell } $f:\xi\Rightarrow \xi'$, for  $\xi,\xi':\lambda\to\lambda'$  1-cells, consists of a family of morphisms
$f_\sigma:\xi_\sigma\to \xi'_\sigma$ in $\m$, one for each arrow $\sigma:j\to i$ in $I$, such that $f_{1}=1_\mathrm{I}$ and
for any two arrows $k\overset{\tau}\to j\overset{\sigma}\rightarrow i$ the following square commutes:
$$
\xymatrix{(\xi_\sigma\otimes \xi_\tau)\otimes \lambda_{\sigma,\tau}\ar[r]^{\textstyle \xi_{\sigma,\tau}}
\ar[d]_{\textstyle
(f_{\sigma}\otimes f_\tau)\otimes 1}&\lambda'_{\sigma,\tau}\otimes \xi_{\sigma\tau}\ar[d]^{\textstyle 1\otimes f_{\sigma\tau}}\\
(\xi'_\sigma\otimes \xi'_\tau)\otimes \lambda_{\sigma,\tau}\ar[r]^{\textstyle \xi'_{\sigma,\tau}}&\lambda'_{\sigma,\tau}\otimes
\xi'_{\sigma\tau}.}
$$

\underline{The vertical composition} of $2$-cells in
$\ctco(I,(\m,\otimes,\boldsymbol{c}))$ is defined by pointwise
composition in $\m$.

\underline{The horizontal composition} of 1-cells  $\xi:\lambda\to \lambda'$ and
$\xi':\lambda'\to \lambda''$ is  ${\xi'\!\otimes\! \xi:\lambda\to \lambda''}$, where
$(\xi'\!\otimes \xi)_\sigma={\xi'}_{\hspace{-3pt}\sigma}\otimes
\xi_\sigma $ and $(\xi'\!\otimes \xi)_{\sigma,\tau}$ is the
composite dotted arrow in the diagram
$$\xymatrix@C=40pt{(({\xi'}_{\hspace{-3pt}\sigma}\otimes \xi_\sigma)\otimes ({\xi'}_{\hspace{-3pt}\tau} \otimes \xi_\tau))\otimes \lambda_{\sigma,\tau}\ar[d]_{\textstyle
\cong }\ar@{.>}[r]^-{\textstyle (\xi'\!\otimes \xi)_{\sigma,\tau}} & \lambda''_{\sigma,\tau}\otimes
(\xi'_{\sigma\tau}\otimes \xi_{\sigma\tau})\\
(({\xi'}_{\hspace{-3pt}\sigma}\otimes ( \xi_\sigma\otimes \xi'_\tau))\otimes \xi_\tau )\otimes \lambda_{\sigma,\tau}
\ar[d]_{\textstyle 1\otimes \boldsymbol{c}\otimes 1\otimes 1} & (\lambda''_{\sigma,\tau}\otimes \xi'_{\sigma\tau})\otimes \xi_{\sigma\tau}\ar[u]_{\textstyle \cong}\\
((\xi'_\sigma\otimes(\xi'_\tau\otimes \xi_\sigma))\otimes \xi_\tau)\otimes \lambda_{\sigma,\tau}\ar[d]_{\textstyle \cong}&((\xi'_\sigma\otimes \xi'_\tau)\otimes \lambda'_{\sigma,\tau})\otimes
\xi_{\sigma\tau}\ar[u]_{\textstyle \xi'_{\sigma,\tau}\otimes 1}\\ (\xi'_\sigma\otimes \xi'_\tau)\otimes(( \xi_\sigma\otimes \xi_\tau)\otimes \lambda_{\sigma,\tau}) \ar[r]^-{\textstyle 1\otimes
\xi_{\sigma,\tau}}&  (\xi'_\sigma\otimes \xi'_\tau)\otimes (\lambda'_{\sigma,\tau}\otimes
\xi_{\sigma\tau})\ar[u]_{\textstyle \cong},}$$
and the horizontal composition of $2$-cells $f'$ and $f$ is $f'\!\otimes f$ where
 $(f'\!\otimes f)_\sigma= f'_\sigma\otimes f_\sigma$, for each arrow $\sigma$ in $I$.

 \underline{The identity} 1-cell of a $3$-cocycle is $1:\lambda\to\lambda$, where $1_\sigma=\mathrm{I}$ for all $\sigma$ in $I$, and each morphism $1_{\sigma,\tau}$ is determined by the commutativity of the square
   $$\xymatrix{(\mathrm{I}\otimes \mathrm{I})\otimes \lambda_{\sigma,\tau}\ar[r]^{\textstyle 1_{\sigma,\tau}}\ar[d]_{\textstyle \boldsymbol{r}\otimes 1}&\lambda_{\sigma,\tau}\otimes  \mathrm{I}   \ar[d]_{\textstyle \boldsymbol{r}} \\ \mathrm{I} \otimes \lambda_{\sigma,\tau}\ar[r]^{\textstyle \boldsymbol{l}}&\lambda_{\sigma,\tau}.}
   $$

\underline{The associativity and identity constraints}  in
$\ctco(I,(\m,\otimes,\boldsymbol{c}))$ are directly obtained from associativity and identity constraints of the braided monoidal category.

\vspace{0.2cm}

The bicategory $\ctco(I,(\m,\otimes,\boldsymbol{c}))$ is pointed by
the {\em trivial} $3$-cocycle, denoted by $\mathrm{I}_0$, which is
defined by the equalities $(\mathrm{I}_0)_{\sigma,\tau}=\mathrm{I}$
and $(\mathrm{I}_0)_{\sigma,\tau,\gamma}=1:\mathrm{I}\otimes
\mathrm{I}\to \mathrm{I}\otimes \mathrm{I}$.

\vspace{0.2cm}
We should note that, for any given category $I$ and braided monoidal
category $(\m,\otimes,\boldsymbol{c})$, there is a tricategory ${\bf
Lax}(I,\Omega^{^{-2}}\hspace{-3pt}\m)$ whose objects are lax
functors, whose 1-cells are lax transformations, whose  2-cells are
lax modifications and whose 3-cells are perturbations. Similarly as
the category of 2-cocycles $\cdco(I,(\m,\otimes))$ is a
subbicategory of ${\bf Lax}(I,\Omega^{^{-1}}\hspace{-3pt}\m)$, our
bicategory $\ctco(I,(\m,\otimes,\boldsymbol{c}))$ introduced above
is precisely the subtricategory of  ${\bf
Lax}(I,\Omega^{^{-2}}\hspace{-3pt}\m)$ given by the normal lax
functors and those lax transformations,  lax modifications and
perturbations  whose components at any objects are identities.

Both constructions ${Z}^3(I,(\m,\otimes,\boldsymbol{c}))$ and
$\ctco(I,(\m,\otimes,\boldsymbol{c}))$ are functorial on $I$, and
they lead to the following definition of geometric nerves for
braided monoidal categories:
\begin{definition} \label{defgnb} The  {\em geometric nerve} of a braided monoidal category
$(\m,\otimes,\boldsymbol{c})$ is  the simplicial set
\begin{equation}\label{gnbmc}Z^3(\m,\otimes,\boldsymbol{c}):\Delta^{\!^{\mathrm{op}}}\
\to \ \set,\hspace{0.4cm}[p]\mapsto
Z^3([p],(\m,\otimes,\boldsymbol{c})).
\end{equation}
This is the simplicial set of objects of the simplicial bicategory
\begin{equation}\label{cgnbmc}\ctco(\m,\otimes,\boldsymbol{c}):\Delta^{\!^{\mathrm{op}}}\
\to \mathbf{Hom}\subset \bicat,\hspace{0.4cm}[p]\mapsto
\ctco([p],(\m,\otimes,\boldsymbol{c})),
\end{equation}
which is called the {\em bicategorical geometric nerve} of the braided monoidal category.
\end{definition}

 \begin{remark}{\em  The geometric nerve $Z^3(\m,\otimes, \boldsymbol{c})$ is a $4$-coskeletal 1-reduced (one vertex, one 1-simplex) simplicial set whose $2$-simplices are the objects $\lambda_{0,1,2}$ of $\m$ and, for $p\geq 3$, the $p$-simplices are families of morphisms
$$ \lambda_{i,j,k,l}:\lambda_{j,k,l}\otimes \lambda_{i,j,l}\to \lambda_{i,j,k}\otimes \lambda_{i,k,l}
,$$
$0\leq i<j<k<l\leq p$, making commutative, for $0\leq i<j<k<l<m\leq p$, the diagrams
$$\xymatrix@C=90pt@R=36pt{(\lambda_{k,l,m}\otimes \lambda_{j,k,m})\otimes \lambda_{i,j,m}\ar[d]_{\textstyle \lambda_{j,k,l,m}\otimes 1}\ar[r]^{\textstyle \boldsymbol{a}(1\otimes \lambda_{i,j,k,m})\boldsymbol{a}^{-1}}&
(\lambda_{k,l,m}\otimes \lambda_{i,j,k})\otimes \lambda_{i,k,m}\ar[d]^{\textstyle \boldsymbol{c}\otimes 1}\\
(\lambda_{j,k,l}\otimes \lambda_{j,l,m})\otimes \lambda_{i,j,m}\ar[d]_{\textstyle (1\otimes \lambda_{i,j,l,m})\boldsymbol{a}}& (\lambda_{i,j,k}\otimes \lambda_{k,l,m})\otimes \lambda_{i,j,m}\ar[d]^{\textstyle (1\otimes \lambda_{i,k,l,m})\boldsymbol{a}}\\
\lambda_{j,k,l}\otimes (\lambda_{i,j,l}\otimes \lambda_{i,l,m})\ar[r]^{\textstyle \boldsymbol{a}^{-1}(\lambda_{i,j,k,l}\otimes 1)\boldsymbol{a}}& \lambda_{i,j,k}\otimes (\lambda_{i,k,l}\otimes \lambda_{i,l,m}).
 }
$$
}
\end{remark}

\vspace{0.2cm}
If $*$ is any object of a bicategory $\c$, then $\c(*,*)$ becomes a monoidal category and there is a
bicategorical embedding $\Omega^{^{-1}}\!\c(*,*)\hookrightarrow \c$. Since, for any braided monoidal
category $(\m,\otimes,\boldsymbol{c})$ and category $I$, there is a quite an obvious monoidal isomorphism
$$(\cdco(I,(\m,\otimes)),\otimes)\cong \ctco(I,(\m,\otimes,\boldsymbol{c}))(\mathrm{I}_0,\mathrm{I}_0),$$ we have a natural (`suspension') homomorphism of bicategories

$$S:\Omega^{^{-1}}\!\cdco(I,(\m,\otimes))\hookrightarrow \ctco(I,(\m,\otimes,\boldsymbol{c})),$$
 that is defined as the composite
 $$\Omega^{^{-1}}\!\cdco(I,(\m,\otimes))\cong \Omega^{^{-1}}\!\ctco(I,(\m,\otimes,\boldsymbol{c}))
 (\mathrm{I}_0,\mathrm{I}_0)\hookrightarrow \ctco(I,(\m,\otimes,\boldsymbol{c})).$$

Hence,  we have a simplicial homomorphism of simplicial bicategories
$$S:\Omega^{^{-1}}\!\cdco(\m,\otimes)\to
\ctco(\m,\otimes,\boldsymbol{c}),$$ whose composition with
(\ref{ipsh}) defines the pseudo-simplicial homomorphism
\begin{equation}\label{J}
E:\ner(\m,\otimes,\boldsymbol{c})\to
\ctco(\m,\otimes,\boldsymbol{c}),
\end{equation}
which, at each label $p\geq 0$, is so given by the commutative square
$$
\xymatrix{\Omega^{^{-1}}\hspace{-3pt}\m^p\ar[r]^-{\textstyle E_p}\ar[d]_{\textstyle \Omega^{^{-1}}
\hspace{-3pt}(~)^\mathrm{e}}&\ctco([p],(\m,\otimes,\boldsymbol{c}))\\
\Omega^{^{-1}}\hspace{-2pt}\cdco([p],(\m,\otimes))\ar[r]^-{\textstyle
\cong}&\Omega^{^{-1}}\ctco([p],(\m,\otimes,\boldsymbol{c}))(\mathrm{I}_0,\mathrm{I}_0).
\ar@{^{(}->}[u]
}
$$

Next Theorem \ref{int1} below states that this pseudo-simplicial
homomorphism (\ref{J}) induces a homotopy equivalence on classifying
spaces so that the simplicial bicategory
$\ctco(\m,\otimes,\boldsymbol{c})$, the bicategorical geometric
nerve,  models the homotopy type of the braided monoidal category
and it can be thought of as a `rectification' of the
pseudo-simplicial  nerve $\ner(\m,\otimes,\boldsymbol{c})$.

\begin{theorem}\label{int1} For any braided monoidal category $(\m,\otimes,\boldsymbol{c})$, the
pseudo-simplicial homomorphism $E:\ner(\m,\otimes,\boldsymbol{c})\to
\ctco(\m,\otimes,\boldsymbol{c})$ induces a homotopy equivalence on
classifying spaces. Thus,
$$\class(\m,\otimes,\boldsymbol{c})\simeq \class \ctco(\m,\otimes,\boldsymbol{c}).$$
\end{theorem}

\begin{proof}
In view of Theorem \ref{p.1.23}, it is sufficient to prove that
every homomorphism of bicategories
$E_n:\Omega^{^{-1}}\hspace{-3pt}\m^n\to
\ctco([n],(\m,\otimes,\boldsymbol{c}))$ induces a homotopy
equivalence on classifying spaces. The result is clear for $n=0$,
since $E_0$ is merely the obvious isomorphism between the two unit
(i.e., with only one $2$-cell) bicategories. For $n=1$, since the
trivial 3-cocycle $\mathrm{I}_0$ is the unique object of the
bicategory $\ctco([1],(\m,\otimes,\boldsymbol{c}))$, it is easy to
see that $E_1$ is actually an isomorphism of bicategories with an
inverse isomorphism
\begin{equation}\label{p1}
P_1: \ctco([1],(\m,\otimes,\boldsymbol{c}))
\overset{\cong}\longrightarrow \Omega^{^{-1}}\hspace{-3pt}\m,
\end{equation}
defined by $$
 P_1:\ \xymatrix@C=0.2pc{\mathrm{I}_0&\ {}_{\textstyle \Downarrow\!f}    &\ar@/_1pc/[ll]_{\textstyle \xi} \ar@/^0.8pc/[ll]^{\textstyle \xi'} \mathrm{I}_0}\hspace{0.2cm}\mapsto\hspace{0.2cm}
 \xymatrix@C=0pc{\ast&\ {}_{\textstyle \downarrow\!f_{0,1}}    &\ar@/_1pc/[ll]_{\textstyle \xi_{0,1}} \ar@/^0.8pc/[ll]^{\textstyle {\xi'}_{\hspace{-3pt}0,1}} \ast.}
 $$

Now, for $n\geq 2$,  our discussion uses the so-called {\em Segal
projections} (see \cite[Definition 1.2]{segal74}) that, on our
simplicial bicategory $\ctco(\m,\otimes,\boldsymbol{c})$, give the
homomorphisms
$$
P_n:\ctco([n],(\m,\otimes,\boldsymbol{c}))\longrightarrow
\Omega^{^{-1}}\hspace{-5pt}\m^n
$$
defined by the commutative triangles

\begin{equation}\label{tp}
\xymatrix{\ctco([n],(\m,\otimes,\boldsymbol{c}))\ar[rr]^{\textstyle
P_n}\ar[rd]\ar@{}@<-4pt>[rd]_<{\textstyle \prod\limits_{k=1}^n
d_0\cdots d_{k-2}
d_{k+1}\cdots d_n}&&\Omega^{^{-1}}\hspace{-5pt}\m^n.\\
&\ctco([1],(\m,\otimes,\boldsymbol{c}))^n\ar[ru]^{\cong}_{\textstyle
P_1^n}& }
\end{equation}

\noindent That is,
$$ P_n:\ \xymatrix@C=0.3pc{\lambda&\
{}_{\textstyle \Downarrow\!f}    &\ar@/_1pc/[ll]_{\textstyle
\,\xi\,} \ar@/^1pc/[ll]^{\textstyle\, \xi'\,} \lambda'}\mapsto
\xymatrix@C=1pc{\ast&
\downarrow\!(f_{0,1},\dots,f_{n-1,n})&\ar@/_1.5pc/[ll]_{\textstyle(\xi_{0,1},
\dots,\xi_{n-1,n})} \ar@/^1.2pc/[ll]^{\textstyle
(\xi'_{0,1},\dots,\xi'_{\hspace{-1pt}n-1,n})} \ast.} $$

\vspace{0.2cm}
For any $n\geq 2$, we have the equality $P_nE_n=1$ and, moreover, there is a oplax transformation,
$$ \Psi:1\Rightarrow E_nP_n:\ctco([n],(\m,\otimes,\boldsymbol{c}))\to \ctco([n],(\m,\otimes,\boldsymbol{c})),$$
whose component at a 3-cocycle $\lambda: [n]\to (\m,\otimes,\boldsymbol{c})$ is the 3-cocycle morphism
$\Psi\lambda=\psi:\lambda\to \mathrm{I}_0$,
where the objects $\psi_{i,j}$ of $\m$, for $i<j$, are inductively determined by the equalities
$$
\psi_{i,j+1}=\left\{\begin{array}{ll}\mathrm{I}& \text{ if } i=j,\\[4pt] \psi_{i,j}\otimes \lambda_{i,j,j+1}&\text{ if } i<j,\end{array}\right.
$$
and the morphisms $\psi_{i,j,k}:(\psi_{i,j}\otimes\psi_{j,k})\otimes \lambda_{i,j,k}\to \mathrm{I}\otimes\psi_{i,k}$, for $i<j<k$, are also inductively defined as follows:
each morphism $\psi_{i,j,j+1}$ is the canonical isomorphism making commutative the triangle
$$
\xymatrix@C=0pt@R=12pt{(\psi_{i,j}\otimes\mathrm{I})\otimes \lambda_{i,j,j+1}\ar[rr]^{\textstyle \psi_{i,j,j+1}}\ar[rd]_{\textstyle \boldsymbol{r}\otimes 1}^\cong &&\mathrm{I}\otimes(\psi_{i,j}\otimes \lambda_{i,j,j+1})\ar[dl]^{\textstyle \boldsymbol{l}}_\cong\\ &\psi_{i,j}\otimes \lambda_{i,j,j+1}&}
$$
and each morphism $\psi_{i,j,k+1}$ is obtained from the morphism $\psi_{i,j,k}$ as the composite dotted arrow
$$\xymatrix{(\psi_{i,j}\otimes(\psi_{j,k}\otimes \lambda_{j,k,k+1}))\otimes \lambda_{i,j,k+1}\ar[d]_{\textstyle \cong}\ar@{.>}[r]^-{\textstyle \psi_{i,j,k+1}}&\mathrm{I}\otimes (\psi_{i,k}\otimes \lambda_{i,k,k+1})\\
(\psi_{i,j}\otimes \psi_{j,k})\otimes (\lambda_{j,k,k+1}\otimes \lambda_{i,j,k+1})\ar[d]_{\textstyle 1\otimes\lambda_{i,j,k,k+1}} & (\mathrm{I}\otimes \psi_{i,k})\otimes \lambda_{i,k,k+1}\ar[u]_{\textstyle \cong}\\
(\psi_{i,j}\otimes \psi_{j,k})\otimes (\lambda_{i,j,k}\otimes
\lambda_{i,k,k+1}) \ar[r]^-{\textstyle \cong}&((\psi_{i,j}\otimes
\psi_{j,k})\otimes \lambda_{i,j,k})\otimes \lambda_{i,k,k+1}.
\ar[u]_{\textstyle \psi_{i,j,k}} }$$

 The component of $\Psi$ at a 3-cocycle morphism $\xi:\lambda\to\lambda'$ is the 2-cell in the bicategory $\ctco([n],(\m,\otimes,\boldsymbol{c}))$

 $$
 \xymatrix@C=30pt{ \lambda\ar[d]_{\textstyle \psi}\ar[r]^{\textstyle \xi}
 \ar@{}[rd]|{\textstyle \overset{\textstyle \widehat{\Psi}}\Rightarrow}&
 \lambda'\ar[d]^{\textstyle \psi'}\\ \mathrm{I}_0\ar[r]_{\textstyle \Omega^{^{-1}}\!X^\mathrm{e}}&\mathrm{I}_0}
 $$

\noindent where $\psi=\Psi\lambda$, $\psi'=\Psi\lambda'$, $X=P_n\xi=(\xi_{0,1},\dots,\xi_{n-1,n})$, and $X^\mathrm{e}$ are given as in (\ref{ef}), defined by the morphisms $\widehat{\Psi}_{i,j}:X^\mathrm{e}_{i,j}\otimes \psi_{i,j}\to \psi'_{i,j}\otimes \xi_{i,j}$  inductively obtained as follows: each morphism $\widehat{\Psi}_{i,i+1}$ is the canonical isomorphism making commutative the triangle
$$\xymatrix@R=10pt{\xi_{i,i+1}\otimes \mathrm{I}\ar[rd]_{\textstyle \boldsymbol{r}}^\cong\ar[rr]^{\textstyle \widehat{\Psi}_{i,i+1}}&&\mathrm{I}\otimes \xi_{i,i+1}\ar[ld]^{\textstyle \boldsymbol{l}}_\cong\\ &\xi_{i,i+1}&}$$
and each morphism $\widehat{\Psi}_{i,j+1}$ is obtained from the morphism $\widehat{\Psi}_{i,j}$ as the composite dotted arrow in the diagram below.
$$\xymatrix@C=60pt{ (X^\mathrm{e}_{i,j}\otimes \xi_{j,j+1})\otimes (\psi_{i,j}\otimes \lambda_{i,j,j+1})\ar[d]_{\textstyle \cong}\ar@{.>}[r]^-{\textstyle \widehat{\Psi}_{i,j+1}}& (\psi'_{i,j}\otimes \lambda'_{i,j,j+1})\otimes \xi_{i,j+1}\\
(X^\mathrm{e}_{i,j}\otimes (\xi_{j,j+1}\otimes \psi_{i,j}))\otimes \lambda_{i,j,j+1}\ar[d]_{\textstyle (1\otimes \boldsymbol{c})\otimes 1}&\psi'_{i,j}\otimes(\lambda'_{i,j,j+1}\otimes \xi_{i,j+1})\ar[u]_{\textstyle \cong} \\
(X^\mathrm{e}_{i,j}\otimes (\psi_{i,j}\otimes \xi_{j,j+1}))\otimes \lambda_{i,j,j+1}\ar[d]_{\textstyle \cong} &\psi'_{i,j}\otimes ((\xi_{i,j}\otimes \xi_{j,j+1})\otimes \lambda_{i,j,j+1})\ar[u]_{\textstyle 1\otimes \xi_{i,j,j+1}} \\
((X^\mathrm{e}_{i,j}\otimes \psi_{i,j})\otimes \xi_{j,j+1})\otimes \lambda_{i,j,j+1}\ar[r]^-{\textstyle (\widehat{\Psi}_{i,j}\otimes 1)\otimes 1}&((\psi'_{i,j}\otimes \xi_{i,j})\otimes \xi_{j,j+1})\otimes \lambda_{i,j,j+1}\ar[u]_{\textstyle \cong} }$$

Hence, by Fact \ref{f2} $(ii)$, every induced map $$\class
E_n:\class (\m,\otimes)^n\to \class
\ctco([n],(\m,\otimes,\boldsymbol{c})),$$  is a homotopy
equivalence (with $\class P_n: \class
\ctco([n],(\m,\otimes,\boldsymbol{c}))\to \class (\m,\otimes)^n$ as a
homotopy-inverse)    and therefore the induced map $\class
E:\class(\m,\otimes,\boldsymbol{c})\to
\class\ctco(\m,\otimes,\boldsymbol{c})$ is also a homotopy
equivalence by Theorem \ref{p.1.23}. \end{proof}

As we show below, Theorem \ref{int1} implies a new proof of a relevant  fact:  {\em The classifying space of the underlying category of a braided monoidal category is, up to group completion, a  double-loop space} \cite{sta,fie,b-f-s-v,berger}. Recall that the loop space of the classifying space of a monoidal category $\Omega\class(\m,\otimes)$ is a group completion of $\class \m$, the classifying space of the underlying category; that is, there is a homotopy  natural map (\ref{hgc}), $\class \m\to \Omega\class(\m,\otimes)$,  which is, up to group completion, a homotopy equivalence.

\begin{theorem} For any braided monoidal category $(\m,\otimes,\boldsymbol{c})$ there is a natural homotopy equivalence $$\class(\m,\otimes)\simeq\Omega\class(\m,\otimes,\boldsymbol{c}).$$
Therefore, the double-loop space $\Omega^2\class(\m,\otimes,\boldsymbol{c})$ is homotopy equivalent to the group completion of $\class\m$.
\end{theorem}
\begin{proof}
By Theorem \ref{int1}, $\class(\m,\otimes,\boldsymbol{c})$ is
homotopy equivalent to $\class\ctco(\m,\otimes,\boldsymbol{c})$,
the classifying space of the simplicial bicategory
$[n]\mapsto
\ctco([n],(\m,\otimes,\boldsymbol{c}))$, which, by the homotopy
equivalences (\ref{cg}), is itself homotopy equivalent to the
realization $|X|$ of the simplicial space $X:[n]\mapsto
\class\ctco([n],(\m,\otimes,\boldsymbol{c}))$.

Now, observe that: 1) the space $X_0$ is a one-point set; 2) the
Segal projection maps $p_n=\prod\limits_{k=1}^n d_0\cdots d_{k-2}
d_{k+1}\cdots d_n:X_n\to (X_1)^n$ are all homotopy equivalences
(since every map $\class P_n: \class
\ctco([n],(\m,\otimes,\boldsymbol{c}))\to \class (\m,\otimes)^n$ is
a homotopy equivalence, as we observed in the proof of Theorem
\ref{int1} above, and the triangles (\ref{tp}) commute); 3)
$X_1\cong \class(\m,\otimes)$ (by the isomorphism (\ref{p1})); and
4) $\pi_0(X_1)=0$, the trivial group (since, by Fact \ref{fm}, the
classifying space of the underlying monoidal category,
$\class(\m,\otimes)$, is homotopy equivalent to the geometric
realization of the  simplicial set with only one vertex
$Z^2(\m,\otimes)$).

Thus, we see that the simplicial space $X:[n]\mapsto
\class\ctco([n],(\m,\otimes,\boldsymbol{c}))$ satisfies the
hypothesis of Segal's Proposition 1.5 in \cite{segal74} (see also
the previous {\em Note} to the proposition). Therefore, the
canonical map $X_1\to \Omega|X|$ is a homotopy equivalence, whence
the homotopy equivalence
$\class(\m,\otimes)\simeq\Omega\class(\m,\otimes,\boldsymbol{c})$
follows.
\end{proof}

Going finally towards our last main result in the paper, let us recall from Definition \ref{defgnb}  that the geometric nerve of a braided monoidal category $Z^3(\m,\otimes,\boldsymbol{c})$ is the simplicial set of objects of the simplicial bicategory  $\ctco(\m,\otimes,\boldsymbol{c})$, so that we have an evident simplicial homomorphism of inclusion $Z^3(\m,\otimes,\boldsymbol{c})\hookrightarrow \ctco(\m,\otimes,\boldsymbol{c})$, where $Z^3(\m,\otimes,\boldsymbol{c})$ is regarded as a simplicial discrete bicategory.

\begin{theorem} For any braided monoidal category $(\m,\otimes,\boldsymbol{c})$, there is a  homotopy equivalence
$$\class(\m,\otimes,\boldsymbol{c})\simeq |Z^3(\m,\otimes,\boldsymbol{c})|.$$
\end{theorem}
\begin{proof} By taking into account Theorem \ref{int1}, it is sufficient to prove that the inclusion simplicial homomorphism $Z^3(\m,\otimes,\boldsymbol{c})\hookrightarrow \ctco(\m,\otimes,\boldsymbol{c})$ induces a homotopy equivalence on classifying spaces.
 To do so, let
$$\Delta^{\hspace{-2pt}^\mathrm{u}}\ctco(\m,\otimes,\boldsymbol{c}):
\Delta^{\!^{\mathrm{op}}}\!\times\!\Delta^{\!^{\mathrm{op}}}\to\set$$
$$([p],[q])\mapsto
\Delta^{\hspace{-2pt}^\mathrm{u}}_p\ctco([q],(\m,\otimes,\boldsymbol{c}))$$
be the bisimplicial set obtained from the simplicial bicategory
$$\ctco(\m,\otimes,\boldsymbol{c}):\Delta^{\!^{\mathrm{op}}}\to
\mathbf{Hom}\subset\bicat$$
 by composing with the unitary geometric nerve  functor (\ref{ngn}).

 Since a $(p,q)$-simplex of $\Delta^{\hspace{-2pt}^\mathrm{u}}\ctco(\m,\otimes,\boldsymbol{c})$ is then a normal lax functor
$$\xi:[p]\to \ctco([q],(\m,\otimes,\boldsymbol{c})),$$
 that consists of 3-cocycles $\xi^u:[q]\to (\m,\otimes,\boldsymbol{c})$, $0\leq u\leq p$, 1-cells
  $$\xi^{u,v}:\xi^v\to\xi^u,$$ $0\leq u\leq v\leq p$ and 2-cells $$\xi^{u,v,w}:\xi^{u,v}\otimes
   \xi^{v,w}\Rightarrow\xi^{u,w},$$ $0\leq u\leq v\leq w\leq p$, in the bicategory $\ctco([q],(\m,\otimes,\boldsymbol{c}))$, satisfying the various conditions, we see that $\xi$ can be described as a list of data
\begin{equation}\label{uxi}\xi=\Big(\xi^u_{i,j,k},\xi^u_{i,j,k,l},\xi^{u,v}_{i,j},\xi^{u,v}_{i,j,k},
\xi^{u,v,w}_{i,j}\Big)_{\hspace{-4pt}^{\scriptsize\begin{array}{l}_{0\leq u\leq v\leq w\leq p}\\ _{0\leq i\leq j\leq k\leq l\leq q}\end{array}}}\end{equation}
where  $$\xi^u_{i,j,k,l}:\xi^u_{j,k,l}\otimes \xi^u_{i,j,l}\to\xi^u_{i,j,k}\otimes \xi^u_{i,k,l}$$ are the morphisms in $\m$ that describe the $3$-cocycles $\xi^u$,
 $$\xi^{u,v}_{i,j,k}: (\xi^{u,v}_{i,j}\otimes\xi^{u,v}_{i,k})\otimes \xi^{v}_{i,j,k}\to \xi^{u}_{i,j,k}\otimes\xi^{u,v}_{i,k}$$ are the morphisms in $\m$ describing the 1-cells $\xi^{u,v}$, and $$\xi^{u,v,w}_{i,j}:\xi^{u,v}_{i,j}\otimes \xi^{v,w}_{i,j}\to \xi^{u,w}_{i,j}$$ are those morphisms in $\m$ that describe the 2-cells $\xi^{u,v,w}$.

Below, we shall interpret the $p$-(resp.\,$q$-)direction as the
horizontal (resp. vertical) one, so that the horizontal face and
degeneracy operators in
$\Delta^{\hspace{-2pt}^\mathrm{u}}\ctco(\m,\otimes,\boldsymbol{c})$
are those of the simplicial sets
$\Delta^{\hspace{-2pt}^\mathrm{u}}\!\ctco([q],(\m,\otimes,\boldsymbol{c}))$,
that is, $d_m^h\xi=(\xi^{d^mu}_{i,j,k},\dots)$, etc., whereas the
vertical ones are induced by those of
$\ctco(\m,\otimes,\boldsymbol{c})$, that is,
$d_m^v\xi=(\xi^{u}_{d^mi,d^mj,d^mk},\dots)$, etc.

Since $Z^3(\m,\otimes,\boldsymbol{c})$
is a simplicial discrete bicategory (i.e., all 1-cells and 2-cells
are identities),
$\Delta^{\hspace{-2pt}^\mathrm{u}}Z^3(\m,\otimes,\boldsymbol{c})$
is a bisimplicial set that is constant in the horizontal direction.
The induced bisimplicial inclusion
$\Delta^{\hspace{-2pt}^\mathrm{u}}Z^3(\m,\otimes,\boldsymbol{c})\hookrightarrow
\Delta^{\hspace{-2pt}^\mathrm{u}}\ctco(\m,\otimes,\boldsymbol{c})$
is then, at each horizontal level  $p\geq 0$, the composite
simplicial map
\begin{equation}\label{1.2.28'}Z^3(\m,\otimes,\boldsymbol{c})
\!
=\!\Delta^{\hspace{-2pt}^\mathrm{u}}_0\ctco(\m,\otimes,\boldsymbol{c})\overset{s_0^h}\hookrightarrow
\Delta^{\hspace{-2pt}^\mathrm{u}}_1\ctco(\m,\otimes,\boldsymbol{c})
\cdots\! \overset{s_{p-1}^h}\hookrightarrow
\Delta^{\hspace{-2pt}^\mathrm{u}}_p\ctco(\m,\otimes,\boldsymbol{c}).\end{equation}

Taking into account now (\ref{ef1}) and  that the classifying space of any diagram of bicategories $\f$ is homotopy equivalent to $|\hoco_I\Delta\f|$, to prove that
$$|Z^3(\m,\otimes,\boldsymbol{c})|
\hookrightarrow \class\ctco(\m,\otimes,\boldsymbol{c})$$ is a
homotopy equivalence, we shall prove that the induced simplicial map
on diagonals $ Z^3(\m,\otimes,\boldsymbol{c})\to \diag\,
\Delta^{\hspace{-2pt}^\mathrm{u}}\ctco(\m,\otimes,\boldsymbol{c})$
is a  weak equivalence. To do so, as every pointwise weak homotopy
equivalence bisimplicial map is a diagonal  weak homotopy
equivalence \cite[IV, Proposition 1.7]{g-j}, it suffices to  prove
that every  one of these simplicial maps (\ref{1.2.28'}) is a weak
homotopy equivalence. In fact, we will prove more: {\em Every
simplicial map
$$
s_{p-1}^h:\Delta^{\hspace{-2pt}^\mathrm{u}}_{p-1}\ctco(\m,\otimes,\boldsymbol{c})\hookrightarrow \Delta^{\hspace{-2pt}^\mathrm{u}}_p\ctco(\m,\otimes,\boldsymbol{c})
$$
embeds the simplicial set $\Delta^{\hspace{-2pt}^\mathrm{u}}_{p-1}\ctco(\m,\otimes,\boldsymbol{c})$ into $\Delta^{\hspace{-2pt}^\mathrm{u}}_p\ctco(\m,\otimes,\boldsymbol{c})$ as a simplicial deformation retract.}

To do so, since $d_p^hs_{p-1}^h=1$, it is enough to exhibit a simplicial homotopy
$$
H:1\Rightarrow s_{p-1}^hd_p^h:\Delta^{\hspace{-2pt}^\mathrm{u}}_p\ctco(\m,\otimes,\boldsymbol{c})\to\Delta^{\hspace{-2pt}^\mathrm{u}}_p\ctco(\m,\otimes,\boldsymbol{c}),
$$
which, for each $p\geq 1$, is given by the maps $h_m$, $0\leq m\leq q$, as in the diagram
$$ \xymatrix@R=20pt{\cdots & \Delta^{\hspace{-2pt}^\mathrm{u}}_p\ctco([q+1],(\m,\otimes,\boldsymbol{c})) \ar@<1.5ex>[rr]^{d_0^v} \ar@<-1.5ex>[rr]_{d_{q+1}^v}
\ar@<1ex>[dd]^{s_{p-1}^hd_p^h} \ar@<-1ex>[dd]_{1} & {}^{^{\vdots}} & \Delta^{\hspace{-2pt}^\mathrm{u}}_p\ctco([q],(\m,\otimes,\boldsymbol{c}))\ar@<1ex>[dd]^{s_{p-1}^hd_p^h}
\ar@<-1ex>[dd]_{1}
\ar@<1ex>[ddll]^{h_q}_{\ \ \dots}\ar@<-1ex>[ddll]_{h_0} & \cdots \\
&&&&\\
\cdots & \Delta^{\hspace{-2pt}^\mathrm{u}}_p\ctco([q+1],(\m,\otimes,\boldsymbol{c}))\ar@<1.5ex>[rr]^{d_0^v} \ar@<-1.5ex>[rr]_{d_{q+1}^v} &{}^{^{\vdots}}&\Delta^{\hspace{-2pt}^\mathrm{u}}_p\ctco([q],(\m,\otimes,\boldsymbol{c})) &
\cdots\\
}
$$
which take a $(p,q)$-simplex (\ref{uxi}) of $\Delta^{\hspace{-2pt}^\mathrm{u}}\ctco(\m,\otimes,\boldsymbol{c})$ to the $(p,q+1)$-simplex

$$h_m\xi=\Big((h_m\xi)^u_{i,j,k},(h_m\xi)^u_{i,j,k,l},(h_m\xi)^{u,v}_{i,j},(h_m\xi)^{u,v}_{i,j,k},
(h_m\xi)^{u,v,w}_{i,j}\Big)_{\hspace{-6pt}^{\scriptsize\begin{array}{l}_{0\leq u\leq v\leq w\leq p}\\_{ 0\leq i\leq j\leq k\leq l\leq q+1}\end{array}}}$$

\noindent defined as follows:

\vspace{0.2cm}
\noindent $\bullet$ The objects $(h_m\xi)^u_{i,j,k}$ are given by the formula

$$
(h_m\xi)^u_{i,j,k}=\left\{\begin{array}{ll}\xi^u_{s^mi,s^mj,s^mk}& \text{ if } \ u<p \text{ or } m<j,\\[6pt] \xi^{p-1}_{i,j,k}& \text{ if } \ u=p \text{ and } k\leq m,\\[6pt] \xi^{p-1,p}_{i,j}\otimes \xi^p_{i,j,k-1}&\text{ if } \ u=p \text{ and } j\leq m<k\,.\end{array}\right.
$$
$\bullet$ The morphisms $$(h_m\xi)^u_{i,j,k,l}:(h_m\xi)^u_{j,k,l}\otimes (h_m\xi)^u_{i,j,l}\to (h_m\xi)^u_{i,j,k}\otimes (h_m\xi)^u_{i,k,l}$$ are
$$
(h_m\xi)^u_{i,j,k,l}=\left\{\begin{array}{ll}\xi^u_{s^mi,s^mj,s^mk,s^ml}& \text{ if } \ u<p \text{ or } m<j,\\[6pt] \xi^{p-1}_{i,j,k,l}& \text{ if } \ u=p \text{ and } l\leq m,\end{array}\right.
$$
while, for $u=p$ and  $j\leq m<k$, the corresponding  $(h_m\xi)^p_{i,j,k,l}$ is defined as the composite dotted morphism
$$\xymatrix@C=55pt{\xi^p_{j,k-1,l-1}\otimes (\xi^{p-1,p}_{i,j}\otimes \xi^{p}_{i,j,l-1})\ar[d]_{\textstyle \cong}\ar@{.>}[r]^-{\textstyle (h_m\xi)^p_{i,j,k,l}}&(\xi^{p-1,p}_{i,j}\otimes \xi^p_{i,j,k-1})\otimes \xi^p_{i,k-1,l-1}\\
(\xi^p_{j,k-1,l-1}\otimes \xi^{p-1,p}_{i,j})\otimes \xi^{p}_{i,j,l-1}
\ar[d]_{\textstyle \boldsymbol{c}\otimes 1}&\xi^{p-1,p}_{i,j}\otimes (\xi^p_{i,j,k-1}\otimes \xi^p_{i,k-1,l-1})\ar[u]_{\textstyle \cong}\\
(\xi^{p-1,p}_{i,j}\otimes \xi^{p}_{j,k-1,l-1})\otimes \xi^{p}_{i,j,l-1}\ar[r]^{\textstyle \cong}&\xi^{p-1,p}_{i,j}\otimes ( \xi^{p}_{j,k-1,l-1}\otimes \xi^{p}_{i,j,l-1})\ar[u]_{\textstyle 1\otimes \xi^p_{i,j,k-1,l-1}} }$$
and for $k\leq m <l$ as the composite dotted morphism
$$\xymatrix@C=40pt{(\xi^{p-1,p}_{j,k}\otimes \xi^{p}_{j,k,l-1})\otimes (\xi^{p-1,p}_{i,j}\otimes \xi^{p}_{i,j,l-1})\ar[d]_{\textstyle \cong}\ar@{.>}[r]^-{\textstyle (h_m\xi)^p_{i,j,k,l}} & \xi^{p-1}_{i,j,k}\otimes(\xi^{p-1,p}_{i,k}\otimes \xi^{p}_{i,k,l-1})\\
((\xi^{p-1,p}_{j,k}\otimes \xi^{p}_{j,k,l-1})\otimes \xi^{p-1,p}_{i,j})\otimes \xi^{p}_{i,j,l-1}\ar[d]_{\textstyle \boldsymbol{c}\otimes 1} &(\xi^{p-1}_{i,j,k}\otimes\xi^{p-1,p}_{i,k})\otimes \xi^{p}_{i,k,l-1}\ar[u]_{\textstyle \cong} \\
(\xi^{p-1,p}_{i,j}\otimes (\xi^{p-1,p}_{j,k}\otimes \xi^{p}_{j,k,l-1}))\otimes \xi^{p}_{i,j,l-1}\ar[d]_{\textstyle \cong}& ((\xi^{p-1,p}_{i,j}\otimes \xi^{p-1,p}_{j,k})\otimes \xi^{p}_{i,j,k})\otimes \xi^{p}_{i,k,l-1}\ar[u]_{\textstyle  \xi^{p-1,p}_{i,j,k}\otimes 1} \\
(\xi^{p-1,p}_{i,j}\otimes \xi^{p-1,p}_{j,k})\otimes (\xi^{p}_{j,k,l-1}\otimes \xi^{p}_{i,j,l-1})\ar[r]^-{\textstyle  1\otimes \xi^p_{i,j,k,l-1}} & (\xi^{p-1,p}_{i,j}\otimes \xi^{p-1,p}_{j,k})\otimes (\xi^{p}_{i,j,k}\otimes \xi^{p}_{i,k,l-1}).\ar[u]_{\textstyle \cong} }$$

\vspace{0.2cm}\noindent $\bullet$  The objects $(h_m\xi)^{u,v}_{i,j}$ are defined by

$$
(h_m\xi)^{u,v}_{i,j}=\left\{\begin{array}{ll}\xi^{u,v}_{s^mi,s^mj}& \text{ if } \ v<p \text{ or } m<j,\\[6pt] \xi^{u,p-1}_{i,j}& \text{ if } \ v=p \text{ and } j\leq m\,.
\end{array}\right.
$$
\noindent$\bullet$ The morphisms $$(h_m\xi)^{u,v}_{i,j,k}:((h_m\xi)^{u,v}_{i,j}\otimes (h_m\xi)^{u,v}_{j,k})\otimes (h_m\xi)^{v}_{i,j,k}\to (h_m\xi)^u_{i,j,k}\otimes (h_m\xi)^{u,v}_{i,k}$$ are
$$
(h_m\xi)^{u,v}_{i,j,k}=\left\{\begin{array}{ll}\xi^{u,v}_{s^mi,s^mj,s^mk}& \text{ if } \ v<p \text{ or } m<j,\\[6pt] \xi^{u,p-1}_{i,j,k}& \text{ if } \ v=p \text{ and } k\leq m,\end{array}\right.
$$
and if $v=p$ and $j\leq m<k$, then the morphism $(h_m\xi)^{u,p}_{i,j,k}$ is the composite
$$\xymatrix@C=15pt{(\xi^{u,p-1}_{i,j}\otimes \xi^{u,p}_{j,k-1})\otimes (\xi^{p-1,p}_{i,j}\otimes\xi^{p}_{i,j,k-1})\ar[d]_{\textstyle \cong}\ar@{.>}[r]^-{\textstyle (h_m\xi)^{u,p}_{i,j,k}}&\xi^u_{i,j,k-1}\otimes \xi^{u,p}_{i,k-1}\\
(\xi^{u,p-1}_{i,j}\otimes (\xi^{u,p}_{j,k-1}\otimes \xi^{p-1,p}_{i,j}))\otimes\xi^{p}_{i,j,k-1}\ar[d]_{\textstyle (1\otimes\boldsymbol{c})\otimes 1}& (\xi^{u,p}_{i,j}\otimes \xi^{u,p}_{j,k-1})\otimes\xi^{p}_{i,j,k-1}\ar[u]_{\textstyle \ \xi^{u,p}_{i,j,k-1}}\\
(\xi^{u,p-1}_{i,j}\otimes(\xi^{p-1,p}_{i,j}\otimes \xi^{u,p}_{j,k-1}))\otimes\xi^{p}_{i,j,k-1}\ar[r]^{\textstyle \cong}&((\xi^{u,p-1}_{i,j}\otimes\xi^{p-1,p}_{i,j})\otimes \xi^{u,p}_{j,k-1})\otimes\xi^{p}_{i,j,k-1}.\ar[u]_{\textstyle (\xi^{u,p-1,p}_{i,j}\otimes 1)\otimes 1} }$$
$\bullet$ The morphisms $$(h_m\xi)^{u,v,w}_{i,j}:(h_m\xi)^{u,v}_{i,j}\otimes (h_m\xi)^{v,w}_{i,j}\to (h_m\xi)^{u,w}_{i,j}$$ are given by
 $$
 (h_m\xi)^{u,v,w}_{i,j}=\left\{\begin{array}{ll}\xi^{u,v,w}_{s^mi,s^mj}& \text{ if } \ w<p \text{ or } m<j,\\[6pt] \xi^{u,v,p-1}_{i,j}& \text{ if } \ w=p \text{ and } j\leq m.\end{array}\right.
 $$

So defined, a straightforward (though quite tedious) verification shows that ${H:1\Rightarrow s_{p-1}^hd_p^h}$ is actually a simplicial homotopy, and this completes the proof.
\end{proof}

\section{Appendix: Coherence conditions}

\noindent {\bf (CC1):}  for $m \overset{d}\to\ell \overset{c}\to k\overset{b}\to j \overset{a}\to i$,
any four composible arrows of $I$, the following equation on modifications holds
$$
\xymatrix@R=-1pt@C=-1pt{&&d^*c^*b^*a^*\ar@2{->}[drr]^{\textstyle d^*c^*\chi} \ar@2{->}[dll]_{\textstyle \chi b^*
a^*}\ar@2{->}[ddd]|-{\textstyle d^*\chi a^*}&& &&&& d^*c^*b^*a^*\ar@2{->}[drr]^{\textstyle d^*c^*\chi}
\ar@2{->}[dll]_{\textstyle \chi b^*a^*}&& \\
(cd)^*b^*a^*\ar@2{->}[ddd]_{\textstyle \chi a^*}& & & & d^*c^*(ab)^* \ar@2{->}[ddd]^{\textstyle d^*\chi} &&(cd)^*b^*a^*\ar@2{->}[ddd]_{\textstyle \chi
a^*}\ar@2{->}[drr]_{(cd)^*\!\chi} &&\overset{(\ref{4})}\cong &&d^*c^*(ab)^* \ar@2{->}[ddd]^{\textstyle d^*\chi}\ar@2{->}[dll]^-{\chi(ab)^*}
\\&\underset{\Rrightarrow}{\omega a^*}&&\underset{\Rrightarrow}{d^*\omega}~~&&&&&
(cd)^*\! (ab)^*\ar@2{->}[ddd]_{\textstyle \chi}&&\\
&&d^*(bc)^* a^*\ar@2{->}[lld]_{\textstyle \chi a^*}
\ar@2{->}[drr]^{\textstyle d^*\chi}&&&=&&\underset{\Rrightarrow}{\omega}&&\underset{\Rrightarrow}{\omega}&\\
(bcd)^*a^*\ar@2{->}[drr]_{\textstyle \chi}&&\underset{\Rrightarrow}{\omega}&&
d^*(abc)^*\ar@2{->}[dll]^{\textstyle \chi}&&(bcd)^*a^*\ar@2{->}[drr]_{\textstyle \chi}&&&&
d^*(abc)^*\ar@2{->}[dll]^{\textstyle \chi} \\
&&(abcd)^*&&&&&&(abcd)^*&&}
$$
\noindent {\bf (CC2):}   for any two composible arrows $k\overset{b}\to j\overset{a}\to i$  of $I$,
$$\label{nor}\xymatrix@R=8pt@C=15pt{&&b^*a^*\ar@2{->}[dddll]_{\textstyle 1_{b^*}a^*}^{\textstyle \overset{\Rrightarrow}{\delta a^*}}
\ar@2{->}[dddrr]^{\textstyle b^*1_{a^*}}_{\textstyle \overset{\Rrightarrow}{b^*\gamma}}
\ar@2{->}[ddd]|-{b^*\iota a^*} &&&&& b^*a^* \ar@2{->}[dddl]_{\textstyle 1_{b^*}a^*}\ar@{}[dd]|>>{\textstyle \cong}\ar@2{->}[dddr]^{\textstyle b^*1_{a^*}} &\\
&&& ~&&&& &\\
&&&&&&&&\\
b^*a^*\ar@2{->}[ddrr]_{\textstyle \chi} && b^*1_j^*a^* \ar@2{->}[ll]^{\textstyle \chi a^*} \ar@2{->}[rr]_{\textstyle b^*\chi} && b^*a^*
\ar@2{->}[ddll]^{\textstyle \chi} & = & b^*a^* \ar@2{->}[ddr]_{\textstyle \chi} \ar@2{->}[rr]^{\textstyle 1_{b^*a^*}}  && b^*a^*\ar@2{->}[ddl]^{\textstyle \chi} \\
&& \underset{\Rrightarrow}{\omega} &&&&& \overset{\textstyle\cong}{~} & \\
&& (ab)^* &&&&& (ab)^* & \\ }
$$

\noindent {\bf (CC3):} for any three composible arrows $\ell \overset{c}\to k\overset{b}\to j \overset{a}\to i$of $I$,
$$
\xymatrix@R=3pt@C=-1pt{&&F_\ell c^*b^*a^*\ar@2{->}[ddd]|-{F_\ell\chi a^*}\ar@2{->}[dddrr] ^{\textstyle F_\ell
c^*\chi}\ar@2{->}[dddll]_{\textstyle \theta b^*a^*}&&&&&&F_\ell c^*b^*a^*\ar@2{->}[dddrr]
^{\textstyle F_\ell c^*\chi}\ar@2{->}[dddll]_{\textstyle \theta b^*a^*}&&\\
&&&&&&&&&&\\
&&&&&&&&\overset{(\ref{4})}\cong&&\\
c^*F_kb^*a^*\ar@2{->}[dd]|-{\textstyle c^*\theta a^*}&&F_\ell(bc)^*a^*\ar@2{->}[dd] |-{ \theta a^*}\ar@2{->}[ddrr]^{F_\ell\chi}
&\underset{\Rrightarrow}{F_\ell\omega}&F_\ell c^*(ab)^*\ar@2{->}[dd]|-{\textstyle F_\ell\chi}&&c^*F_k
b^*a^*\ar@2{->}[rr]^{c^*F_k\chi}\ar@2{->}[dd]|-{\textstyle c^*\theta a^*}&& c^*F_k(ab)^*\ar@2{->}[dd]|-{c^*\theta}&&
F_\ell c^*(ab)^*\ar@2{->}[dd]|-{\textstyle F_\ell\chi}\ar@2{->}[ll]_(0.4){\theta(ab)^*}\\
&\underset{\Rrightarrow}{\Pi a^*}&&&&&&\underset{\Rrightarrow}{c^*\Pi}& &\underset{\Rrightarrow}{\Pi}
&\\
c^*b^*F_ja^*\ar@2{->}[rr]_{\chi' F_ja^*}\ar@2{->}[ddd]|-{\textstyle c^*b^*\theta}&&(bc)^*F_ja^*\ar@2{->}[ddd] |-{(bc)^*\theta}&
&F_\ell(abc)^*\ar@2{->}[ddd]^-{\textstyle \theta}&=&c^*b^*F_ja^* \ar@2{->}[ddd]|-{\textstyle c^*b^*\theta}&&c^*(ab)^*F_i\ar@2{->}[dddrr]^{\chi'
F_i}&&
F_\ell(abc)^*\ar@2{->}[ddd]^-{\textstyle \theta}\\
&\overset{(\ref{4})}\cong&&\underset{\Rrightarrow}\Pi&&&&&\\
&&&&&&&&\ar@3{->}[u]{\omega'F_i}&&\\
c^*b^*a^*F_i\ar@2{->}[rr]_{\textstyle \chi'\!a^*\!F_i}&& (bc)^*a^*\!F_i\ar@2{->}[rr]_{\textstyle \chi'\!
F_i}&&(abc)^*F_i&&c^*b^*a^*F_i\ar@2{->}[rr]_{\textstyle \chi'\!a^*\!F_i} \ar@2{->}[uuurr]|-{c^*\chi'\! F_i}&&
(bc)^*a^*F_i\ar@2{->}[rr]_{\textstyle \chi'\! F_i}&&(abc)^*F_i }
$$

\noindent {\bf (CC4):} for  $a:j\to i$ any arrow in  $I$,
$$
\xymatrix@R=4pt@C=4pt{
F_j1_j^*a^*\ar@{=>}[rr]^{\textstyle \theta a^*}\ar@{=>}[ddddrr]^{F_j\chi}&&1_j^*F_ja^*\ar@{=>}[rr]^{\textstyle 1_j^*\theta}&&1_j^*a^*F_i\ar@{=>}[dddd]^(.3){\textstyle \chi'\!F_i}&&&F_j1_j^*a^*\ar@{=>}[rr]^{\textstyle \theta a^*}&&1_j^*F_ja^*\ar@{=>}[rr]^{\textstyle 1_j^*\theta}&&1_j^*a^*F_i\ar@{=>}[dddd]^(.3){\textstyle \chi'\!F_i}\\
&&&\ar@3{->}[d]_{\textstyle \Pi}&&&&\ar@3{->}@<3ex>[d]&\hspace{-7pt}\Gamma a^*\,\,&\overset{(\ref{4})}\cong &&\\
\ar@3{->}@<3ex>[d]&&&&&=&&&&a^*F_i\ar@{=>}[uurr]|{\iota'\!a^*\!F_i}\ar@{=>}[ddrr]|{1_{a^*}\! F_i}&\,\,\,\gamma'F_i\hspace{-9pt}&\ar@3{->}@<-3ex>[d]\\
&\hspace{-9pt}F_j\gamma\,\,\,&&&&&&&&\overset{(\ref{4})}\cong&&\\
F_ja^*\ar@{=>}[uuuu]^(.3){\textstyle F_j\iota a^*}\ar@{=>}[rr]_{\textstyle F_j1_{a^*}}&&F_ja^*\ar@{=>}[rr]_{\textstyle \theta}&&a^*F_i&&&F_ja^*\ar@{=>}[uuuu]^(.3){\textstyle F_j\iota a^*}\ar@{=>}[rr]_{\textstyle F_j1_{a^*}}\ar@{=>}[uuuurr]|(.3){\iota'\! F_j\!a^*}\ar@{=>}[uurr]_{\theta}&&F_ja^*\ar@{=>}[rr]_{\textstyle \theta}
&&a^*F_i
}
$$
$$
\xymatrix@R=3pt@C=3pt{
F_ja^*1_i^*\ar@{=>}[rr]^{\textstyle \theta 1_i^*}\ar@{=>}[dddddrr]^{F_j\chi}&&a^*F_i1_i^*\ar@{=>}[rr]^{\textstyle a^*\theta}&&a^*1_i^*F_i\ar@{=>}[ddddd]^(.3){\textstyle \chi'\!F_i}&&&F_ja^*1_i^*\ar@{=>}[rr]^{\textstyle \theta 1_i^*}&\ar@{}[dddddr]_(.42){\textstyle \underset{(\ref{4})}\cong}&a^*F_i1_i^*\ar@{=>}[rr]^{\textstyle a^*\theta}\ar@3{->}@<4ex>[dd]^{\textstyle a^*\Gamma}&&a^*1_i^*F_i\ar@{=>}[ddddd]^(.3){\textstyle \chi'\!F_i}\\
&&&&&&&&&&&\\
&&\ar@3{->}[d]{\textstyle \Pi}&&&&&& &&&\\
\ar@3{->}@<3ex>[d]&&&&&=&&&&\ar@{=>}[uuu]^{a^*\!F_i\iota}a^*F_i\ar@{=>}[uuurr]_{a^*\!\iota'\!F_i}\ar@{=>}[ddrr]|{1_{a^*}\! F_i}&\,\,\,\delta'^{-1}F_i&\ar@3{->}@<-5ex>[d]\\
&\hspace{-7pt}F_j\delta^{-1}\,\,\,&&&&&&&&\overset{(\ref{4})}\cong&&\\
F_ja^*\ar@{=>}[uuuuu]^(.3){\textstyle F_ja^*\iota}\ar@{=>}[rr]_{\textstyle F_j1_{a^*}}&&F_ja^*\ar@{=>}[rr]_{\textstyle \theta}&&a^*F_i&&&F_ja^*\ar@{=>}[uuuuu]^(.3){\textstyle F_ja^*\iota}\ar@{=>}[rr]_{\textstyle F_j1_{a^*}}\ar@{=>}[uurr]_{\theta}&&F_ja^*\ar@{=>}[rr]_{\textstyle \theta}
&&a^*F_i
}
$$
\noindent {\bf (CC5):} for any two composible arrows of $I$, $k\overset{b}\to j \overset{a}\to i$,
$$\xymatrix@C=5pt@R=1pt{ &&F_kb^*a^*\ar@2{->}[drr]^{\textstyle mb^*\!a^*} \ar@2{->}[dll]_{\textstyle \theta a^*} &&&&&& F_kb^*a^*
\ar@2{->}[drr]^{\textstyle mb^*\!a^*}\ar@2{->}[dll]_{\textstyle \theta a^*}\ar@2{->}[dddddd]_{ F_k\chi}&&\\
b^*F_ja^*\ar@2{->}[ddd]|-{\textstyle b^*\theta}\ar@2{->}[drr]_{b^*ma^*}&& \underset{\Rrightarrow}{\mathrm{M}a^*}&
&F'_kb^*a^*\ar@2{->}[dll]^{\theta'a^*}\ar@2{->}[dddddd]|-{\textstyle F'_k\chi}&&b^*F_ja^* \ar@2{->}[ddd]|-{\textstyle b^*\theta}
&&&&F'_kb^*a^*\ar@2{->}[dddddd]|-{\textstyle F'_k\chi}\\
&&b^*F'_ja^*\ar@2{->}[ddd]|-{b^*\theta'}&&&&&&&&\\
&\underset{\Rrightarrow}{b^*\mathrm{M}}&&&&&&&&&\\
b^*a^*F_i\ar@2{->}[ddd]|-{\textstyle \chi'F_i}\ar@2{->}[drr]^{b^*a^*m}
&&&\underset{\Rrightarrow}{\Pi'}&&\hspace{-7pt}=\hspace{-7pt}&b^*a^*F_i\ar@2{->}[ddd]|-{\textstyle \chi'F_i}& \underset{\Rrightarrow}\Pi
&&\overset{(\ref{4})}\cong &\\
&&b^*a^*F'_i\ar@2{->}[ddd]|-{\chi'F'_i}&&&&&&&&\\
&\overset{(\ref{4})}\cong &&&&&&&F_k(ab)^*
\ar@2{->}[drr]^{m(ab)^*}\ar@2{->}[dll]_{\theta}&&\\
(ab)^*F_i\ar@2{->}[drr]_{\textstyle (ab)^*m}&&&&F'_k(ab)^*
\ar@2{->}[dll]^{\textstyle \theta'}&&(ab)^*F_i\ar@2{->}[drr]_{\textstyle (ab)^*m}
&&\underset{\Rrightarrow}{\mathrm{M}} &&F'_k(ab)^*\ar@2{->}[dll]^{\textstyle \theta'}\\
&&(ab)^*F'_i&&&&&&(ab)^*F'_i&& }$$
\noindent {\bf (CC6):}  for any object $i$ of $I$,
$$ \xymatrix@R=4pt@C=2pt{&& F_i\ar@{=>}[lldd]_{\textstyle F_i\iota} \ar@{=>}[rrrr]^{\textstyle m}&&&& F'_i\ar@{=>}[lldd]_{F'_i\iota}\ar@{=>}[dddd]^{\textstyle \iota' F'_i} &&&&&&&&F_i\ar@{=>}[rrr]^{\textstyle m}\ar@{=>}@<-2pt>[dddd]^(0.3){\iota'\! F_i}\ar@{=>}[lllldd]_{\textstyle F_i\iota} &&& F'_i\ar@{=>}[dddd]^{\textstyle \iota'F'_i}\\
&& &\overset{(\ref{4})}\cong &&&&&&&&&&&&&&\\
F_i1_i^* \ar@{=>}[rrrr]^{m1_i^*}\ar@{=>}[rrdd]_{\textstyle \theta} &&&& F'_i1_i^* \ar@{=>}[rrdd]_{\theta'}\ar@{}[rr]|(0.6){{\textstyle \overset{\textstyle \Gamma'}\Rrightarrow}}&&&&~~~=&& F_i1_i^*\ar@{=>}[rrrrdd]_{\textstyle \theta}\ar@{}[rrrr]|(0.6){{\textstyle \overset{\textstyle \Gamma}\Rrightarrow}} &&&& \ar@{}[rr]|(0.8){{\textstyle \overset{(\ref{4})}\cong}}&&&\\
&& &  \overset{\textstyle \mathrm{M}}\Rrightarrow &&&&& &&&&& &&&& \\
&& 1_i^*F_i\ar@{=>}[rrrr]_{\textstyle 1_i^*m} &&&& 1_i^*F'_i &&&&&&&& 1_i^*F_i \ar@{=>}[rrr]_{\textstyle 1_i^*m} && &1_i^*F'_i }$$

\noindent {\bf (CC7):} for any arrow $a: j\to i$ of the category $I$, the square below commutes.
$$
\xymatrix{ a^*m_i\circ \theta_a\ar@3{->}[r]^{\textstyle \mathrm{M}}
\ar@3{->}[d]_{\textstyle a^*\!\sigma_i\circ 1} &\theta'_a\circ m_ja^* \ar@3{->}[d]^{\textstyle 1\circ \sigma_{\!j}a^*} \\
a^*m'_i\circ \theta_a\ar@3{->}[r]^{\textstyle \mathrm{M}'}&\theta'_a\circ m'_ja^* }
$$

\end{document}